\pdfoutput=1

\documentclass{amsart}
\usepackage{amsmath,amsthm}
\usepackage{amsfonts,amssymb}

\usepackage{graphicx}
\usepackage{ifpdf}
\usepackage{bbm}

\newtheorem{thm}{Theorem}[section]

\newtheorem{lem}[thm]{Lemma}
\newtheorem{prop}[thm]{Proposition}

\newtheorem{defn}[thm]{Definition}

 \usepackage{ifpdf}
    \numberwithin{equation}{section}
    \numberwithin{table}{section}
    \numberwithin{figure}{section}

\theoremstyle{remark}

\makeatletter
\newcommand{\bddots}{%
  \mathinner{\mkern1mu\raise\p@\vbox{\kern7\p@\hbox{.}}\mkern2mu
    \raise4\p@\hbox{.}\mkern2mu\raise7\p@\hbox{.}\mkern1mu}}
\makeatother

\def\ONE{{\mathbbm 1}}

\def\a{{\alpha}}

\def\eb{{\mathbf e}}
\def\m{{\mathbf m}}

\def\ib{{\mathbf i}}
\def\jb{{\mathbf j}}
\def\kb{{\mathbf k}}
\def\lb{{\mathbf l}}
\def\mb{{\mathbf m}}
\def\vb{{\mathbf v}}
\def\sb{{\mathbf s}}
\def\tb{{\mathbf t}}

 \def\NN{{\mathbb N}}
 \def\CA{{\mathcal A}}
 
 \def\CG{{\mathcal G}}
 
 \def\CP{{\mathcal P}}
 \def\CK{{\mathcal K}}
 \def\CI{{\mathcal I}}
 \def\CL{{\mathcal L}}
 \def\CS{{\mathcal S}}
 \def\CT{{\mathcal T}}
 \def\CTC{{\mathcal T\!C}}
 \def\CTS{{\mathcal T\!S}}
 
 \def\CH{{\mathcal H}}
 \def\CP{{\mathcal P}}
 
 \def\CC{{\mathbb C}}
 
 \def\HH{{\mathbb H}}
 \def\PP{{\mathbb P}}
 
 \def\RR{{\mathbb R}}
 
 \def\ZZ{{\mathbb Z}}
 \def\VV{{\mathbb V}}
 
 \def\A{{\mathcal A}}
 \def\G{{\mathcal G}}
 
\newcommand{\e}{\mathrm{e}}

\newcommand{\tr}{{\mathsf {tr}}}
\newcommand{\TC}{{\mathsf {TC}}}
\newcommand{\TS}{{\mathsf {TS}}}

\def \la {\langle}
\def \ra {\rangle}
\newcommand{\wt}{\widetilde}
\newcommand{\wh}{\widehat}

 \ifpdf
 \DeclareGraphicsExtensions{.pdf,.png,.jpg}
 \else
 \DeclareGraphicsExtensions{.eps}
 \fi
 \graphicspath{{images/}}

\begin{document}

\title[Discrete Fourier analysis in $d$-variables]
{Discrete Fourier analysis on fundamental domain of  $A_d$ lattice 
and on simplex in $d$-variables}

\author{Huiyuan Li}
\address{Institute of Software\\
Chinese Academy of Sciences\\ Beijing 100080,China}
\email{hynli@mail.rdcps.ac.cn}
\author{Yuan Xu}
\address{  Department of Mathematics\\ University of Oregon\\
    Eugene, Oregon 97403-1222.}
\email{yuan@math.uoregon.edu}

\date{\today}

\thanks{The first author was supported by NSFC Grants 10601056,
10431050 and 60573023. The second author was supported by NSF
Grant DMS-0604056}

\begin{abstract}
A discrete Fourier analysis on the fundamental domain $\Omega_d$ of
the $d$-dimensional lattice of type $A_d$ is studied, where
$\Omega_2$ is the regular hexagon and $\Omega_3$ is the rhombic
dodecahedron, and analogous results on $d$-dimensional simplex are
derived by considering invariant and anti-invariant elements. Our
main results include Fourier analysis in trigonometric functions,
interpolation and cubature formulas on these domains. In particular,
a trigonometric Lagrange interpolation on the simplex is shown to
satisfy an explicit compact formula and the Lebesgue constant of the
interpolation is shown to be in the order of $(\log n)^d$. The basic
trigonometric functions on the simplex can be identified with
Chebyshev polynomials in several variables already appeared in
literature. We study common zeros of these polynomials and show that
they are nodes for a family of Gaussian cubature formulas, which
provides only the second known example of such formulas.
\end{abstract}

\maketitle

 \tableofcontents


\section{Introduction}
\setcounter{equation}{0}

The classical discrete Fourier analysis works with  the periodic
exponential functions $\e^{2\pi i k x}$ for $k$ in the set $\ZZ_n :
= \{0,1,\ldots, n-1\}$. The central piece is the discrete Fourier
transform defined via an inner product over the set of equally
spaced points $\{k/n: k \in \ZZ_n\}$ in $[0,1)$, which agrees with
the continuous Fourier transform for trigonometric polynomials of
degree at most $n$. This is equivalent to Gaussian quadrature
formula for trigonometric polynomials, and it can be used to define
trigonometric interpolation based on the equally spaced points. All
three quantities, discrete Fourier transform, quadrature, and
interpolation, are important tools in numerous applications and form
an integrated part of the discrete Fourier analysis.

In several variables, periodicity of functions can be defined via a
lattice $L$, which is a discrete group described by $A\ZZ^d$, where
$A$ is a nonsingular matrix. A function $f$ is periodic with respect
to $L=A\ZZ^d$ if $f(x) = f(x+A k)$ for any $k \in \ZZ^d$. The usual
multiple Fourier analysis of $\RR^d$ uses periodicity defined by the
lattice $\ZZ^d$. One can develop a discrete Fourier analysis in
several variables with a lattice periodicity, which has been defined
in connection with signal processing and sampling theory (see, for
example, \cite{DM, Hi, Ma}). The domain on which the analysis takes
place is the fundamental domain of the lattice, which is a bounded
set $\Omega$ that tiles $\RR^d$ in the sense that $\Omega + A \ZZ^d
= \RR^d$. It is known that the family of exponentials $\{\e^{2\pi i
\alpha \cdot x }: \alpha \in  L^\perp\}$, in which $L^\perp$ is the
dual lattice, forms an orthonormal basis for $L^2(\Omega)$ and these
exponentials are periodic with respect to $L$. These two facts make
it possible to establish a complete analogue of the classical
Fourier analysis. Recent in \cite{LSX}, a discrete Fourier transform
is defined and used to study cubature and trigonometric
interpolation on the domain $\Omega$, which provides a general
framework for the discrete Fourier analysis with lattice tiling.
Detail studies are carried out in \cite{LSX} for the hexagon
lattice, for which the fundamental domain is a regular hexagon, when
$d=2$, and in \cite{LX08} for the face-centered cubic lattice, for
which the fundamental domain is a rhombic dodecahedron, when $d =
3$. The regular hexagon and the rhombic dodecahedron are invariant
under the reflection group $\CA_2$ and $\CA_3$, respectively. In the
present paper we shall consider the lattice of type $A_d$ for all
$d \ge 2$, for which the fundamental domain is the union of the
images of a regular simplex in $\RR^d$ under the group action.

Our goal is then to establish a discrete Fourier analysis on the
fundamental domain of the $A_d$ lattice in $\RR^d$ that is
comparable with the classical theory in one variable as outlined in
the first paragraph of this introduction. However, this is only the
first step. The frame work developed in \cite{LSX} applies to the
fundamental domain $\Omega$ that tiles $\RR^d$ without overlapping,
which means that $\Omega$ contains only part of its own boundary.
For example, if the lattice is $\ZZ^d$, then the domain is the half
open and half closed $[-\frac 12, \frac12)^d$. The points on which
the discrete Fourier transform is defined consist of lattice points
in $n \Omega$, which is not symmetric under the reflection group
$\CA_d$, as part of the
boundary points are not included. Thus, as a second step, we will
develop a discrete Fourier analysis that uses a symmetric set of
points by working with the congruence of the boundary points under
translation by the lattice and under the action of the group
$\CA_d$, which requires delicate analysis. Furthermore, by
restricting to functions that are invariant under the group, the
results on the set of symmetric points in the step 2 can be
transformed to results on a simplex that makes up the fundamental
domain. This is our third step, which establishes a Fourier analysis
on trigonometric functions on a regular simplex in $\RR^d$. For the
classical Fourier analysis, this step amounts to a discrete analysis
on cosine and sine functions. We will define analogues of cosine and
sine functions on a simplex in the step 3. The classical Chebyshev
polynomials arise from cosine and sine functions. As our fourth
step, we define generalized Chebyshev polynomials of the first and
the second kind from those generalized cosine and sine functions,
respectively, and study the discrete analysis of these algebraic polynomials.

In each of the steps outlined above, we will develop a concrete discrete
Fourier analysis associated with the $A_d$ lattice in $\RR^d$, providing
explicit formulas and detailed analysis. The cases $d=2$ and $d=3$
studied in \cite{LSX} and \cite{LX08} provide a road map for our study.
However, the extension from the cases of two and three variables to
the general $d$-variables is far from a trivial exercise. In order to carry
out our program on the discrete analysis, it is essential that we have a
complete understanding of the boundary of the fundamental domain
and its congruent relations under the group. For the step 3, we also
need to understand the boundary of the simplex and its relation with
the fundamental domain. In $\RR^2$ and $\RR^3$, we can grasp the
geometry of the hexagon and the rhombic dodecahedron by looking at
their graphs and understand the relative positions between faces, edges
and vertices, which are all parts of the boundary of the domain. This is
not going to be the case in $\RR^d$ for $d > 3$. For $\RR^d$ we will
have to rely on the group theory and describe the boundary elements
by working with subgroups. The study in the cases of $d =2$ and $d=3$
has revealed much of what can be done and they will serve as examples
throughout the paper.

The generalized cosine and sine functions on the simplex, as well as
the generalized Chebyshev polynomials, have been studied before in
the literature. In fact, in the case of $d =2$, they are first studied by
Koornwinder in \cite{K}, who started with the eigenfunctions of the
Laplace operator and showed that the generalized Chebyshev
polynomials are orthogonal on the region bounded by the Steiner's
hypocycloid. Later these polynomials are generated and studied by
several authors (see \cite{Ba,Ba2,Be, De, DL, DL2, EL, R} and the
reference therein), often as eigenfunctions of certain differential operators.
In \cite{Be}, they appear as eigenfunctions of the radial part of the
Laplacian-Beltrami operator on certain symmetric spaces. We highly
recommend the article \cite{Be}, which also includes a detailed account
on the history of these polynomials. Furthermore, the generalized
Chebyshev polynomials belong to large families of special functions
(see, for example, \cite{BeOp}), just as the classical Chebyshev
polynomials are the special cases of Jacobi polynomials. This, however,
appears to have little bearing with our study in this paper. It should
be mentioned, however, that our aim is to develop a discrete Fourier
analysis on the generalized Chebyshev polynomials, which come
down to study common zeros of these polynomials, a topic that does
not seem to have been studied systematically before.

There are other types of discrete Fourier analysis in several variables.
For the multiple Fourier analysis, with lattice $\ZZ^d$, there is a rich
collection of literature as it is essentially the tensor product of one
dimensional results. Among others, we mention the recent work on
anti-symmetric exponential and trigonometric functions by Klimyk
and Patera, see \cite{KP1, KP2} and the references therein.

The paper is organized as follows. In Section 2 we shall recall the
basic facts on lattice and tiling, and describe the framework for
the discrete Fourier analysis established in \cite{LSX}. In Section
3 we study the case of $A_d$ lattice and carry out both the step 2
and the step 3 of our program. The best way to describe the $A_d$
lattice and its fundamental domain appears to be using homogeneous
coordinates, which amounts to regard $\RR^d$ as the hyperplane $t_1
+\ldots + t_{d+1} =0$ of $\RR^{d+1}$. In fact, the group $\CA_d$
becomes permutation group $S_{d+1}$ in homogeneous coordinates. In
Section 4 we establish a discrete Fourier analysis on the simplex,
which studies trigonometric polynomials on the standard simplex in
$\RR^d$. One of the concrete result is a Lagrange interpolation by
trigonometric polynomial on equal spaced points on the simplex,
which can be computed by a compact formula and has an operator norm,
called the Lebesgue constant, in the order of $(\log n)^d$. Finally,
in Section 5, we derive the basic properties of the generalized
Chebyshev polynomials of the first and the second kind, and study
their common zeros, which are  related intrinsically to the cubature
formula for algebraic polynomials. In particular, it is known that a
Gaussian cubature exists if and only if the set of its nodes is
exactly the variety of the ideal generated by the orthogonal
polynomials of degree $n$, which however happens rarely and only one
single family of examples is known in the literature. We will show
that this is the case for the generalized Chebyshev polynomials of
the second kind, so that a Gaussian cubature formula exists for the
weight function with respect to which the Chebyshev polynomials of
the second kind are orthogonal, which provides the second family of
examples for such a cubature.

\section{Discrete Fourier analysis with lattice}
\setcounter{equation}{0}

In this section we recall results on discrete Fourier analysis associated
with lattice. Background and the content of the first subsection can be
found in \cite{CS, DM, Hi, Ma}. Main results, stated in the second subsection,
are developed in \cite{LSX} but will be reformulated somewhat, for which
explanation will be given. We will be brief and refer the proof and further
discussions to \cite{LSX}.


\subsection{ Lattice and Fourier series}

We will consider lattices on $\RR^d$ but will use homogeneous coordinates,
which amounts to identify $\RR^d$ with the hyperplane $x_1+\ldots + x_{d+1}=0$
in $\RR^{d+1}$. For convenience, we will  state the results on the lattices
on a $d$-dimensional subspace $\VV^d$ of the Euclidean space $\RR^m$
with $m \ge d$. Since $\VV^d$ can be identified with $\RR^d$, all results
proved for $\RR^d$ can be restated for $\VV^d$.

Let $\VV^d$ be a $d$-dimensional subspace of $\RR^m$ with $m\ge d$.
A lattice $L$ of $\VV^d$ is a discrete subgroup that contains $d$ linearly
independent vectors,
\begin{align*}
   L:=\left\{ k_1 a_1 +k_2a_2+\cdots+k_da_d:\ k_i \in \ZZ,\ i=1,2,\cdots,d \right\},
 \end{align*}
where $a_1,\cdots,a_d$ are linearly independent column vectors in
$\VV^d$. Let $A $ be the matrix of size $m\times d$, which has
column vectors $a_1,\cdots,a_d$. Then $A$ is called a {\it generator
matrix} of the lattice $L$. We can write $L$ as $A\ZZ^d$; that is
 \begin{align*}
        L = A\ZZ^d := \left\{A k:\ k\in \ZZ^d \right\}.
 \end{align*}
The dual lattice $L^{\perp}$ of $L$ is given by
 \begin{align*}
 L^{\perp}:=\ &\left\{ x\in \VV^d: \ x \cdot y \in \ZZ \text{ for all } y \in L\right\},
 \end{align*}
where  $x\cdot y$ denotes the usual Euclidean inner product of $x$
and $y$. The generator matrix of $L^\perp$, denoted by $A^\perp$, is given
by $A^\perp= A (A^\tr A)^{-1}$, which is the transpose of the Moor-Penrose
inverse of $A$. In particular, if $m = d$, then the generator matrix of
$L^\perp$  is simply $A^{-\tr}$ and $L^{\perp}= A^{-\tr}\ZZ^d
=\left\{A^{-\tr}k:\ k \in \ZZ^{d} \right\}$.

A bounded set $\Omega\subset \VV^d$ is said to tile  $\VV^d$ with
the lattice $L$ if
\begin{align*}
 \sum_{\alpha\in L} \chi_{\Omega} (x+\alpha) = 1 \quad \text{ for almost all }
     x \in \VV^d,
\end{align*}
where $\chi_{\Omega}$ denotes the characteristic function of
$\Omega$, which we write as $\Omega+L= \VV^d$. Tiling and Fourier
analysis are closely related as demonstrated by the Fuglede theorem.
Let $\int_{\Omega} f(x) dx$ denote the integration of the function
$f$ over $\Omega$. Let $\langle \cdot, \cdot \rangle_{\Omega}$
denote the inner product in $L^2({\Omega})$,
\begin{align}
\label{eq:innerproduct} \langle f,g\rangle_{\Omega} :=
\frac{1}{|\Omega|} \int_{\Omega} f(x) \overline{g(x)} dx,
\end{align}
where $|\Omega|$ denotes the measure of $\Omega$ and the bar denotes
the complex conjugation. The following
fundamental result was proved by Fuglede in \cite{F}.

\begin{thm}  \label{Fuglede}
Let $\Omega\subset\VV^d$ be a bounded domain and $L$ be a lattice of
$\VV^d$. Then $\Omega+L= \VV^d$ if and only if $\left\{\e^{2\pi i
\kappa \cdot x}: \kappa \in L^{\perp} \right\}$ is an orthonormal
basis with respect to the inner product \eqref{eq:innerproduct}.
\end{thm}

Written explicitly, the orthonormal property states that
\begin{align}
\label{eq:orthonormal}
  \frac{1}{|\Omega|} \int_{\Omega} \e^{2\pi i\, \kappa \cdot x} dx =
      \delta_{\kappa,0}, \qquad \kappa\in L^{\perp}.
\end{align}
Moreover, if $L = A\ZZ^d$, then the measure
$|\Omega|=\sqrt{\det (A^{\tr}A)}$.

The set $\Omega$ is called a {\it spectral set}  for the lattice $L$. If
$L=A\ZZ^d$ we also write $\Omega_A$ instead of $\Omega$.

Because of Theorem \ref{Fuglede}, a function $f\in L^1(\Omega)$ can be
expanded into a Fourier series
\begin{align*}
 f(x) \sim \sum_{\kappa\in L^{\perp}} c_{\kappa} \e^{2\pi i\, \kappa \cdot x }, \qquad
 c_{\kappa} := \frac{1}{|\Omega|} \int_{\Omega} f(x) \e^{-2\pi i\, \kappa \cdot x} dx.
\end{align*}
The Fourier transform $\widehat{f}$ of a function defined on
$L^1({\VV^d})$ and its inversion are defined by
\begin{align*}
\widehat{f}(\xi): = \int_{\VV^d} f(x) \e^{-2\pi i\, \xi \cdot x }
dx, \qquad {f}(x) := \int_{\VV^d} \widehat{f}(\xi) \e^{2\pi i\, \xi
\cdot x } d\xi.
\end{align*}
Our first result is the following sampling theorem (see, for example, \cite{Hi,Ma}).

\begin{prop} \label{pro:interpolation}
Let $\Omega$
 be the spectral set of the lattice $L$. Assume that $\widehat{f}$ is
supported on $\Omega$ and $\widehat{f}\in L^2(\Omega)$. Then
\begin{align*}
   f(x) = \sum_{\kappa\in L^{\perp}} f(\kappa) \Phi_{\Omega} (x-\kappa)
\end{align*}
in $L^2(\Omega)$, where
\begin{align*}
\Phi_{\Omega} (x) = \frac{1}{|\Omega|} \int_{\Omega} \e^{2\pi i \xi
\cdot x} d\xi.
\end{align*}
\end{prop}

This theorem is a consequence of the Poisson summation formula. We
notice that
\begin{align*}
\Phi_{\Omega} (\kappa) = \delta_{0,\kappa} \qquad \text{ for all }
\kappa\in L^{\perp}
\end{align*}
by Theorem \ref{Fuglede}, so that $\Phi_{\Omega}$ can be considered
as a cardinal interpolation function.


\subsection{Discrete Fourier analysis and interpolation}
A function $f$ defined on $\VV^d$ is called {\it periodic} with
respect to the lattice $L$ if
\begin{align*}
          f(x+\alpha) = f(x)\qquad \text{ for all } \alpha\in L.
\end{align*}
The spectral set $\Omega$ of the lattice $L$ is not unique. In order
to carry out the discrete Fourier analysis with respect to the lattice,
we shall fix an $\Omega$ such that $\Omega$ contains $0$ in
its interior and we further require that $\Omega$ tiles $\VV^d$ with
$L$ without overlapping and without gap. In other words, we require
that
\begin{align}
\label{eq:tiling2}
 \sum_{\alpha\in L}\chi_{\Omega} (x+\alpha) = 1 \qquad \text{ for all }  x\in \VV^d.
\end{align}
For example,  for the standard cubic lattice $\ZZ^d$ we can take
$\Omega=[-\frac12,\frac12)^d$. Correspondingly a spectral set
satisfying \eqref{eq:tiling2} is called a {\it fundamental domain}
of Lattice $L$.

\begin{defn} \label{def:N}
Let $L_A=A\ZZ^d$ and $L_B=B\ZZ^d$ be two lattices of $\VV^{d}$,
$\Omega_A$ and $\Omega_B$ satisfy \eqref{eq:tiling2}. Assume $N:=
B^\tr A$ is a nonsingular matrix that has all entries being integers.
Define
\begin{align*}
 \Lambda_N    := \Omega_A \cap L_B^{\perp}\quad \hbox{and} \quad
 \Lambda_N^{\dag} := \Omega_B \cap L_A^{\perp}.
\end{align*}
\end{defn}

Two points $x,y\in \VV^d$ are said to be congruent with respect to
the lattice $A\ZZ^d$, if $x-y\in A\ZZ^d$, and we write $x\equiv y
\pmod{A}$ or $x\equiv y \pmod{A\ZZ^d}$. The following two theorems
are the central results for the discrete Fourier transform.

\begin{thm}  \label{th:DFT}
Let $A,\,B$ and $N$ be as in Definition \ref{def:N}. Then for any
$\kappa\in L_A^{\perp}$,
\begin{align}
\label{eq:DFT} \frac{1}{|\det(N)|} \sum_{\alpha\in \Lambda_N}
\e^{2\pi i\, \kappa \cdot \alpha}
 =\begin{cases} 1, & \text{if }\ \kappa \equiv 0 \pmod{B}, \\
  0, & \text{otherwise},
   \end{cases}
\end{align}
and
\begin{align}
\label{eq:IDFT} \frac{1}{|\det(N)|} \sum_{\kappa\in
\Lambda_N^{\dag}} \e^{-2\pi i\, \kappa \cdot \alpha}
 =\begin{cases} 1, & \text{if }\ \alpha \equiv 0 \pmod{A}, \\
  0, & \text{otherwise}.
   \end{cases}
\end{align}
\end{thm}

\begin{thm} \label{thm:2.4}
Let $A,\, B$ and $N$ be as in Definition \ref{def:N}. Define the
discrete inner product
\begin{align*}
  \langle f,g \rangle_N = \frac{1}{|\det(N)|} \sum_{\alpha\in \Lambda_N} f(\alpha)
    \overline{ g(\alpha)}
\end{align*}
for $f,\, g \in C(\Omega_A)$, the space of continuous functions on
$\Omega_A$. Then
\begin{align}
\label{eq:equiv_inner}
 \langle f,\,g \rangle_N = \langle f,\,g \rangle
\end{align}
for all $f,\,g$ in the finite dimensional subspace
\begin{align*}
  \mathcal{H}_N := \mathrm{span}\left\{\phi_{\kappa}:\ \phi_{\kappa}(x) = \e^{2\pi i\, \kappa \cdot x},\ \kappa \in \Lambda_N^\dag \right\}.
\end{align*}
\end{thm}

Let $|E|$ denote the cardinality of the set $E$. Setting $\kappa =0$ or
$\alpha=0$ in \eqref{eq:DFT} or \eqref{eq:IDFT}, respectively, we see
that
\begin{equation} \label{L=L^d}
     | \Lambda_N| = | \Lambda_{N^\dag}| = |\det(N)|.
\end{equation}
In particular, the dimension of $\mathcal{H}_N$ is $|\Lambda_N^\dag| = |\det (N)|$.

Let $\mathcal{I}_N f$ denote the Fourier expansion of $f\in
C(\Omega_A)$ in $\mathcal{H}_N$ with respect to the inner product
$\langle\cdot, \cdot\rangle_N$. Then, analogous to the sampling
theorem in Proposition \ref{pro:interpolation}, $\CI_Nf$ satisfies
the following formula
\begin{align*}
 \mathcal{I}_N f(x) =\sum_{\alpha\in \Lambda_N} f(\alpha) \Phi^A_{\Omega_B} (x-\alpha), \quad f\in C(\Omega_A),
\end{align*}
where
\begin{align*}
\Phi^A_{\Omega_B} (x) = \frac{1}{|\det(N)|} \sum_{\kappa\in
\Lambda_N^\dag}
      \e^{2\pi i\, \kappa \cdot x}.
\end{align*}
The following theorem shows that $\CI_N f$ is an interpolation function.

\begin{thm} \label{thm:interpolation}
Let $A$, $B$ and $N$ be as in Definition \ref{def:N}. Then $\mathcal{I}_N f$ is the
unique interpolation operator on $N$ in $\mathcal{H}_N$; that is
\begin{align*}
  \mathcal{I}_N f(\alpha) = f(\alpha), \quad \forall \alpha\in \Lambda_N.
\end{align*}
In particular, $|\Lambda_N | = |\Lambda_N^\dag|$. Furthermore, the
fundamental interpolation function $\Phi^A_{\Omega_B}$ satisfies
\begin{align}
\label{eq:LagInter}
\Phi^A_{\Omega_B}(x) = \sum_{\kappa\in L_A}
\Phi_{\Omega_B}(x+\kappa).
\end{align}
\end{thm}

\begin{proof} Equation \eqref{eq:LagInter} was proved in \cite{LSX}. Using
\eqref{eq:IDFT} gives immediately $\Psi^A_{\Omega_B}=\delta_{\alpha,0}$
for $\alpha\in \Lambda_N$, so that the interpolation holds. This also shows
that $\{\Psi_{\Omega_B}^A(x- \alpha): \alpha \in \Lambda_N\}$  is linearly
independent. By \eqref{L=L^d}, $| \Lambda_N| = |\det (N)| = |\Lambda_N^{\dag}|$.
Furthermore, \eqref{eq:DFT} and \eqref{eq:IDFT} show that the interpolation
matrix $M = ( \phi_{\kappa}(\alpha) )_{\kappa\in \Lambda^{\dag},\alpha
\in \Lambda_N}$ is invertible. Consequently the interpolation on points in
$\Lambda_N$  is unique.
\end{proof}

This theorem is stated in \cite{LSX} under the additional requirement that
$ \Lambda_{N}^{\dag}=\Lambda_{N^{\tr}}$, which holds in particular in the
case that $A$ is a constant multiple of $B$. We prove the more general
version here for future references.

The results state in this subsection provide a framework for the
discrete Fourier analysis on the spectral set of a lattice. We will
apply it to the $d$-dimensional lattice associated with the group
$\A_d$ in the following section. As mentioned in the introduction,
the case $d=2$ and $d=3$ have been considered in detail in
\cite{LSX} and \cite{LX08}, respectively. These lower dimensional
cases provide a roadmap for the $d$-dimensional results.


\section{Discrete Fourier analysis on the fundamental domain}
\setcounter{equation}{0}

\subsection{Lattice $A_d$ and Fourier analysis}
For $d\ge 1$, we identify $\RR^d$ with the hyperplane
\begin{align*}
 \RR_H^{d+1} := \left\{(t_1,t_2,\cdots,t_{d+1})\in \RR^{d+1}:
 t_1+t_2+\cdots+t_{d+1}=0\right\}
\end{align*}
of $\RR^{d+1}$. The lattice $A_d$ that we will consider in this
paper is simply constructed by
\begin{align*}
\ZZ^{d+1}_H := \ZZ^{d+1} \cap \RR_H^{d+1}=
 \left\{(k_1,k_2,\cdots,k_{d+1})\in \ZZ^{d+1}:k_1+k_2+\cdots+k_{d+1}=0\right\}.
\end{align*}
In other words, we will use {\it homogeneous coordinates} $\tb \in
\RR^{d+1}_H$ to describe our results in $d$-variables and
$A_d=\ZZ^{d+1}_H$ in our homogeneous coordinates.

Throughout this paper, we adopt the convention of using bold letters, such as
$\tb$ and $\kb$, to denote homogeneous coordinates. The advantage of
homogeneous coordinates lies in preservation of symmetry. In fact, many
of our formulas are symmetric and more transparent under homogeneous
coordinates \cite{LSX, LX08, Sun, Sun2}.

The lattice $A_d$ is the root lattice of the reflection group $\A_d$
\cite[Chapter 4]{CS}. Under homogeneous coordinates, the group
$\A_d$ is generated by the reflections $\{\sigma_{ij}: 1 \le i < j
\le d+1\}$, where $\sigma_{ij}$ is defined by
$$
 \tb \sigma_{ij} := \tb - 2 \frac{\la \tb, \eb_{i,j}\ra}{\la \eb_{i,j},\eb_{i,j}\ra} \eb_{i,j}
     = \tb - (t_i-t_j)\eb_{i,j},
\quad \hbox{where $\eb_{i,j} := e_i- e_j$}.
$$
The last equation shows that $\sigma_{i,j}$ is a transposition.
Thus, the group $\A_d$ is exactly the permutation group $S_{d+1}$ of
$d+1$ elements. Let $\mathcal{G} := S_{d+1}$. For $\tb \in
\RR^{d+1}_H$ and $\sigma \in \CG$, the action of $\sigma$ on $\tb$
is denoted by $\tb \sigma$, which means permutation of the element
in $\tb$ by $\sigma$.

\begin{figure}[ht]
\centering
\hfill\includegraphics[width=0.5\textwidth]{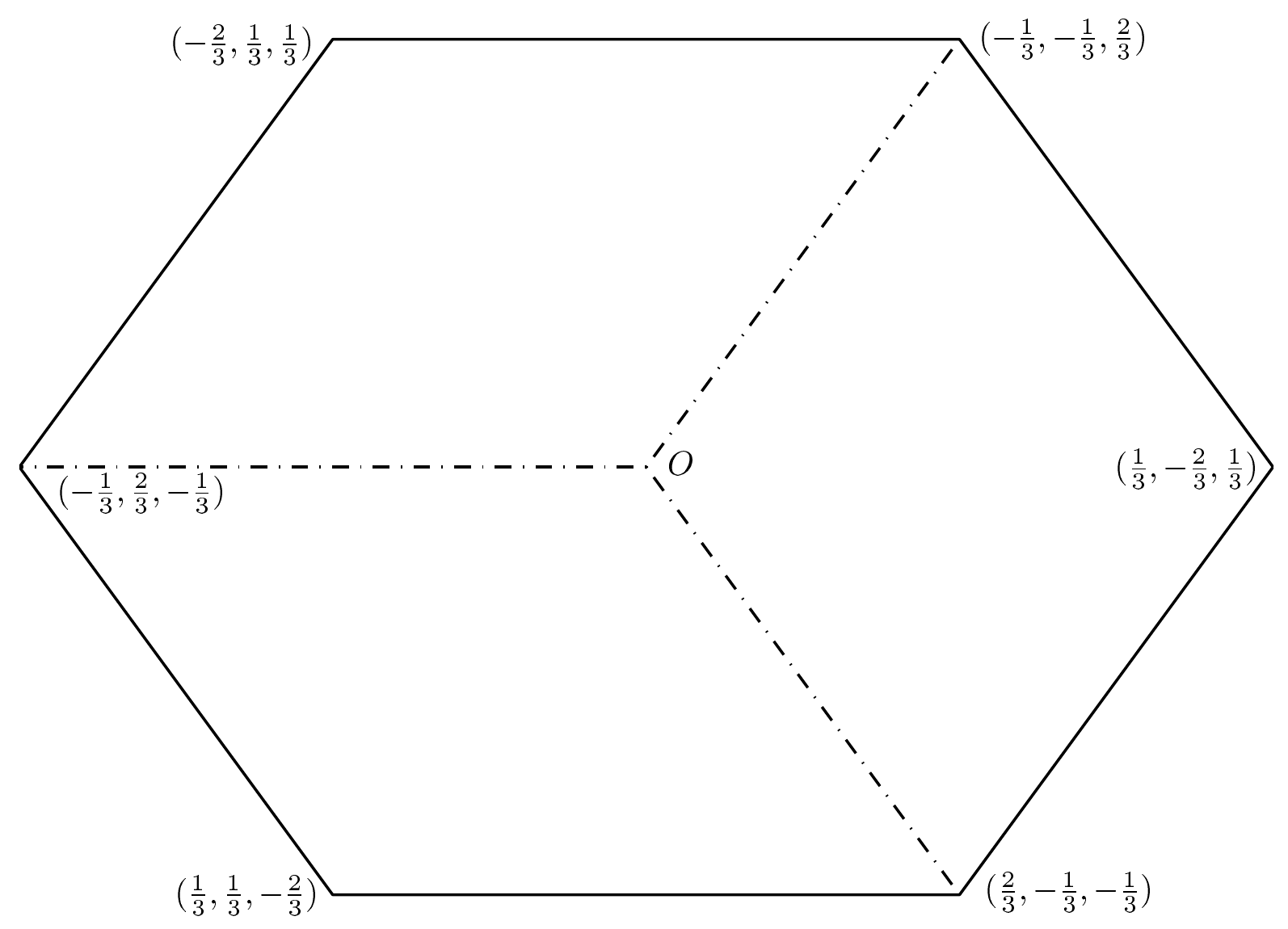}\hspace*{\fill}
\caption{Fundamental domain --- the regular hexagon for the lattice
$A_2$.}\label{hexagon}
\end{figure}

\begin{figure}[ht]
\centering
\hfill\includegraphics[width=0.7\textwidth]{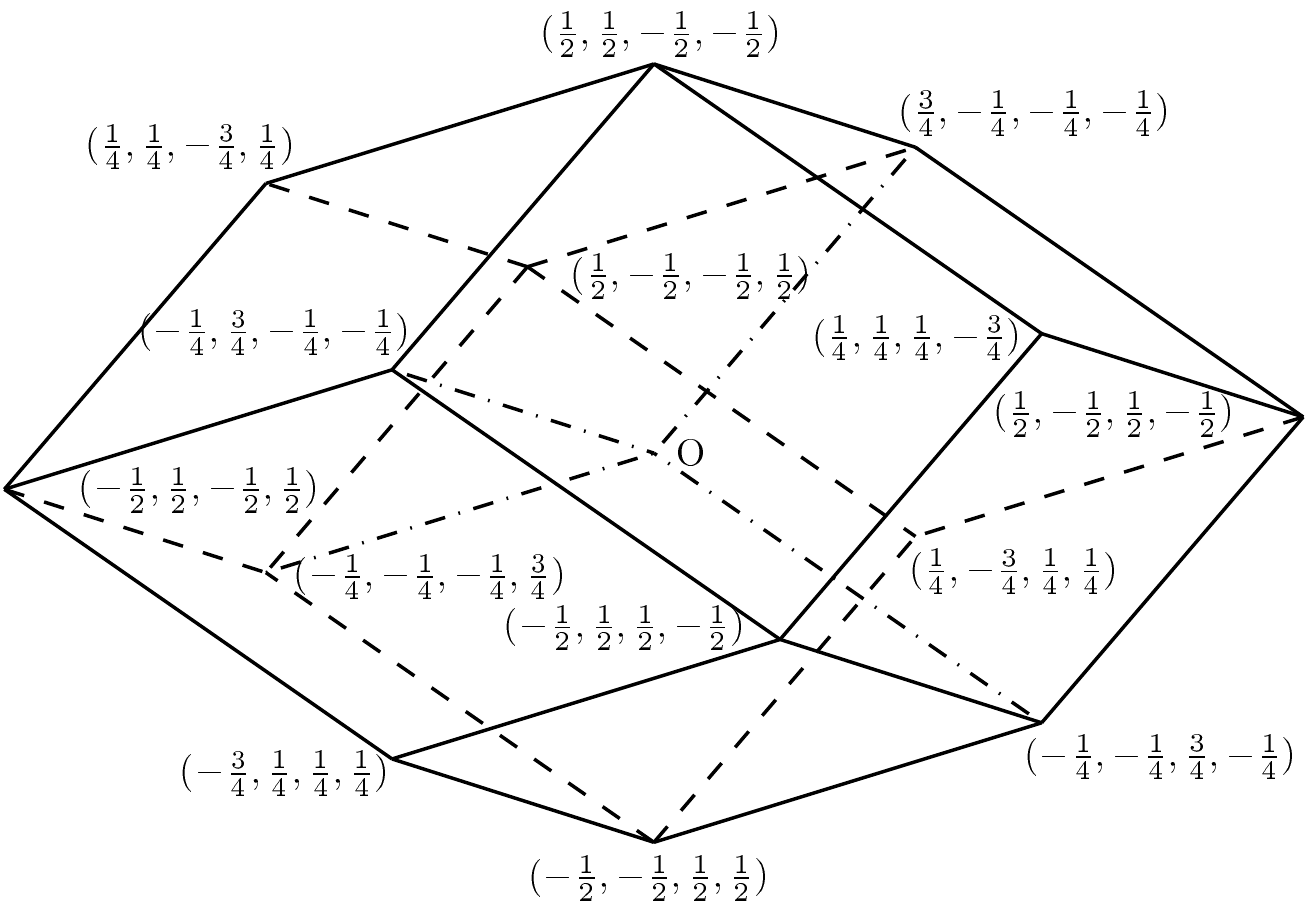}\hspace*{\fill}
\caption{Fundamental domain --- the rhombic dodecahedron for the
lattice $A_3$.}\label{dodecahedron}
\end{figure}

The fundamental domain that tiles $\RR^{d+1}_H$ with lattice $A_d$
is the spectral set of $A_d$,
\begin{equation}\label{Omega}
   \Omega_H := \left\{\tb \in \RR_H^{d+1}:
             -1< t_i-t_j  \leq  1,\,  1\leq i < j \leq d+1\right\}.
\end{equation}
The strict inequality in the definition of $\Omega_H$ reflects our requirement
in \eqref{eq:tiling2}. In the case of $d=2$ and $d=3$,
the fundamental domain are hexagon and rhombic dodecahedron, respectively,
and they are depicted in Figures \ref{hexagon}-\ref{dodecahedron}, in which
vertices are labeled in homogeneous coordinates.

As described in \cite[Chapter 21]{CS}, the spectral set
$\overline{\Omega}_H$ is the union of its {\it fundamental simplex} $\triangle_H$,
defined by
\begin{align}\label{triangle}
\triangle_H: =\left\{\tb\in \RR^{d+1}_H: 0\le t_i-t_j\le 1,\ 1\le i< j \le d+1\right\}.
\end{align}
Let $\underbrace{t,t,\dots,t}_k$ be denoted by $\left\{ t \right\}^k$. The $d$
non-zero vertices of $\triangle_H$ are given by
\begin{align*}
\frac{\vb^{k}}{d+1}, \quad 1\le k \le d, \qquad \hbox{where
       $\vb^k:=\big(\{d+1-k\}^k, \{-k\}^{d+1-k}\big)$}.
\end{align*}
Let $\tb\mathcal{G}$ denote the orbit of $\tb$ under $\mathcal{G}$,
that is, $\tb \G:= \left\{\tb\sigma :\ \sigma\in \mathcal{G}\right\}$.  Then
$\overline{\Omega}_H$ is the union of the images of $\triangle_H$ under
$\CG$, that is
\begin{align}
\label{eq:SimpDecom}
\overline{ \Omega}_H = \triangle_H\CG:
= \bigcup_{\tb\in \triangle_H} \tb \CG = \bigcup_{\sigma\in \CG}
    \left\{\tb\sigma:\tb\in \triangle_H\right\}.
\end{align}
Furthermore, the partition is non-overlapping, i.e.,  for any
$\tb,\sb\in \triangle_H$ and $\tb\neq \sb$, $\tb \CG \cap \sb \CG=
\emptyset$.

For any $\tb\in \RR_H^{d+1}$, the stabilizer of $\tb$ is denoted by
$\CG_{\tb}:=\left\{ \sigma \in \CG:  \tb\sigma=\tb\right\}$. Since $\CG$ is
finite, $|\CG_{\tb}| \times |\tb\CG| = |\CG|$. For the vertices of $\triangle_H$,
$|\mathbf{v}^k\CG|=\frac{(d+1)!}{k!(d+1-k)!}$ and $|\CG_{\mathbf{v}^k}| = k!(d+1-k)!$.
It is worthwhile to note that $\overline{ \Omega}_H$ is the convex hull of the
$\sum_{k=1}^{d}|\mathbf{v}^k\CG|= 2^{d+1}-2 $ vertices in
\begin{align*}
  \left \{ \frac{\vb^{k}\sigma}{d+1}:  \sigma\in \CG,  1\le k\le d \right\}.
\end{align*}

The generator matrix of the lattice  $A_d$ is given by the
$(d+1)\times d$ matrix
\begin{align*}
   A  : = \begin{pmatrix}
               1   &  0  & \dots & 0 &0\\
               0   &  1  & \dots & 0 &0\\
               \vdots & \vdots & \ddots & \vdots &\vdots\\
               0   &  0  & \dots & 1 & 0\\
               0   &  0  & \dots & 0 & 1\\
               -1  &  -1 & \dots & -1 &-1
              \end{pmatrix}
\end{align*}
and the generator matrix of the dual lattice $L^{\perp}_d$ is the $(d+1) \times d$
matrix
\begin{align*}
 A^{\perp} = A (A^\tr A)^{-1} =  \frac{1}{d+1}\begin{pmatrix}
               d   &  -1  & \dots & -1 &-1\\
               -1  &   d  & \dots & -1 &-1\\
               \vdots & \vdots & \ddots & \vdots &\vdots\\
               -1   & -1  & \dots &  d & -1\\
               -1   & -1  & \dots & -1 & d\\
               -1  &  -1  & \dots & -1 &-1
              \end{pmatrix}.
\end{align*}
To describe $L^\perp_d$ explicitly, we introduce the notation
\begin{align*}
 \HH : =\left\{\kb\in \ZZ^{d+1}_H:\   k_1\equiv k_2\equiv\cdots \equiv k_{d+1}
    \!\!\!\!\pmod{d+1}\right\}.
\end{align*}
For $j \in \ZZ^d$ set $\kb = (d+1) A^\perp j$. It follows readily that
$\kb \in \ZZ_H^{d+1}$ and a quick computation shows $\kb \in \HH$.
On the other kind, the definition of $A^\perp$ leads to $j = A^\tr \kb/(d+1)$,
which shows $j \in \ZZ^3$ if $\kb \in \HH$. Consequently, the dual lattice is
given by
\begin{align*}
    L^{\perp}_d =\left\{ \frac{\kb}{d+1}:\   \kb \in \HH  \right\}.
\end{align*}

For $\e^{2\pi i \kb \cdot \tb}$ with $\kb \in L^\perp_d$, we
introduce the new notation
\begin{align}\label{phi}
   \phi_{\jb}(\tb) : =  \e^{\frac{2\pi i}{d+1}  \jb\cdot \tb}  \qquad \hbox{for} \quad
     \jb \in \HH.
\end{align}
By the Fuglede theorem, $\left\{ \phi_{\jb}:\ \jb \in \HH  \right\}$
forms an orthonormal basis in $L^2(\Omega_H)$; that is,
\begin{align}
\label{eq:orth-H}
    \langle \phi_{\kb},  \phi_{\jb}\rangle
     = \delta_{\kb, \jb}, \quad \kb, \jb\in \HH,
\end{align}
where the inner product is defined by
\begin{align}
\label{eq:ip-H}
\langle f,g\rangle := \frac{1}{|\Omega_H|} \int_{\Omega_H} f(\tb) \overline{g(\tb)} d\tb
= \frac{1}{\sqrt{d+1}} \int_{\Omega_H} f(\tb) \overline{g(\tb)} d\tb.
\end{align}

\begin{defn} \label{H-periodic}
A function $f$ is $H$-periodic if it is periodic with respect to the
lattice $A_d$; that is, $f(\tb) = f(\tb + \kb )$ for all $\kb\in
\ZZ^{d+1}_H$.
\end{defn}

Evidently, the functions $\phi_\jb(\tb)$ in \eqref{phi} are $H$-periodic.
Furthermore, \eqref{eq:orth-H} shows that an $H$-periodic function $f$
can be expanded into a Fourier series
\begin{equation}\label{H-Fourier}
   f \sim \sum_{\kb \in \HH } \wh f_\kb \phi_\kb(\tb), \quad\hbox{where} \quad
       \wh f_\kb := \frac1{\sqrt{d+1}} \int_{\Omega_H} f(\tb) \phi_{-\kb}(\tb) d \tb.
\end{equation}

Recall that $\sb \equiv \tb \pmod{H}$ means $\sb - \tb \in \ZZ_H^{d+1}$.
The following lemma will be useful later.

\begin{lem} \label{lem:s=equiv=t} 
If $\tb,\sb \in \Omega_H$  and $\sb \equiv \tb \pmod{H}$, then $\tb=\sb$.
\end{lem}

\begin{proof}
If $\tb, \sb \in \Omega_H$ and $\tb \equiv \sb \pmod{H}$, then $\sb
- \tb \in \ZZ_H^{d+1}$ and, set $\kb := \sb-\tb$, $-1\leq  k_i  -k_j
\leq 1$  for all $1\le i,j\le {d+1}$. The last condition means that
either $k_i \in \left\{0,1 \right\}$ for all $ 1\le i \le {d+1}$ or
$ k_i \in \left\{0,-1\right\}$ for all $1 \le i \le {d+1}$. The
homogenous condition $\sum_{i=1}^{d+1}k_i=0$ then shows that $\kb=0$
or $\sb=\tb$.
\end{proof}


\subsection{Structure of the fundamental domain}
In this subsection, we concentrate on the structure of $\overline{\Omega}_H$.
For this purpose, we set $\NN_{d+1}:=\left\{ 1,2,\cdots,d+1  \right\}$ and start
with the observation that $\overline{\Omega}_H$ can be partitioned into $d+1$
congruent parts.

\begin{lem}\label{lm:Decomp}
For $1\le j\le {d+1}$ define
\begin{align*}
   \Omega_H^{\{j\}} := \left\{ \tb\in \RR_H^{d+1}:   0< t_i-t_j \le 1 \text{ and }
   0\le t_l-t_j < 1  \text{ for }  1\le  i< j <l \le d+1
   \right\}.
\end{align*}
Then
\begin{align*}
{\Omega}_H =  \bigcup_{1\le j\le {d+1}} {\Omega}_H^{\{j\}}
\qquad \hbox{and}\qquad
{\Omega}_H^{\{i\}}\cap {\Omega}_H^{\{j\}} = \emptyset \quad \text{for }
1\le i\neq j \le d+1.
\end{align*}
\end{lem}

\begin{proof}
For $\tb \in {\Omega}_H$, let $ J= \left\{ k\in \NN_{d+1}:   t_k\le
t_i \text{ for } i\in\NN_{d+1} \right\}$ and $j=\min_{k\in J} k$.
The definition shows that $t_i - t_j \ge 0$ for all $i \in
\NN_{d+1}$; furthermore, if $i<j < k$ then, as $\tb \in \Omega_H$, $
0 < t_i - t_j \le 1$ and $-1 < t_j-t_k \le 0$, so that $\tb \in
\Omega^{\{j\}}_H$, which implies that ${\Omega}_H \subseteq
\bigcup_{j\in \NN_{d+1}} {\Omega}_H^{\{j\}}$. Since
${\Omega}_H^{\{j\}} \subset {\Omega}_H$ for each $j\in \NN_{d+1}$,
it follows that ${\Omega}_H =\bigcup_{j\in \NN_{d+1}}
{\Omega}_H^{\{j\}}$.

If $\tb \in {\Omega}_H^{\{i\}}\cap{\Omega}_H^{\{j\}}$ with $1\le
i<j\le d+1$, then $\tb \in {\Omega}_H^{\{i\}}$ so that $0\leq
t_j-t_i$ and $\tb \in {\Omega}_H^{\{j\}}$ so that $0\le   t_i-t_j$.
As a consequence, $t_j<t_i\le t_j$, which implies
${\Omega}_H^{\{i\}}\cap{\Omega}_H^{\{j\}}=\emptyset$.
\end{proof}

\begin{figure}[htb]
\hfill\begin{minipage}[b]{0.48\textwidth} \centering
\includegraphics[width=0.9\textwidth]{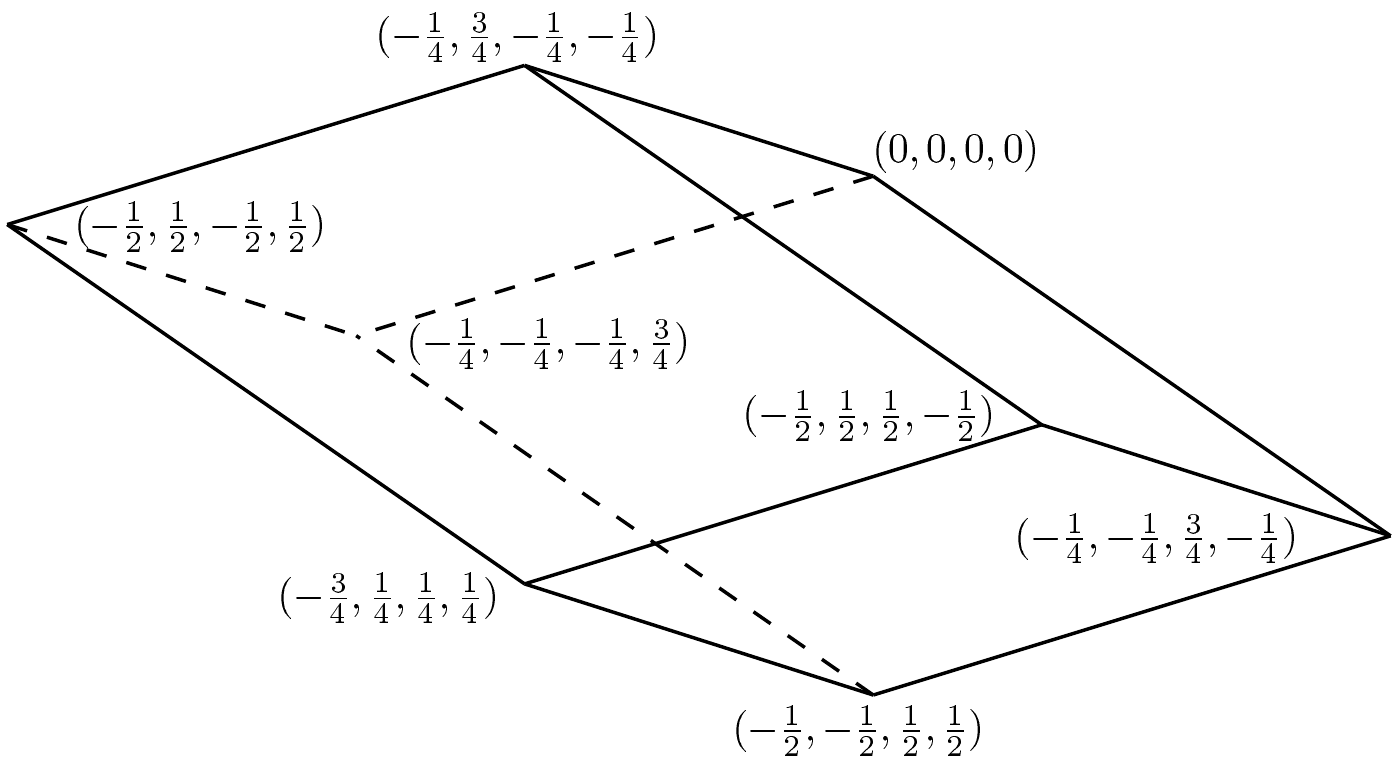}
\end{minipage}\hfill\begin{minipage}[b]{0.48\textwidth}
\centering
\includegraphics[width=0.75\textwidth]{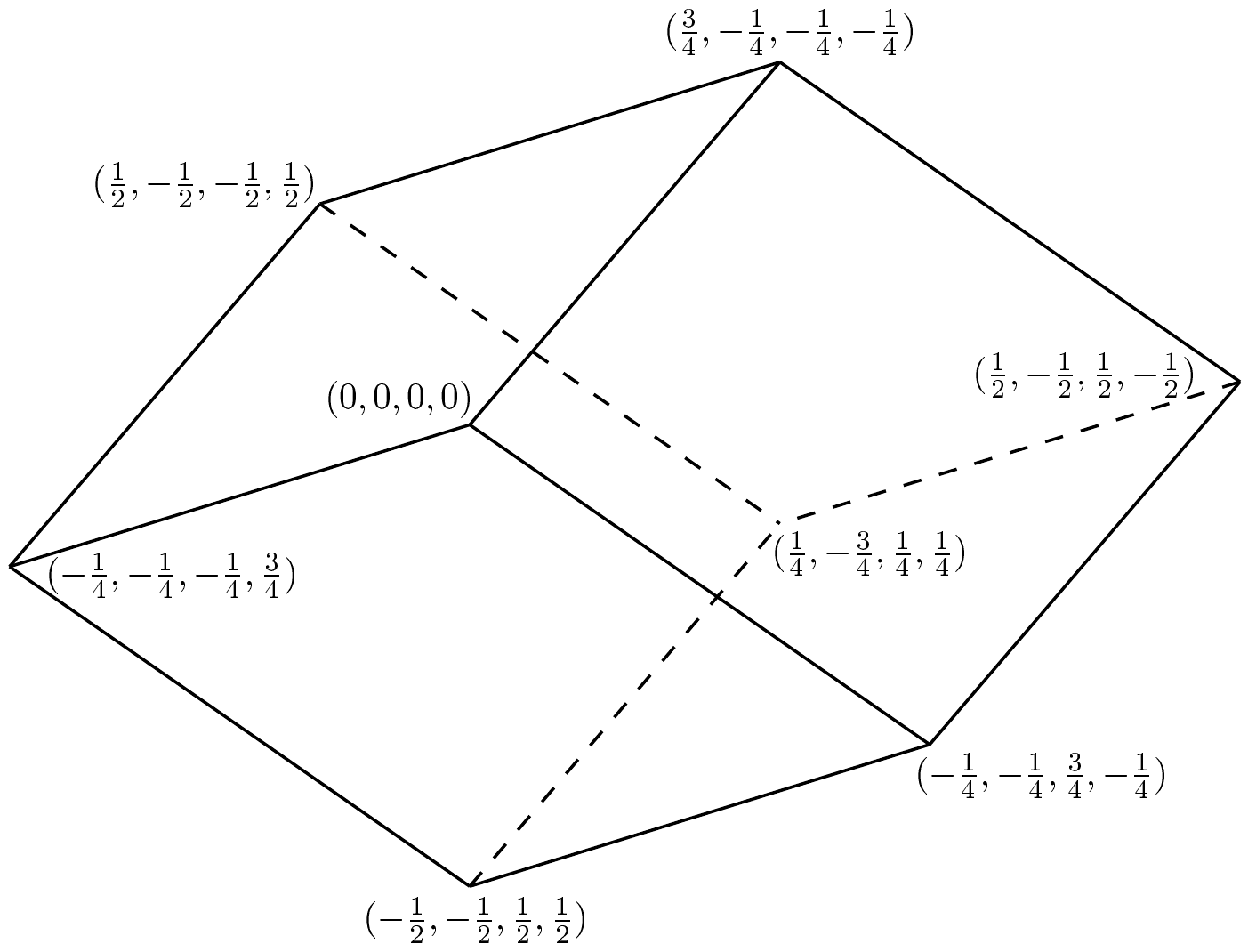}
\end{minipage}\hspace*{\fill}
\vspace*{1em} \hfill\begin{minipage}[b]{0.48\textwidth} \centering
\includegraphics[width=0.9\textwidth]{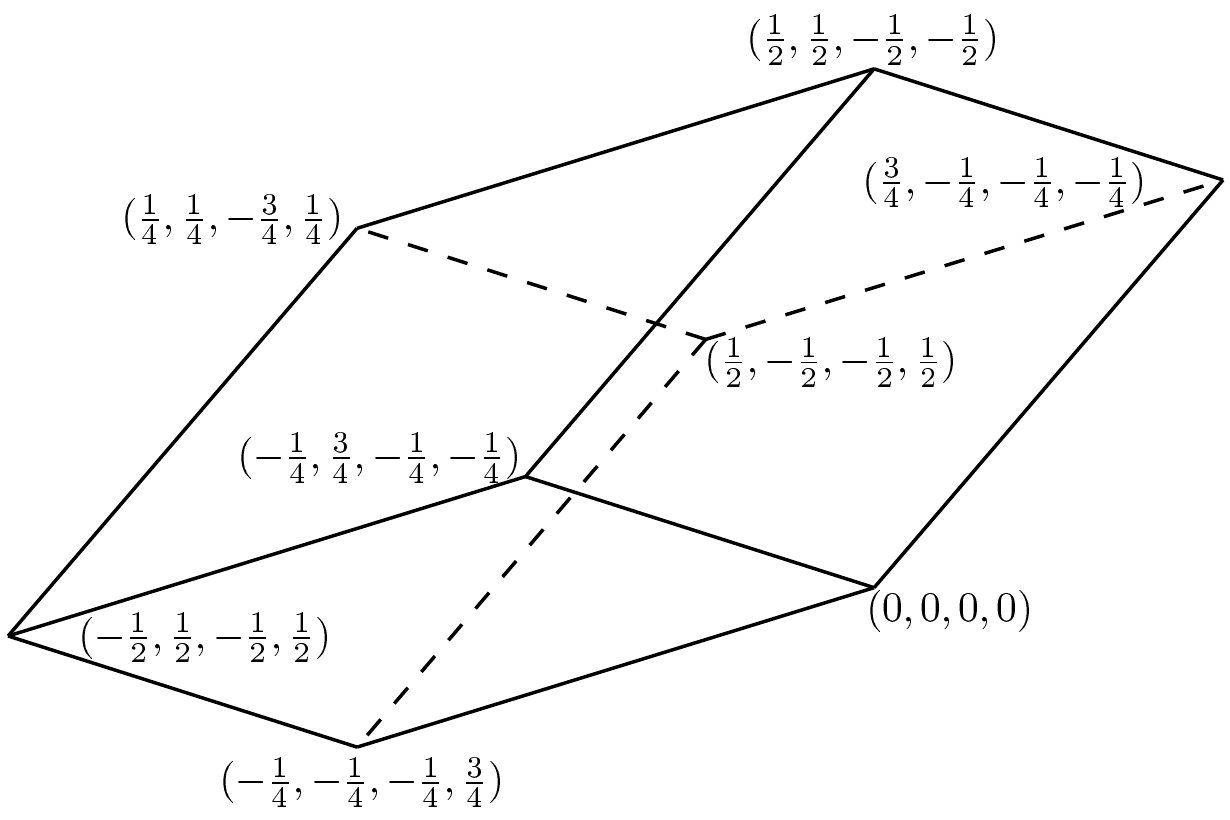}
\end{minipage}\hfill\begin{minipage}[b]{0.48\textwidth}
\centering
\includegraphics[width=0.75\textwidth]{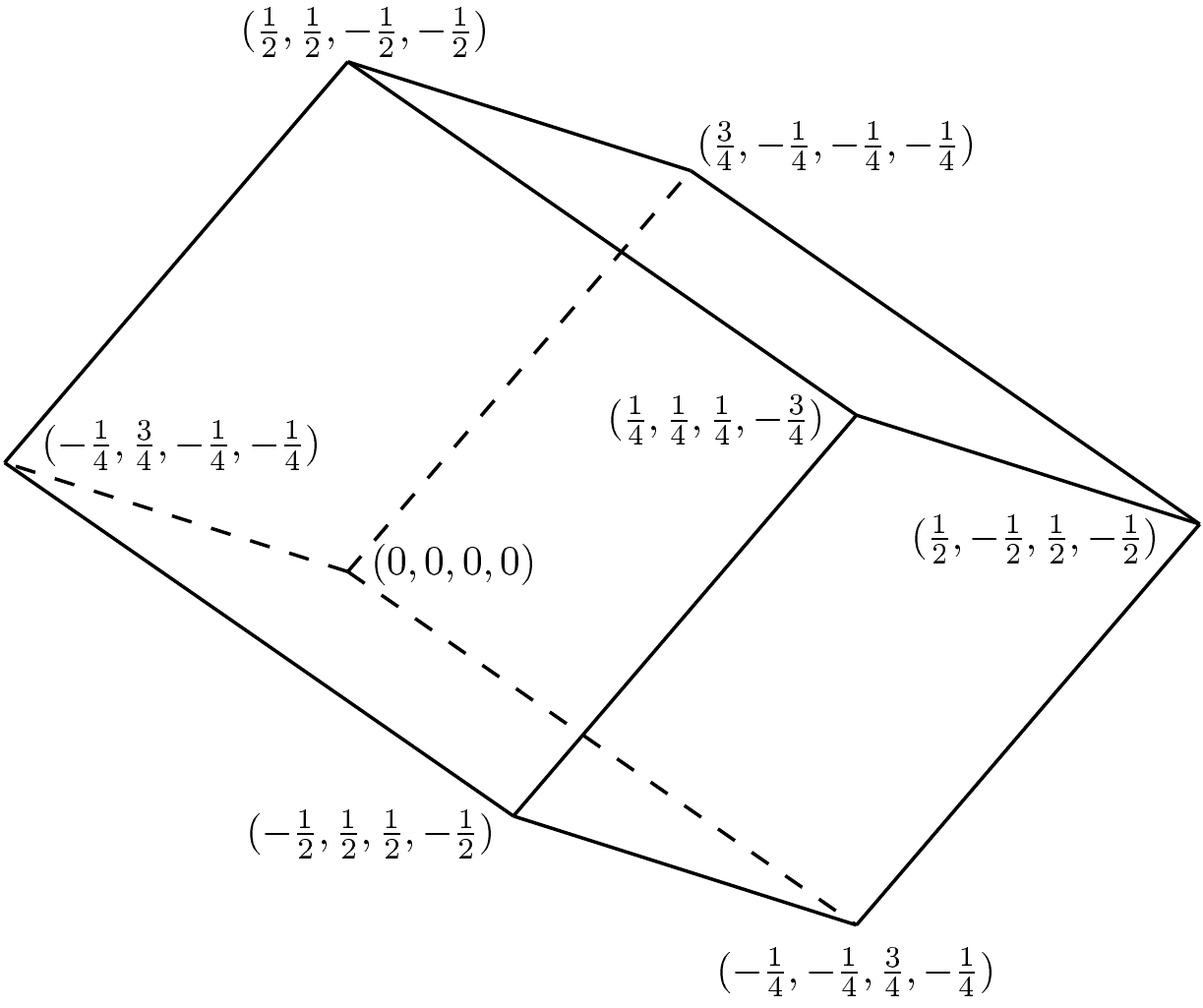}
\end{minipage}\hspace*{\fill}
\caption{Rhombohedral partition of the rhombic dodecahedron.
Upper-Left: $\Omega_H^{\{1\}}$; Upper-Right: $\Omega_H^{\{2\}}$;
Lower-Left: $\Omega_H^{\{3\}}$; Lower-Right: $\Omega_H^{\{4\}}$.}
\label{Partition3d}
\end{figure}

For $d=2$, the partition of the hexagon is evident from Figure \ref{hexagon}:
the sets $\Omega_H^{\{1\}},  \Omega_H^{\{2\}},  \Omega_H^{\{3\}}$ are the
three parallelograms inside the hexagon, starting from the one in the upper
left rotating clockwise. In the case $d=3$, the rhombic dodecahedron is more complicated (see Figure \ref{dodecahedron}), its decomposition is
depicted in the Figure \ref{Partition3d}.

The next lemma clarifies the intersection of the closure of $\Omega_n^{\{j\}}$,
for which we need to define $\overline{\Omega}_H^J$ for a set $J \subseteq
\NN_{d+1}$.

\begin{lem} \label{lm:Decomp*}
For $\emptyset\subset J\subseteq  \NN_{d+1}$, define
\begin{align*}
\overline{\Omega}_H^{J}:= \left\{\tb \in \RR_H^{d+1}: t_i = t_j, \, \forall i,j
\in J; \text{ and } \, 0\leq  t_i - t_j \leq 1,\, \forall j\in J, \
\forall i\in \NN_{d+1}\setminus J\right\}.
\end{align*}
Then
\begin{align*}
\overline{\Omega}_H =  \bigcup_{j\in \NN_{d+1}} \overline{\Omega}_H^{\{j\}}
\qquad \hbox{and}\qquad \overline{\Omega}_H^{J} = \bigcap_{j \in J}
\overline{\Omega}_H^{\{j\}}.
\end{align*}
\end{lem}

\begin{proof}
The first equation follows easily from the previous lemma by taking
closure.

If $\tb \in \overline{\Omega}^{\{i\}}_H\cap\overline{\Omega}^{\{j\}}_H$,
then $0\leq  t_i-t_j\leq 0$; that is, $t_i=t_j$.  Hence, if $\tb \in \bigcap_{j\in
J} \overline{\Omega}_H^{\{j\}}$ then $t_i = t_j, \, \forall i,j\in J$, which
implies $\bigcap_{j\in J} \overline{\Omega}_H^{\{j\}} \subseteq
\overline{\Omega}_H^{J}$. Since
$\overline{\Omega}_H^J \subset \overline{\Omega}_H^{\{j\}}$ whenever
$j\in J$ by definition, we conclude that $\overline{\Omega}_H^{J} =
\bigcap_{j\in J} \overline{\Omega}_H^{\{j\}}$.
\end{proof}

\subsection{Boundary of the fundamental domain}
In order to carry out the discrete Fourier analysis on the
fundamental domain  $\Omega_H$, we need to have a detailed knowledge
of the boundary of the region. Much of this subsection is parallel
to the study in \cite{LX08}, where the case $d =3$ is discussed in
detail. Some of the proofs, in fact, are essentially the same as in
the case of $d=3$ and often only minor adjustment is needed. In such
cases, we shall point out the necessary adjustment and refer the
proof to \cite{LX08}.

We use the standard set theoretic notations $\partial \Omega$,  $\Omega^\circ$
and $\overline{\Omega}$ to denote the boundary, the interior and the closure
of $\Omega$, respectively. Clearly $\overline{\Omega} = \Omega^\circ \cup
\partial \Omega$.  For $i,j \in \NN_{d+1}$ and $i\ne j$, define
$$
         F_{i,j}  = \{ \tb \in \overline{ \Omega}_H: t_i - t_j =1\}.
$$
There are a total $ d(d+1)$ distinct $F_{i,j}$, each stands for one
facet of $d-1$ dimension, together with its boundary, of the
fundamental domain.

The boundary of our $d$-dimensional polytope consists of lower
dimensional sets, which can be obtained from the intersections of
$d-1$ dimensional facets. For example, for $d=3$, the boundary
consists of faces, edges, and vertices; edges are intersections of
faces and vertices are intersection of edges. In order to describe the
intersections of facets, we define, for nonempty subsets $I, J$ of $\NN_{d+1}$,
\begin{align*}
 \Omega_{I,J} :=  \bigcap_{i\in I, j\in J} F_{i,j} =
   \left\{ \tb \in \overline{ \Omega}_H:\  t_j = t_i-1 \text{ for all } i\in I,  j\in J\right\}.
\end{align*}
The main properties of these sets are summarized in the following lemma,
the proof of which is given in \cite{LX08}.

\begin{lem}  \label{lm:boundary}
Let $I, J, I_i, J_i$ be nonempty subsets of $\NN_{d+1}$. Then
\begin{enumerate}
\item[(i)] $\Omega_{I,J} = \emptyset$ if and only if $I\cap J \neq \emptyset$.
\item[(ii)]  $ \Omega_{I_1,J_1} \cap \Omega_{I_2,J_2} =\Omega_{I, J}$  if
  $I_1 \cup I_2=I$ and $J_1\cup J_2=J$.
\end{enumerate}
\end{lem}

%

To carry out a discrete Fourier analysis, we need to distinguish points on
the boundary elements of different dimensions. It is also necessary to
distinguish a closed boundary element and an open one.
For example, for $d=3$, we will distinguish a face with its boundary edges
and a face without its edges. To make these more precise, we introduce the
notation
\begin{align*}
  &\CK = \left\{(I,J):  I, J \subset \NN_{d+1}; \   I\cap J = \emptyset \right\},\\
  &\CK_0 = \left\{ (I,J)\in \CK:\ i<j, \,\, \hbox{for all}\,\, (i,j) \in I\times J \right\}.
\end{align*}

\begin{defn}
For $(I,J)\in \CK$, the boundary element $B_{I,J}$ of the
fundamental domain,
\begin{align*}
  B_{I,J}: = \left\{\tb \in \Omega_{I,J}: \  \tb \not \in \Omega_{I_1,J_1}
\text{ for all }
  (I_1,J_1)\in \CK \text{ with } |I|+|J| < |I_1| + |J_1| \right\},
\end{align*}
is called a face of dimension $k$, or simply a $k$-face  if  $k=d+1-|I|-|J|$.
In particular, it is called a facet if $k=d-1$ $($or $|I|+|J|=2)$, a ridge if
$k = d-2$ $(|I|+|J|=3)$,
a face if $k=2$ $(|I|+|J|=d-1)$, an edge if $k=1$ $(|I|+|J|=d)$
and a vertex if $k=0$ $(|I|+|J|=d+1)$.
\end{defn}

\begin{figure}[ht]
\centering
\hfill\includegraphics[width=0.5\textwidth]{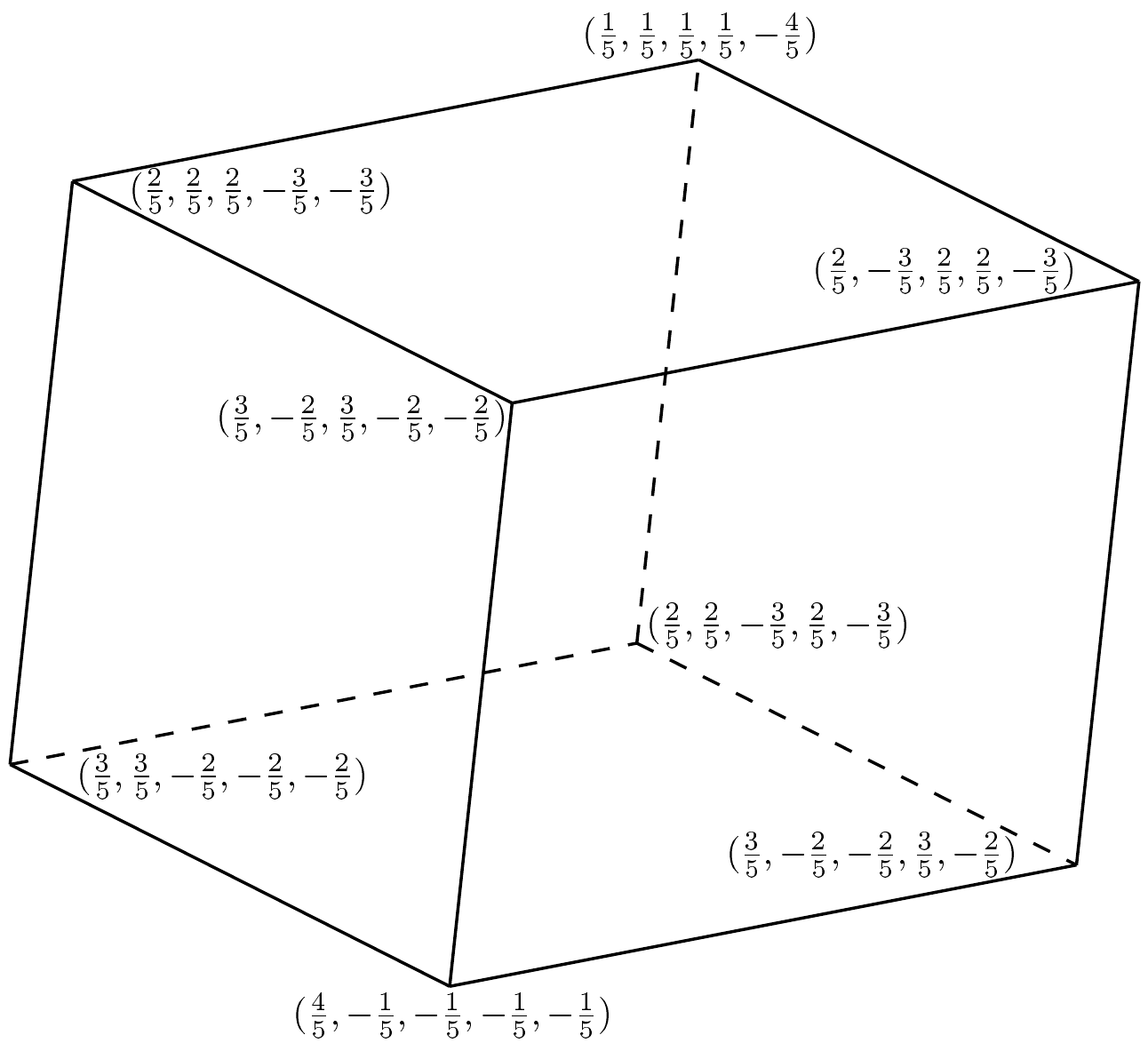}\hspace*{\fill}
\caption{The 3-face (facet) $t_1-t_5=1$ of the fundamental domain
when $d=4$.}\label{cell4d}
\end{figure}

Figure \ref{cell4d} shows a 3-face, or a facet, of the fundamental
domain when $d=4$.

In the following, when we say a $k$-face we mean the {\it
open} set, that is, without any of its boundaries. For  $k$-faces
with $k\ge 1$, the boundary elements $B_{I,J}$ represent the interiors.
In fact, it is easy to see that  $B_{\{i\},\{j\}} = F_{i,j}^\circ$ for
distinct integers $i, j \in \NN_{d+1}$ and
$B_{I,J}=\Omega_{I,J}^{\circ}$ for $(I,J)\in \CK$.
By the definition of $\Omega_H$ and $\Omega_{I,J}$, we can write
$B_{I,J}$ more explicitly as
\begin{align}
\label{eq:BIJ}
 B_{I,J} = \left\{\tb\in \RR_H^{d+1}:\
   t_i>t_l>t_j=t_i-1\text{ for all } i\in I,\, j\in J \text{ and }
l\in \NN_{d+1}\setminus (I\cup J)  \right\}.
\end{align}

The boundary elements of the same type will often be considered as
a group. For this purpose, we further define,  for $0 < i,j < i+j \le d+1$,
\begin{align} \label{CK}
\begin{split}
   \CK^{i,j}: =&\, \left\{(I,J) \in \CK:\  |I| = i,\ |J| =j  \right\}, \quad
       B^{i,j}: = \bigcup_{(I, J) \in \CK^{i,j}}  B_{I,J} \\
   \CK^{i,j}_0: = &\, \left\{(I,J) \in \CK_0:\  |I| = i,\ |J| =j  \right\}, \quad
       B^{i,j}_0: =\bigcup_{(I, J) \in \CK_0^{i,j}}  B_{I,J}.
\end{split}
\end{align}

\begin{prop} \label{prop:BIJ} 
Let $(I,J) \in \CK$ and $(I_1,J_1)\in \CK$.
\begin{enumerate}
\item[(i)]  $ B_{I,J} \cap B_{I_1,J_1} = \emptyset$,  if $I \neq I_1$ and $J\neq J_1$.
\item[(ii)] $\overline{\Omega}_H\setminus \Omega_H^{\circ}
 =  \bigcup_{(I,J)\in \CK}  B_{I,J}  = \bigcup_{0<i,j<i+j\le d+1} B^{i,j}$.
 \item[(iii)] $\Omega_H\setminus \Omega_H^{\circ}
 =  \bigcup_{(I,J)\in \CK_0}  B_{I,J}  = \bigcup_{0<i,j<i+j\le d+1} B_0^{i,j}$.
\end{enumerate}
\end{prop}


The proof of this proposition is essentially the same as the one given in
\cite{LX08} for $d=3$.

The boundary elements can also be described by using symmetry. In fact,
boundary elements of the same type can be transformed by the group
$\CG = S_{d+1}$. To make it more precise, we define, for $(I,J)\in \CK$
and $\sigma \in \CG$,
\begin{align*}
   B_{I,J}\sigma := \left\{\tb\sigma:\ \tb\in B_{I,J} \right\}.
\end{align*}
Then it follows readily that
\begin{align*}
     B^{|I|,|J|} =\bigcup_{\tb \in \G} B_{I,J}\sigma
     : = \left\{ \tb \sigma: \tb \in B_{I,J},  \sigma\in \mathcal{G} \right\}.
\end{align*}

Recall that $\sigma_{ij}= \tb - (t_i-t_j)\eb_{i,j}$ denotes the
transposition in $\mathcal{G}$ that interchanges $i$ and $j$. We
clearly have $\sigma_{ij} = \sigma_{ji}$ and $\sigma_{jj}$ is the
identity element of the group. For a nonempty set $I\subset
\NN_{d+1}$, define $\mathcal{G}_{I} := \left\{\sigma_{ij}:  i,j\in
I\right\}$. It is easy to verify that $\mathcal{G}_{I} $ forms a
subgroup of $\mathcal{G}= \mathcal{S}_{d+1}$ of order $|I|$.

\begin{lem} \label{lem:|BIJ|}
Let $\# B^{i,j}$ and $\# B^{i,j}_0$ denote the number of distinct $(d+1-i-j)$-faces
contained in $B^{i,j}$ and $B_0^{i,j}$, respectively. Then
$$
     \# B^{i,j} = \frac{(d+1)!}{i!j!(d+1-i-j)!} \qquad \hbox{and} \qquad
     \# B^{i,j}_0 = \frac{(d+1)!}{(i+j)!(d+1-i-j)!}.
$$
\end{lem}

\begin{proof}
By definition, $B^{i,j}$ and $B_0^{i,j}$ are unions of
$(d+1-i-j)$-faces of the fundamental domain.  The first formula
follows from the fact that $B_{I,J}\sigma=B_{I,J}$ if $\sigma\in
\CG_I \cup \CG_J \cup \CG_{\NN_{d+1}\setminus (I\cup J)}$. The
second one follows from the fact that , for any $K\subseteq
\NN_{d+1}$ and $i,j\ge 1$ such that $|K|= i+j$, there is a unique
$(I,J)\in \CK^{i,j}_0$ such that $I \cup J = K$.
\end{proof}

Using these formulas, we can derive by setting $k =  d+1-i-j$ that
the fundamental domain $\overline{\Omega}_H$ has a total
$$
\sum_{j = 1 }^{d-k} \frac{(d+1)!}{j!k!(d+1-j-k)!}=\big(2^{d+1-k}-2\big)\binom{d+1}{k},
      \qquad 0 \le k \le d-1,
$$
distinct $k-$faces, while $\Omega_H$ has a total $ (d-k)
\binom{d+1}{k}$ distinct $k-$faces.

For the periodic functions, we will need to consider points on the boundary
elements that are congruent modulus $\ZZ_H^{d+1}$.  For $(I,J) \subset \CK$
we further  define
\begin{align*}
   [B_{I,J}]:= \left\{ B_{I,J} + \kb : \kb \in \ZZ_H^{d+1}\right\}\cap
       \overline{\Omega}_H
     = \left\{ \tb+\kb\in \overline{\Omega}_H:  \tb \in B_{I,J},  \kb \in \ZZ_H^{d+1} \right\}.
\end{align*}
By lemma \ref{lem:s=equiv=t}, the set $[B_{I,J}]$ consists of exactly those
boundary elements that can be obtained from $B_{I,J}$ by congruent modulus
$\ZZ_H^{d+1}$.  As an example, in the case of $d=3$, we have
\begin{align} \label{B_{1,{2,3}}}
\begin{split}
& B_{\{1\},\{2,3\}} = \left \{(t, t-1, t-1, 2-3t):  \tfrac12 < t <\tfrac34  \right \}, \\
& B_{\{1,2\},\{3\}}= \left\{ (1-t, 1-t, -t, 3t-2):  \   \tfrac12 < t <\tfrac34  \right\},
\end{split}
\end{align}
and from this explicit description we can deduce, for example,
\begin{align*} 
  [B_{\{1\},\{2,3\}}] & = B_{\{1\},\{2,3\}} \cup \left (B_{\{1\},\{2,3\}} +(-1,1,0,0)\right)
     \cup \left (B_{\{1\},\{2,3\}} +(-1,0,1,0)\right)  \notag \\
       & = B_{\{1\},\{2,3\}} \cup B_{\{2\},\{1,3\}} \cup B_{\{3\},\{1,2\}}.
\end{align*}
The last equation indicates that $[B_{I,J}]$ is a union of $B_{I',J'}$, which
is stated and proved in \cite{LX08} for $d =3$ and the proof can be adopted
with obvious replacement of $d=3$ by general $d$.

\begin{lem}  \label{congruence} 
Let $(I,J)\in \mathcal{K}$. Then
\begin{align} \label{eq:[B]}
   [B_{I,J}] = \bigcup_{\sigma\in \mathcal{G}_{I\cup J}} B_{I,J}\sigma.
\end{align}
\end{lem}

Since $\CK^{i,j}$ can be obtained from $\CK_0^{i,j}$ from the action of $\CG$,
it follows that
\begin{equation} \label{eq:Bij=[B]}
  B^{i,j} = \bigcup_{(I,J) \in \CK_0^{i,j}} [B_{I,J}]= \bigcup_{B \in B_0^{i,j}} [B],
        \qquad  0 < i,j < i+j \le d+1.
\end{equation}
We also note that $[B_{I,J}]\cap [B_{I_1,J_1}] = \emptyset$ if $(I,J) \neq (I_1,J_1)$
for $(I,J)\in \mathcal{K}_0$ and $(I_1,J_1)\in \mathcal{K}_0$, which shows that
\eqref{eq:Bij=[B]} is a non-overlapping partition.

To illustrate the above partitions, let us consider the case $d =2$, for
which
 \begin{align*}
   &  B^{1,1} = [B_{\{1\},\{2\}}] \cup [B_{\{1\},\{3\}}] \cup [B_{\{2\},\{3\}}],
\\ &  B^{1,2} = [B_{\{1\},\{2,3\}}],\quad  B^{2,1} = [B_{\{1,2\},\{3\}}],
 \end{align*}
where
 \begin{align*}
  & B_{\{1\},\{2\}} = \left\{(1-t, 2t-1, -t):   \tfrac13< t<\tfrac23   \right\}, \\
  & B_{\{1\},\{3\}} =  B_{\{1\},\{2\}}\sigma_{2,3}, \quad   B_{\{2\},\{3\}}
       =  B_{\{1\},\{2\}}    \sigma_{1,2}, \\
  & B_{\{1\},\{2,3\}} = \left\{\left(\tfrac23, -\tfrac13,-\tfrac13\right)\right\},\quad
     B_{\{1,2\},\{3\}} = \left\{\left(\tfrac13,\tfrac13,-\tfrac23\right)\right\}.
 \end{align*}
In the case $d =3$, the boundary elements are given explicitly in \cite{LX08}.



\subsection{Point sets for discrete Fourier analysis for the lattice $A_d$}

We apply the general result on discrete Fourier analysis in Section
2 to the lattice $A_d$ by choosing $A=A$ and $B= nA$ with $n$ being
a positive integer.  Let $\ONE=(\{1\}^d)$. Then the matrix
\begin{align*}
  N : = B^{\tr}A = n I + n \ONE^\tr  \ONE =
    \begin{pmatrix} 2n & n & \dots & n & n\\ n & 2n  &\dots & n& n\\
  \vdots &\vdots  &\ddots &\vdots &\vdots \\
  n & n & \dots & 2n & n\\
  n & n & \dots & n & 2n
         \end{pmatrix}
\end{align*}
has integer entries. Since $N$ is now a  symmetric matrix, $\Lambda_N
 = \Lambda_{N^{\tr}}$; moreover, it is easy to see that $\Lambda_N^\dag =
 \Lambda_N$. We denote $\Lambda_N$ by $\HH_n$, which is given by
$$
\HH_n := \left\{ \kb \in \HH: \tfrac{\kb}{(d+1)n}\in
\Omega_H\right\} = \left\{\kb \in \HH:   -n < \frac{k_i -k_j}{d+1} \leq
n,\, 1\leq i < j \leq d+1 \right\}.
$$
The finite dimensional space $\mathcal{H}_N$ of exponentials in
Theorem \ref{thm:2.4} becomes
$$
\mathcal{H}_{n} :=\mathrm{span} \left\{\phi_{\kb}:\  \kb \in \HH_n\right\}
  \qquad \text{with}\quad
\mathrm{dim} \, \mathcal{H}_{n} = \det (N) = (d+1)n^{d}.
$$

\begin{figure}[htb]
\centering
\includegraphics[width=0.5\textwidth]{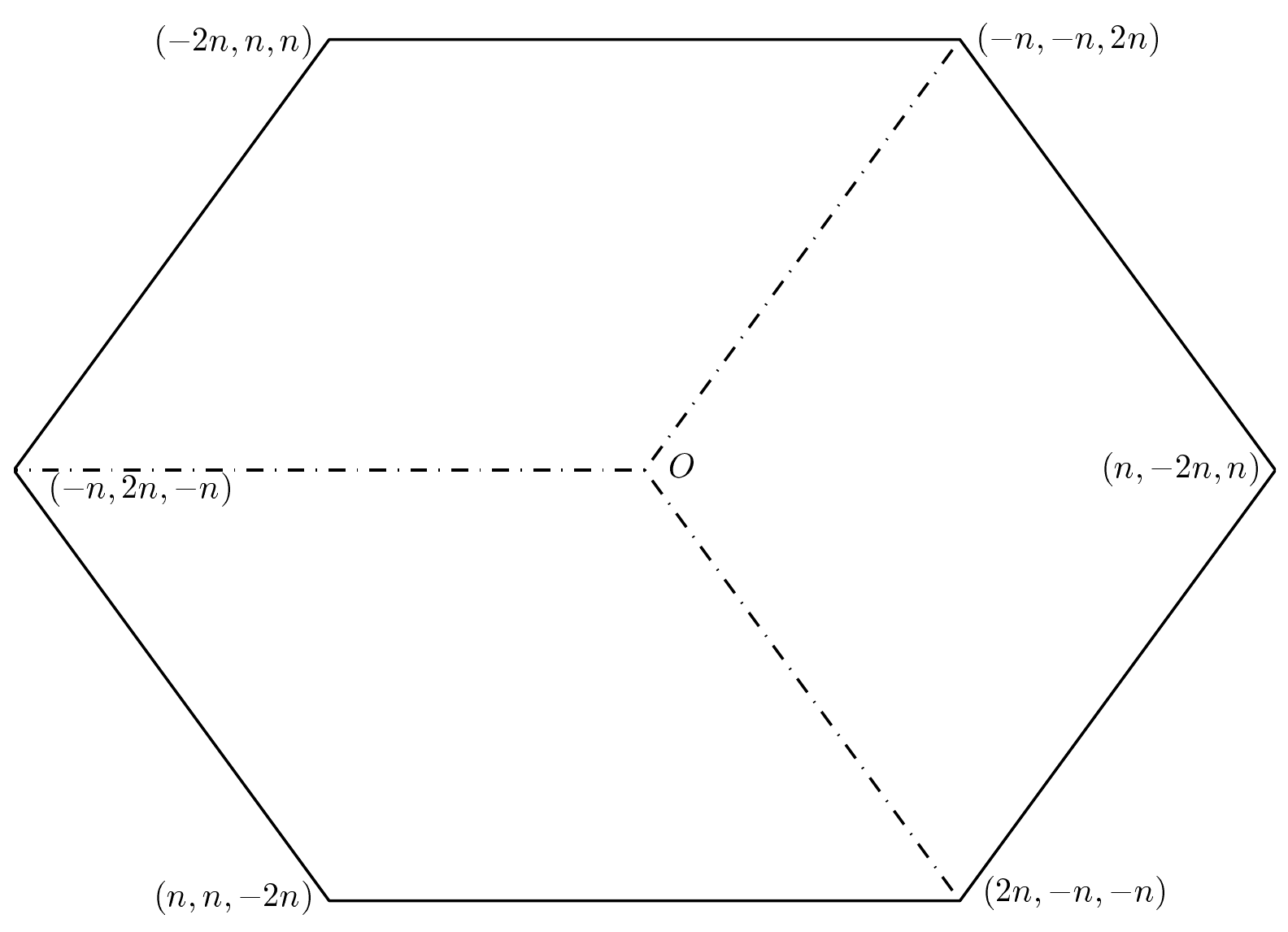}
\caption{$\HH_n^*$ for $d=2$.} \label{fig:hex2}
\end{figure}

\begin{figure}[htb]
\centering
\includegraphics[width=0.6\textwidth]{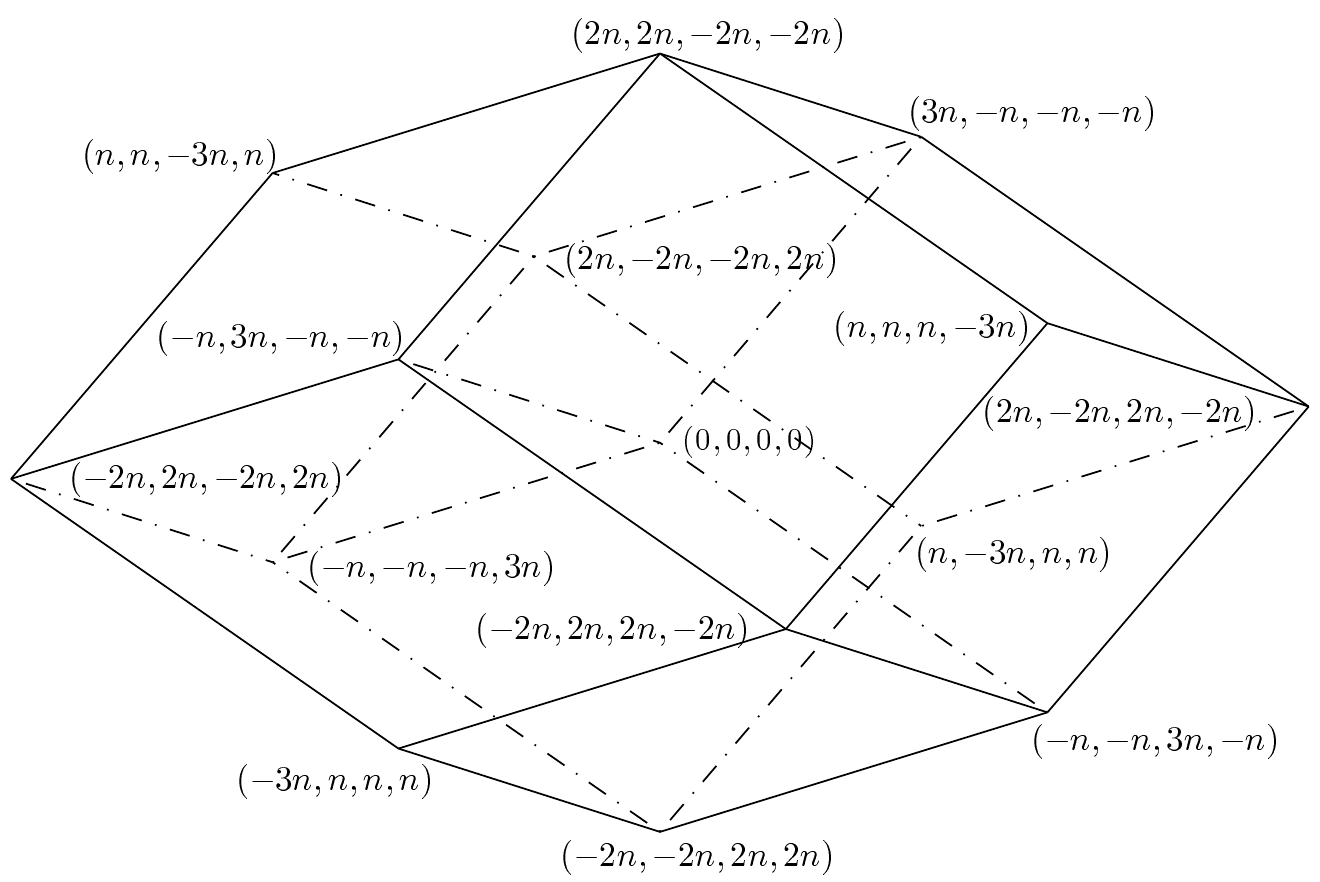}
\caption{$\HH_n^*$ for $d=3$.} \label{hfcc0}
\end{figure}

The points in $\HH_n$  are not symmetric under $\CG$, since points on
half of the boundary are not included. We further define the symmetric
extension of $\HH_n$ by
$$
      \HH^*_n :=\left\{ \kb \in \HH: \tfrac{\kb}{(d+1)n}\in \overline{\Omega}_H \right\}
= \left\{\kb \in \HH:   -n \leq \frac{k_i -k_j}{d+1} \leq
n,\, 1\leq i < j \leq d+1 \right\}.
$$
For $d =2$ and $3$, the scope of $\HH_n^*$ is depicted in Figure \ref{fig:hex2}
and Figure \ref{hfcc0}, in which the vertices are labeled in homogeneous
coordinates.

For discrete Fourier analysis, it is essential to understand the
sets $\HH_n$ and $\HH_n^*$, for which the main task lies in studying
the structure of points on the boundary; in other words, we need to
understand exactly when $\frac{\kb}{(d+1)n}$ belongs to a specific
boundary elements of $\Omega_H$. Much of the work in the previous
subsection is a preparation for this task. In the rest of this
subsection, we develop results in this regard.

We start with the partition of $\HH_n$ and $\HH_n^*$ into $d+1$ congruent
parts, each within a parallelepiped, corresponding to the partition of
$\Omega_H$ and $\overline{\Omega}_H$, as shown in Figures 3.3 for
$d =3$.

\begin{lem}\label{lm:DecompH}
For $1 \le j \le d+1$, define
\begin{align*}
   \HH_n^{\{ j\}} = \left\{ \kb\in \HH:   1\le \frac{k_i-k_j}{d+1} \le n \text{ and }
   0\le \frac{k_l-k_j}{d+1} \le n-1  \text{ for }  1\le  i< j <l \le d+1
   \right\}.
\end{align*}
Then
\begin{align*}
\HH_n=  \bigcup_{j\in \NN_{d+1}} \HH_n^{\{j\}} \qquad
\hbox{and}\qquad  \HH_n^{\{i\}} \cap \HH_n^{\{j\}} = \emptyset \quad
   \hbox{for $1 \le i \ne j \le d+1$}.
\end{align*}
\end{lem}

\begin{proof}
By definition, $\HH_n^{\{j\}} = \left\{\kb\in \HH:\   \frac{\kb}{(d+1)n}
\in \Omega_H^{\{j\}}  \right\}$. Hence, this lemma is an immediate
consequence of Lemma \ref{lm:Decomp}.
\end{proof}

Similarly, as a consequence of Lemma \ref{lm:Decomp*}, we have
the following lemma.

\begin{lem}\label{lm:Decomp*H}
For $\emptyset\subset J\subseteq  \NN_{d+1}$, define
\begin{align*}
\overline \HH_n^J:= \left\{\kb \in \HH: k_i = k_j, \, \forall i,j \in J;
\text{ and } \, 0\leq  k_i - k_j \leq (d+1)n,\, \forall j\in J, \
\forall i\in \NN_{d+1}\setminus J\right\}.
\end{align*}
Then
\begin{align*}
\HH_n^* =  \bigcup_{j\in \NN_{d+1}} \overline \HH_n^{\{j\}} \qquad
\hbox{and}\qquad \overline \HH_n^{J} = \bigcap_{j \in J} \overline \HH_n^{\{j\}}.
\end{align*}
\end{lem}

Next we consider further partitions of $\HH_n^*$.  The set of interior points
is defined by
\begin{align*}
  \HH_n^{\circ} := \left\{ \jb \in \HH:  \tfrac{\jb}{(d+1)n} \in  \Omega_H^{\circ} \right\}.
\end{align*}
The set of points on the boundary of $\overline \Omega_H$ is then
$\HH_n^* \setminus \HH_n^\circ$. Recall \eqref{CK}, we define, for $0<i,j<i+j\le d+1$, \begin{align} \label{eq:Hn^ij}
  \HH_n^{i,j} := \left\{ \kb \in \HH:\  \tfrac{\kb}{(d+1)n} \in B^{i,j}  \right\},\qquad
  \HH_{n,0}^{i,j} := \left\{ \kb \in \HH:\  \tfrac{\kb}{(d+1)n} \in B^{i,j}_0  \right\}.
\end{align}
The index set $\HH_n^{i,j}$ consists of those points $\jb$ in $\HH_n$
such that $\frac{\jb}{(d+1)n}$ are in the boundary element $B^{i,j}$ of
$\partial \Omega_H$. Using Proposition \ref{prop:BIJ}, it is easy to see that
$\HH_n^{i,j} \cap \HH_n^{i_1,j_1} = \emptyset$ if $i\ne i_1$ or $j \ne j_1$,
\begin{align} \label{eq:Hnij}
\bigcup_{0<i,j<i+j\le d+1} \HH_n^{i,j}   = \HH_n^*\setminus \HH_n^{\circ} \qquad
\hbox{and} \qquad
\bigcup_{0<i,j<i+j\le d+1} \HH_{n,0}^{i,j}   = \HH_{n}\setminus \HH_n^{\circ}.
\end{align}

\begin{lem} 
For $n \ge 1$, $0<i,j<i+j\le d+1 $,
$$
  |\HH_n^\circ| = n^{d+1} - (n-1)^{d+1}, \qquad |\HH_n^{i,j}|=
     \frac{(d+1)!}{i!j!(d+1-j)!} (n-1)^{d+1-i-j}.
$$
\end{lem}

\begin{proof}
The first equation follows from $|\HH_n^\circ| = |\HH_{n-1}^*|$ and
\eqref{|Hn*|} below. As shown in Lemma \ref{lem:|BIJ|}, $B^{i,j}$
consists of $\frac{(d+1)!}{i!j!(d+1-j)!}$ distinct boundary elements
$B_{I,J}$ with $(I,J)\in \CK$, $|I|=i$ and $|J|=j$, which are
$(d+1-i-j)$-faces. Each of these faces contains $(n-1)^{d+1-i-j}$
distinct points in $\HH_n^*$ as seen by \eqref{eq:BIJ}, from which
the second equation follows.
\end{proof}

We will not need to define the point set that corresponds to $B^{I,J}$ for
each pair $(I, J) \subset \CK^{i,j}$. However, we do need to understand the
congruent relation modulus $\ZZ_H^{d+1}$ among points in $\HH_n^*$. For each
$\jb \in \HH_{n,0}^{i,j}$, we define
\begin{equation}\label{CSk}
  \CS_\jb := \left \{\kb \in \HH_n^*:  \tfrac{\kb}{(d+1)n} \equiv
     \tfrac{\jb}{(d+1)n} \pmod{\ZZ_H^{d+1}} \right \}.
\end{equation}
Since, by definition, $[B_{I,J}] = \{\sb : \sb \equiv \tb \pmod{\ZZ_H^{d+1}},
\tb \in B_{I,J}\}$,
it follows from \eqref{eq:Bij=[B]} that
$$
\HH_n^{i,j} = \bigcup_{\kb\in \HH^{i,j}_{n,0}}\CS_{\kb}.
$$
By \eqref{eq:[B]}, $[B_{I,J}]$ is the union of $\binom{|I|+|J|}{|I|}$ components
of $B_{I,J}\sigma$. It follows that $|\CS_\jb| = \binom{k+l}{l}$ for
$\jb \in \HH_{n,0}^{k,l}$.

\subsection{$A_d$ Fourier partial sum}
For the Fourier expansion \eqref{H-Fourier} of an $H$-periodic
function, we define its partial sum with respect to the fundamental
domain of $A_d$ as
\begin{align} \label{PartialSum}
S_n f (\tb) := \sum_{\kb\in \HH_n^*} \langle f, \phi_{\kb}  \rangle
 \phi_{\kb}(\tb)
= \frac{1}{\sqrt{d+1}}\int_{\Omega_H} f(\sb) D_n^H(\tb-\sb) d\sb,
\end{align}
where $D_n^H$ is the Dirichlet kernel for the Fourier partial sum
\begin{align} \label{def:Dn^H}
D_n^H(\tb) := \sum_{\kb\in \HH_n^*}  \e^{\frac{2\pi i}{d+1}\,
\kb\cdot \tb}.
\end{align}
The following theorem gives a compact formula for the Dirichlet
kernel.

\begin{thm} \label{H-partial} 
For $ n \ge 0$,
\begin{align} \label{DnTheta-n}
    D_n^H(\tb) =\Theta_{n+1}(\tb) - \Theta_{n}(\tb), \qquad \hbox{where}\quad
      \Theta_{n}(\tb): = \prod_{j=1}^{d+1} \frac{\sin \pi n t_j}{\sin \pi t_j}.
\end{align}
\end{thm}

This theorem is stated and proved in \cite{LX08} for $d =3$ and the
proof carries over to general $d$ with the obvious replacement of
$d+1=4$ by general $d+1$.  Recently a different derivation of this formula
is given in \cite{SunX}.


As one consequence of the compact expression in \eqref{DnTheta-n}, we obtain
\begin{equation} \label{|Hn*|}
           |\HH^*_n| = D_n^H(0) = (n+1)^{d+1}-n^{d+1}.
\end{equation}
Another result that follows from the compact formula of the Dirichlet kernel
is an upper bound for the Lebesgue constant, which is the norm of the partial
sum  $S_n f$ in \eqref{PartialSum}.  Let $\|f\|_\infty$ denote the uniform norm
of $f \in C(\overline{\Omega}_H)$ and let $\|S_n\|_\infty$ denote the operator
norm of $S_n: C(\overline{\Omega}_H) \mapsto C(\overline{\Omega}_H)$.

\begin{thm} \label{SnNorm} 
There are positive constants $c$ and $C$ independent of $f$ and $n$
such that
\begin{equation} \label{eq:SnNorm}
   c(\log n)^d \le \|S_n \|_{\infty} \le C (\log n)^d.
\end{equation}
\end{thm}

The upper bound $\|S_n\|_\infty \le C (\log n)^d$ was proved in
\cite{LX08} for $d =3$, the proof extends to the general $d$. It
turns out, however, the upper bound can be derived from a general
result in \cite{Po}, which establishes the inequality
\begin{equation} \label{eq:SnNorm2}
     \int_{[-\pi,\pi]^n}\bigg|\sum_{k\colon k/n\in E}e^{i k \cdot x}\bigg| dx
       \leq C  (\log n)^d
\end{equation}
for $E$ being a polyhedron in $\RR^d$. Indeed, choosing a constant
$\alpha$, if necessary, the set $\alpha \Omega$ can be imbedded inside
$[- \pi,\pi]^d$ and $\alpha$ can be absorbed into the index set, so that
\eqref{eq:SnNorm} follows from the general result. The lower bound that
matches \eqref{eq:SnNorm2} was established in \cite{V} for $E$ being a
convex polyhedron in $\RR^d$. Upong choosing an appropriate dilation
constant $\alpha$ so that $\alpha \Omega$ contains $[- \pi,\pi]^d$, the
lower bound in \eqref{eq:SnNorm} follows.

\subsection{Discrete Fourier analysis on the fundamental domain}
For the lattice $A_d$, the general result on the discrete inner
product, Theorem \ref{thm:2.4}, gives the following proposition in
homogeneous coordinates:

\begin{prop}\label{pro:inner}
For $n\geq 0$, define
\begin{align*}
     \langle f,\, g \rangle_n := \frac{1}{(d+1)n^{d}} \sum_{j\in \HH_n} f(\tfrac{\jb}{(d+1)n})
        \, \overline{g(\tfrac{\jb}{(d+1)n})},      \quad f,\, g\in C(\overline{\Omega}_H).
\end{align*}
Then
$$
\langle f,\, g \rangle = \langle f,\, g \rangle_n, \qquad f,\, g \in \CH_n.
$$
\end{prop}

The inner product $\langle\cdot, \cdot\rangle_n$ is defined as a sum over
the index set $\HH_n$, which is not symmetric over $\overline{\Omega}_H$
as only part of the points in the boundary of $\overline{\Omega}_H$ are
accounted for. More interesting to us is to consider an inner product based
on the symmetric index set $\HH^*_n$, which turns out to be equivalent to
$\langle \cdot, \cdot \rangle_n$ for $H$-periodic functions.

\begin{defn} \label{defn:Sym-ipd} 
For $n \ge 0$ define the symmetric discrete inner product
\begin{align*}
 \langle f,\, g \rangle_n^* : =  \frac{1}{(d+1)n^{d}} \sum_{\jb\in \HH^*_n} c_{\jb}^{(n)}
 f(\tfrac{\jb}{(d+1)n})\, \overline{g(\tfrac{\jb}{(d+1)n})}, \qquad
    f,\, g \in C(\overline{\Omega}_H),
\end{align*}
where $c_{\jb}^{(n)}=1$ if $ \jb \in \HH_n^{\circ}$, and
$c_{\jb}^{(n)}=\frac{1}{\binom{i+j}{i}}$ if $ \jb \in \HH_n^{i,j}$.
\end{defn}

For instance, if $d=2$ then $\Omega_H$ is a regular hexagon and we have
\begin{align*}
    c_{\jb}^{(n)} =\begin{cases} 1, & \jb \in \HH_n^{\circ}, \quad
               \qquad \qquad \!\!\! (\hbox{$n^{3}-(n-1)^{3}$ points in the interior}),\\
       \frac12, & \jb \in \HH_n^{1,1}, \quad
             \quad\qquad (\hbox{$6\times (n-1)$ points on the edges}),\\
       \frac13, & \jb \in \HH_n^{1,2}\cup \HH_n^{2,1}, \quad
             (\hbox{$6$ vertices});
\end{cases}
\end{align*}
if  $d=3$, then $\Omega_H$ is the rhombic dodecahedron and we have
\begin{align*}
    c_{\jb}^{(n)} =\begin{cases} 1, & \jb \in \HH_n^{\circ}, \quad
               \qquad \qquad \!\!\! (\hbox{$n^{4}-(n-1)^{4}$ points in the interior}),\\
       \frac12, & \jb \in \HH_n^{1,1}, \quad
             \quad\qquad (\hbox{$12 (n-1)^2$ points on the faces}),\\
       \frac13, & \jb \in \HH_n^{1,2}\cup \HH_n^{2,1}, \quad
             (\hbox{$2 \times 12 (n-1) $ points on the edges}),\\
       \frac14, & \jb \in \HH_n^{1,3}\cup \HH_n^{3,1}, \quad
             (\hbox{$2\times 4$  points on the vertices}),\\
       \frac16, & \jb \in \HH_n^{2,2}, \qquad  \!  \qquad (\hbox{$6$ points on the vertices}).
\end{cases}
\end{align*}

It is easy to verify that
\begin{align*}
&\sum_{\jb \in \HH_n^*\backslash \HH_n^{\circ}} c_{\jb}^{(n)} = \sum_{0<i,j<i+j\le d+1}  \frac{\big|\HH^{i,j}_n\big|}{\binom{i+j}{i}}
\\ &\qquad =  \sum_{0<i,j<i+j\le d+1} \binom{d+1}{i+j}(n-1)^{d+1-i-j}
 =  \sum_{k=2}^{d+1} (k-1)\binom{d+1}{k}(n-1)^{d+1-k}
\\ &\qquad =  (d+1)\sum_{k=1}^{d}\binom{d}{k}(n-1)^{d-k} -
   \sum_{k=2}^{d+1} \binom{d+1}{k} (n-1)^{d+1-k}
\\ &\qquad =  (d+1)\big(n^{d} - (n-1)^{d}\big)  - \big(n^{d+1}-(n-1)^{d+1}-(d+1)(n-1)^{d}   \big)
\\ &\qquad = (d+1)n^{d} - n^{d+1}+(n-1)^{d+1},
\end{align*}
so that $\la 1, 1 \ra_n^*=\frac{1}{(d+1)n^{d}} \big(|\HH_n^{\circ}| + (d+1)n^{d} - n^{d+1}+(n-1)^{d+1}\big) = 1$ by \eqref{|Hn*|}.

\begin{thm} \label{ipdH}
 For $n\geq 0$,
\begin{align*}
   \langle f,\, g\rangle = \langle f,\,g \rangle_n = \langle f,\,g \rangle_n^*,
     \quad f,\,g\in \mathcal{H}_n.
 \end{align*}
\end{thm}

\begin{proof}
 Let $f$ be  an  $H$-periodic function. Then
\begin{align*}
\sum_{\jb\in \HH^*_n\setminus \HH_n^{\circ}} c_{\jb}^{(n)}
f(\tfrac{\jb}{(d+1)n}) &= \sum_{0<i,k<i+k\leq d+1}
\frac{1}{\binom{i+k}{i}} \sum_{\jb \in \HH_n^{i,k} }
      f(\tfrac{\jb}{(d+1)n}) \\
& =  \sum_{0<i,k<i+k\leq d+1} \frac{1}{\binom{i+k}{i}} \sum_{\jb \in
\HH_{n,0}^{i,k} }
     \sum_{\kb \in \CS_\jb} f(\tfrac{\kb}{(d+1)n}).
\end{align*}
Since $|\CS_\jb| = \binom{i+k}{i}$ for $\jb \in \HH_{n,0}^{i,k}$,  using the
invariance of $f$, we then conclude that
$$
 \sum_{\jb\in \HH^*_n\setminus \HH_n^{\circ}} c_{\jb}^{(n)} f(\tfrac{\jb}{(d+1)n})
  =  \sum_{0<i,k<i+k\leq d+1} \frac{1}{\binom{i+k}{i}} \sum_{\jb \in \HH_{n,0}^{i,k} }
       \binom{i+k}{i} f(\tfrac{\jb}{(d+1)n}) =
        \sum_{\jb\in \HH_n\setminus \HH_n^{\circ}}  f(\tfrac{\jb}{(d+1)n}).
$$
Since $c_\jb^{(n)}=1$ if $\jb \in \HH_n^\circ$, the proof is completed.
\end{proof}

Since Theorem \ref{ipdH} shows that the integral of $f \overline{g} \in \CH_n$
agrees with the discrete sum over $\HH_n^*$, it is not surprising that we have
a cubature formula, which turns out to have a high degree of precision. To be
more precise, we define by $\CT_n$ the space of generalized trigonometric
polynomials,
$$
\CT_n : = \mathrm{span} \left\{\phi_{\kb}:   \kb\in \HH^*_n\right\}.
$$

\begin{thm}\label{th:cubature-H} 
For $n\ge 0$, the cubature formula
\begin{align} \label{cuba-trig}
   \frac{1}{\sqrt{d+1}} \int_{\Omega} f(\tb) d{\tb} = \frac{1}{(d+1)n^{d}} \sum_{\jb \in
     \HH_n^*}   c_{\jb}^{(n)} f(\tfrac{\jb}{(d+1)n})
\end{align}
is exact for all $f\in \CT_{2n-1}$.
\end{thm}

\begin{proof} It suffices to prove that \eqref{cuba-trig} is exact
for $f(\tb)=\phi_{\kb}(\tb)$ with any $\kb\in \HH^*_{2n-1}$. Since
$\Omega_H$ tiles $\RR_H^{d+1}$ with the lattice $A_d$, there exit
$\sb\in \Omega_H$ and $\lb\in \ZZ^{d+1}_H$ such that
$\frac{\kb}{(d+1)n} = \sb+\lb$. Noting that $\kb\in \HH^*_{2n-1}$
and $\lb\in \ZZ^{d+1}_H$ are vectors with integer entries, we
further deduce that $\sb=\frac{\mb}{(d+1)n}$ for certain $\mb\in
\HH_n$, which states that $\kb = \mb+(d+1)n\lb$ with certain $\mb\in
\HH_n$ and $\lb\in \ZZ^{d+1}_H$. Next we show that $\mb=0$ only if
$\kb=0$. In fact $\kb=(d+1)n\lb\in \HH^*_{2n-1}$ clearly states that
$\frac{n\lb}{2n-1}\in \Omega_H$. Then the restriction $\frac1n-2
<l_i-l_j \le 2-\frac1n$ for $1\le i< j\le d+1$ and the homogeneity
of $\lb$ imply that $\lb=0$, which gives $\kb=0$ in return. Now
applying the periodicity of $\phi_{\jb}$, Theorem \ref{ipdH} and
\eqref{eq:orth-H} yields that
\begin{align*}
\frac{1}{(d+1)n^{d}} \sum_{\jb \in
     \HH_n^*}  c_{\jb}^{(n)} \phi_{\kb}(\tfrac{\jb}{(d+1)n})
     & = \frac{1}{(d+1)n^{d}} \sum_{\jb \in
     \HH_n^*}   c_{\jb}^{(n)} \phi_{\mb}(\tfrac{\jb}{(d+1)n})
     \\ 
     & =  \langle \phi_{\mb}, \phi_{0} \rangle
     = \delta_{\mb,0} = \delta_{\kb,0}
     = \frac{1}{\sqrt{d+1}} \int_{\Omega} \phi_{\kb}(\tb) d{\tb}.
\end{align*}
This completes the proof.
\end{proof}

The cubature \eqref{cuba-trig} is an analogue of the classical quadrature
formula for trigonometric polynomials of one variable (see, for example,
\cite[Vol. II, p. 8]{Z}), which holds for trigonometric polynomials of one
variable and has vast applications.


\subsection{Interpolation on the fundamental domain}
For the fundamental domain of $A_d$, the general result on the
interpolation, Theorem \ref{thm:interpolation}, becomes  the
following:

\begin{prop}
For $n> 0$, define
\begin{align*}
   \mathcal{I}_n f(\tb) := \sum_{\jb\in \HH_n} f(\tfrac{\jb}{(d+1)n})
           \Phi_n (\tb-\tfrac{\jb}{(d+1)n}) ,\quad
   \text{ where }\quad  \Phi_n (\tb) = \frac{1}{(d+1)n^{d}} \sum_{\kb \in \HH_n} \phi_k (\tb),
\end{align*}
for $f\in C(\overline{\Omega}_H)$. Then $\mathcal{I}_nf\in
\mathcal{H}_n$ and
\begin{align*}
    \mathcal{I}_nf(\tfrac{\jb}{(d+1)n}) = f(\tfrac{\jb}{(d+1)n}), \quad \forall \jb \in \HH_n.
\end{align*}
\end{prop}

The function $\CI_n f$ interpolates on points in $\HH_n$, which is not symmetric
as noted before. We are more interested in interpolation on all points in $\HH_n^*$.
The operator $\CI_n^*$ that we are able to define, however, does not interpolate
at all points in $\HH_n^*$.  On the other hand, $\CI_n^* f $ can be used to
derive an truly interpolation function for points on the fundamental simplex, which
will be developed in the next section. Recall $\CS_\kb $ defined in \eqref{CSk}.

\begin{thm} \label{interp-H} 
For $n\geq 0$ and $f\in C(\overline{\Omega}_H)$, define
\begin{align*}
    \mathcal{I}_n^* f(\tb) := \sum_{\jb\in \HH^*_n} f(\tfrac{\jb}{(d+1)n}) \ell_{\jb,n} (\tb),
\end{align*}
where
\begin{align*}
  \ell_{\jb,n}(\tb) = \Phi^*_n(\tb-\tfrac{\jb}{(d+1)n}) \quad \text{ and }\quad
  \Phi^*_n(\tb) = \frac{1}{(d+1)n^{d}} \sum_{\kb \in \HH^*_n} c_\kb^{(n)} \phi_\kb (\tb).
\end{align*}
Then $\mathcal{I}^*_nf\in \mathcal{T}_n$ and it satisfies
\begin{align}    \label{eq:sym-interp}
         \mathcal{I}^*_n f(\tfrac{\jb}{(d+1)n}) = \begin{cases}
         f(\frac{\jb}{(d+1)n}), & \jb\in \HH^{\circ}_n,\\
         \\
         \displaystyle \sum_{\kb\in S_{\jb}}  f(\tfrac{\kb}{(d+1)n}), & \jb\in \HH^*_n
                 \setminus \HH_n^{\circ}.
\end{cases}
\end{align}
Furthermore, $\Phi_n^*(\tb)$ is a real function and it satisfies
\begin{align} \label{eq:Phi-n*}
 \Phi_n^*(\tb) = & \frac{1}{(d+1)n^{d}}    \sum_{j=1}^{d+1}
\bigg(\prod_{i=1\atop i\neq j }^{ d+1}
    \frac{\sin \pi n t_i}{\sin \pi t_i}  \bigg)      \cos \pi n t_j
   \sum_{I\subseteq\NN_{d+1}^{\{j\}} } \frac{|I|!(d-|I|)!}{(d+1)!}
\, \cos \pi \big( t_j+2\sum_{i\in I} t_i\big),
\end{align}
where $\NN_{d+1}^{\left\{j\right\}} := \NN_{d+1}\setminus\{j \}$.
\end{thm}

\begin{proof}
By definition,
\begin{align*}
\ell_{\jb,n}(\tfrac{\kb}{(d+1)n}) =
\Phi_n^*(\tfrac{\kb-\jb}{(d+1)n}) =
    \frac{1}{(d+1)n^{d}} \sum_{\lb \in \HH_n^*} c_\lb^{(n)} \phi_\lb(\tfrac{\kb-\jb}{(d+1)n}).
\end{align*}
Since $\Omega_H$ tiles $\RR_H^{d+1}$, there exist
$\mb,\lb\in\ZZ^{d+1}_H$ such that $\frac{\mb}{(d+1)n} \in \Omega_H$
and $\kb-\jb=\mb+(d+1)n\lb$. Thus, by Theorem \ref{ipdH},
\begin{align*}
\ell_{\jb,n}(\tfrac{\kb}{(d+1)n}) =&\frac{1}{(d+1)n^{d}}\sum_{\ib\in
\HH^*_n}
 c_{\ib}^{(n)} \phi_{\ib}(\tfrac{\mb}{(d+1)n}) = \frac{1}{(d+1)n^{d}}\sum_{\ib\in \HH^*_n}
    c_{\ib}^{(n)} \phi_{\mb}(\tfrac{\ib}{(d+1)n}) \\
 =& \langle \phi_{\mb}, \phi_{0}\rangle^*_n =
 \langle \phi_{\m}, \phi_{0}\rangle = \delta_{\mb,0}.
\end{align*}
Equivalently we can write the above equation as
\begin{align} \label{k=jmod}
\begin{split}
\ell_{\jb,n}(\tfrac{\kb}{(d+1)n}) = \langle \phi_{\kb},
\phi_{\jb}\rangle_{n}^* = \begin{cases}
 1, & \kb = \jb + (d+1)n \lb,\, \lb\in \ZZ_\HH^{d+1},\\
 0, & \text{otherwise},
\end{cases}
\end{split}
\end{align}
which proves \eqref{eq:sym-interp}.

To derive the compact formula for $\Phi_n^*(\tb)$ we start with the
following symmetry argument,
\begin{align*}
   \frac{1}{|\CG|}\sum_{\sigma\in \CG} \sum_{\kb\in \HH_n} \phi_{\kb}(\tb\sigma) & =
   \frac{1}{|\CG|}\sum_{\sigma\in \CG} \sum_{\kb\in \HH_n} \phi_{\kb\sigma^{-1}}(\tb) =
   \frac{1}{|\CG|}\sum_{\sigma\in \CG} \sum_{\kb\in \HH_n} \phi_{\kb\sigma}(\tb) \\
 &=  \frac{1}{|\CG|}\sum_{\sigma\in \CG} \sum_{\kb\in \HH_n\sigma} \phi_{\kb}(\tb)
 = \frac{1}{|\CG|}\sum_{\sigma\in \CG} \sum_{\kb\in (\HH_n\setminus\HH_n^\circ)\sigma} \phi_{\kb}(\tb)
 + \sum_{\kb\in \HH_n^\circ} \phi_{\kb}(\tb).
 \end{align*}
Using the fact that $\bigcup_{\sigma\in\CG} B^{i,j}_0\sigma= B^{i,j}$ and
\eqref{eq:Hnij}, we deduce that
\begin{align*}
\sum_{\sigma\in \CG} \sum_{\kb\in (\HH_n\setminus\HH_n^\circ)\sigma} \phi_{\kb}(\tb)
&  =  \sum_{0<i,j<i+j\le d+1} \sum_{\sigma\in \CG}\sum_{\kb\in \HH_{n,0}^{i,j}\sigma} \phi_{\kb}(\tb)
 =  \sum_{0<i,j<i+j\le d+1} \sum_{\kb\in \HH_{n}^{i,j}}\frac{(d+1)! i!j!}{(i+j)!} \phi_{\kb}(\tb)
\\
& =  (d+1)!  \sum_{0<i,j<i+j\le d+1} \sum_{\kb\in \HH_{n}^{i,j}}c_{\kb}^{(n)} \phi_{\kb}(\tb)
 = (d+1)! \sum_{\kb\in \HH_{n}^{*}\setminus\HH^{\circ}_n }c_{\kb}^{(n)} \phi_{\kb}(\tb),
\end{align*}
where the factor $\frac{(d+1)! i!j!}{(i+j)!} $ in the third summand comes from
the fact that, by Lemma \ref{lem:|BIJ|},  $\sum_{\sigma\in \CG}
\sum_{\kb\in \HH_{n,0}^{i,j}\sigma}$ summarizes over
$(d+1)!\times \frac{(d+1)!}{(i+j)!(d+1-i-j)!}$ terms whereas
$\sum_{\kb\in \HH_{n}^{i,j}}$ summarizes over $\frac{(d+1)!}{i!j!(d+1-i-j)!}$
terms. As a result, we obtain that
\begin{align*}
\frac{1}{|\CG|}&\sum_{\sigma\in \CG} \sum_{\kb\in \HH_n} \phi_{\kb}(\tb\sigma) =
\sum_{\kb\in \HH_{n}^{*}\setminus\HH^{\circ}_n }c_{\kb}^{(n)} \phi_{\kb}(\tb)
+ \sum_{\kb\in \HH_n^\circ} \phi_{\kb}(\tb)
 = \sum_{\kb\in \HH_n^*}c_{\kb}^{(n)} \phi_{\kb}(\tb).
\end{align*}
In other words, we have shown that
$$
  \Phi_n^*(\tb) = \frac{1}{|\CG|} \sum_{\sigma\in \CG} \wt D_n(\tb)
   \qquad \hbox{with} \qquad \wt D_n(\tb)= \sum_{\kb\in \HH_n} \phi_{\kb}(\tb\sigma).
$$

We now evaluate the partial sum $\wt D_n(\tb)$. Using Lemma
\ref{lm:DecompH} and the fact that $t_1 + \ldots + t_{d+1} =0$, we
obtain
\begin{align*}
\wt D_n(\tb) & =  \sum_{j=1}^{d+1} \sum_{\kb\in \HH^{\{j\}}_n}
\phi_{\kb}(\tb)
      = \sum_{j=1}^{d+1} \sum_{\kb \in \HH^{\{j\}}_n} \e^{\frac{2\pi i}{d+1}
          (  \sum_{1\le i<j} (k_i-k_j)t_i +  \sum_{j< l\le d+1} (k_l-k_j)t_l )} \\
   & = \sum_{j=1}^{d+1} \prod_{i=1}^{j-1}\sum_{\nu=1}^{n} \e^{2\pi i\nu t_i}
        \prod_{l=j+1}^{d+1}\sum_{\nu=0}^{n-1} \e^{2\pi i\nu t_l}
= \sum_{j=1}^{d+1} \prod_{i=1}^{j-1}\frac{\sin \pi n t_i}{\sin \pi
t_i} \e^{\pi i (n+1) t_i}
       \prod_{l=j+1}^{d+1} \frac{\sin \pi n t_l}{\sin \pi t_l} \e^{\pi i (n-1) t_l}\\
        & 
= \prod_{i=1}^{d+1}\frac{\sin \pi n t_i}{\sin \pi t_i}
   \sum_{j=1}^{d+1} \e^{\pi i (t_j+2\sum_{1\le i < j}t_i) }
        \frac{\sin \pi t_j}{\sin \pi n t_j}\e^{-\pi i n t_j }.
\end{align*}
Using symmetry and  $t_1 + \ldots + t_{d+1} =0$ again, we further derive that
\begin{align*}
   \frac{1}{|\CG|}\sum_{\sigma \in \CG} \wt D_n(\tb\sigma)
   = \frac{1}{|\CG|} \prod_{i=1}^{d+1}\frac{\sin \pi n t_i}{\sin \pi t_i}   \sum_{j=1}^{d+1}
    \frac{\sin \pi t_j}{\sin \pi n t_j}\e^{-\pi i n t_j}
 \sum_{I\subseteq\NN_{d+1}^{\{j\}} } |\CG_{I}|\times |\CG_{\NN_{d+1}^{\{j\}}\setminus I }|\,
  \e^{\pi i ( t_j+2\sum_{i\in I} t_i)}.
\end{align*}
Since $\HH_n^*$ is symmetric under the mapping $\tau_i: \tb \mapsto
  (t_1,\ldots, t_{i-1},-t_i,t_{i+1}, \ldots, t_{d+1})$ and the definition of $c_\kb^{(n)}$
implies that $c_{\kb \tau_i}^{(n)} = c_\kb^{(n)}$, we see that the function
$\sum_{\kb \in \HH_n^*} c_\kb^{(n)}\phi_\kb(\tb)$ is invariant under the action
of $\tau_i$. As a consequence, we this function must be real valued. Hence,
taking the real part in the above expression, we conclude that
\begin{align*}
\sum_{\kb \in \HH_n^*} c_\kb^{(n)}\phi_\kb(\tb) = \sum_{j=1}^{d+1}
\bigg(\prod_{1\le i\le d+1\atop i\neq j }
    \frac{\sin \pi n t_i}{\sin \pi t_i}  \bigg)      \cos \pi n t_j
   \sum_{I\subseteq\NN_{d+1}^{\{j\}} } \frac{|I|!(d-|I|)!}{(d+1)!}
\, \cos \pi \big( t_j+2\sum_{i\in I} t_i\big),
\end{align*}
which completes the proof.
\end{proof}

We note that the condition $t_1 + \ldots + t_{d+1} =0$ implies various
relations between trigonometric functions $\sin \pi t_i$ and $\cos \pi t_i$.
For example, in the case of $d =2$, we have that
\begin{align*}
  \sin 2 t_1+ \sin 2 t_2 + \sin 2 t_3 & = -4 \sin  t_1\sin  t_2\sin  t_3 \\
  \cos 2 t_1+ \cos 2 t_2 + \cos 2 t_3 & = 4 \cos t_1\cos  t_2\cos  t_3 -1.
\end{align*}
Using these relations, it is possible to rewrite the compact formula in
\eqref{eq:Phi-n*} in other forms. For instance, in the case of $d =2$, the
compact formula \eqref{eq:Phi-n*} becomes
\begin{align*}
\Phi^{*}_n(\tb) = \frac{1}{3n^d}\frac{\sin \pi n t_1}{\sin \pi t_1 }
  \frac{\sin \pi n t_2}{\sin \pi t_2 } \cos \pi n t_3
   \left(\tfrac{2}{3}\cos \pi t_3 + \tfrac{1}{3}\cos\pi(t_1-t_2) \right)
\\+\frac{1}{3 n^d}\frac{\sin \pi n t_2}{\sin \pi t_2 } \frac{\sin \pi n t_3}
 {\sin \pi t_3 } \cos \pi n t_1
\left(\tfrac{2}{3} \cos \pi t_1 + \tfrac{1}{3} \cos\pi(t_2-t_3) \right)\\
+\frac{1}{3 n^d}\frac{\sin \pi n t_3}{\sin \pi t_3 } \frac{\sin \pi n t_1}{\sin \pi t_1 }
  \cos \pi n t_2 \left(\tfrac{2}{3} \cos \pi t_2 + \tfrac{1}{3}\cos\pi(t_3-t_1) \right),
\end{align*}
whereas a different formula is derived in \cite{LSX}, given in terms
of $\Theta_n(\tb)$.

Using the compact formula of $\Phi_n^*$, we can estimate the operator norm
of $I_n^*$.

\begin{thm} \label{LebesgueH}
Let $\|I_n^*\|_\infty$ denote the operator norm of $I_n^*:
C(\overline{\Omega}_H)
 \mapsto C(\overline{\Omega}_H)$. Then there is a constant $c$, independent of
$n$, such that
$$
       \|I_n^*\|_\infty  \le c (\log n)^d.
$$
\end{thm}

\begin{proof}
A standard procedure shows that
$$
   \|I_n^*\|_\infty = \max_{\tb \in \overline{\Omega}_H}  \sum_{\kb \in \HH_n^*}
         \left|\Phi_n^*(\tb-\tfrac{\kb}{(d+1)n})\right|.
$$
Using the compact formula of $\Phi_n^*$ in Theorem \ref{interp-H},
it is easy to see that it suffices to prove that
$$
I_{j}:=  \frac{1}{(d+1)n^{d}}  \max_{\tb\in \overline{\Omega}_H}  \sum_{\kb
\in \HH_n^*}
         \bigg|\prod_{1\le i\le d+1\atop i\neq j}\frac{\sin \pi n (t_i- \frac{k_i}{(d+1)n} )}{\sin \pi (t_i- \frac{k_i}{(d+1)n} ) }
         \bigg| \le c (\log n)^d, \qquad  1\le j\le d+1,\ n \ge 0.
$$
For $\tb \in \RR_H^{d+1}$, let $\tb^{(j)}:= (t_1,\ldots,t_{j-1},t_{j+1}, \ldots, t_{d+1})
\in \RR^d$. It is easy to see that $\{\tb^{(j)}: \tb \in \overline \Omega_H\}
\subset [-1,1]^d$. Hence, enlarging the domain $\HH_n^*$ to $\{\kb \in \ZZ_H^{d+1}:
-(d+1)n \le k_i \le (d+1)n, \,k_i \equiv 0 \pmod{d+1}, \,1 \le i \le d\}$, we see that
\begin{align*}
   I_j & \le \frac{1}{(d+1)n^{d}}  \max_{t \in {[-1,1]^d}}  \prod_{l=1}^{d} \sum_{k_l=-n}^{n}
      \left|\frac{\sin \pi n (t_l- \frac{k_l}{n} )}{\sin \pi (t_l- \frac{k_l}{n} ) } \right|\\
   &   \le    \frac1{d+1}  \max_{t \in {[-1,1]}} \left (\frac{1}{n}
        \sum_{k = -n}^{n} \left | \frac{\sin n \pi (t -\tfrac{k}{n})}
          {\sin\pi (t -\tfrac{k}{n})} \right|  \right)^d   \le c (\log n)^d,
\end{align*}
where the last step follows from the usual estimate of one
variable (cf. \cite[Vol. II, p. 19]{Z}).
\end{proof}

We expect that the estimate is sharp, that is, $\|I_n^*\| \ge c (\log n)^d$, but
do not have a proof at this point.


\section{Discrete Fourier analysis on the simplex}
\setcounter{equation}{0}

The fundamental domain of the lattice $A_d$ is the union of the
images of the fundamental simplex $\triangle_H$ under the group
$\CG$, as shown in \eqref{eq:SimpDecom}. Hence, if we consider
invariant functions under the group $\CG$, then the discrete Fourier
analysis on the fundamental domain in the previous section can be
carried over to the analysis on the simplex $\triangle_H$, which is
developed below.

\subsection{Generalized sine and cosine functions}
In the case of one-variable, the invariant and anti-invariant sums
of the exponential functions are cosine and sine functions. We now
consider their analogous in our setting. It turns out that these
functions have already appeared in the literature, as mentioned in
the introduction.  It should be pointed out, however, that our study
is on the discrete Fourier analysis, which has little overlap with
the previous study in the literature. Most of the overlap will
appear in the next section, when these trigonometric functions are
transformed to Chebyshev polynomials.

Recall that the reflection group $\A_d$ is the permutation group $\CG$.
Denote the identity element in $\CG$ by $1$. It is easy to see that
\begin{align}
  \label{eq:sigma}
  \sigma_{ij}^2 =1, \qquad \sigma_{ij}\sigma_{jk}\sigma_{ij} = \sigma_{ik},
    \quad i,j,k \in \NN_{d+1}.
\end{align}
For $\sigma \in \G$, let $|\sigma|$ denote the number of
inversions in $\sigma$. The group $\G$ is naturally divided into two
parts, $\G^+: = \left\{\sigma\in \mathcal{G}:\  |\sigma| \equiv 0
\pmod{2} \right\}$ of elements with even inversions, and $\G^{-}: =
\left\{\sigma\in \mathcal{G}:\  |\sigma| \equiv 1 \pmod{2} \right\}$
of elements with odd inversions.
The action of $\sigma \in \G$ on the function $f: \RR_H^{d+1}
\mapsto \RR$ is defined by $\sigma f (\tb):= f(\tb \sigma)$. A
function $f$ in homogeneous coordinates is called {\it invariant}
under $\G$ if $\sigma f = f$  for all $\sigma \in \G$, and it is
called {\it anti-invariant} under $\G$ if $\sigma f = \rho(\sigma)f$
with $\rho(\sigma) =1$ if $\sigma \in \G^+$ and $\rho(\sigma) = -1$
if $\sigma \in \G^-$.
The following proposition follows immediately from the definition.

\begin{prop}
Define two operator $\mathcal{P}^{+}$ and $\mathcal{P}^{-}$ acting
on $f(\tb)$ by
\begin{align} \label{eq:Ppm}
\mathcal{P}^{\pm}f(\tb) := \frac1{(d+1)!} \bigg[ \sum_{\sigma\in
\mathcal{G}^{+}}
   f(\tb \sigma)\pm \sum_{\sigma\in \mathcal{G}^{-}} f(\tb \sigma) \bigg].
\end{align}
Then the operators $\mathcal{P}^{+}$ and $\mathcal{P}^{-}$ are
projections from the class of $H$-periodic functions onto the class
of invariant, and respectively anti-invariant functions.
\end{prop}


Recall that $\phi_\kb(\tb) = \e^{\frac{2\pi i}{d+1} \kb \cdot \tb}$.
Applying the operators $\CP^\pm$ to these exponential functions
gives the basic invariant and anti-invariant functions.

\begin{defn} 
For $\kb \in \HH$ define
\begin{align*}
&\TC_{\kb}(\tb):=\mathcal{P}^{+}\phi_{\kb}(\tb) = \frac1{(d+1)! }
\bigg[ \sum_{\sigma\in \mathcal{G}^{+}}\phi_{\kb}(\tb \sigma)
+ \sum_{\sigma\in \mathcal{G}^{-}} \phi_{\kb}(\tb \sigma) \bigg],\\
&\TS_{\kb}(\tb):= \mathcal{P}^{-}\phi_{\kb}(\tb) = \frac{1}{(d+1)!}
 \bigg[ \sum_{\sigma\in \mathcal{G}^{+}}\phi_{\kb}(\tb \sigma)
- \sum_{\sigma\in \mathcal{G}^{-}} \phi_{\kb}(\tb \sigma) \bigg],
\end{align*}
and call them \emph{generalized cosine} and \emph{generalized sine},
respectively.
\end{defn}

By definition, $\TC_{\kb}$ is invariant and $\TS_{\kb}$ is anti-invariant.
Because of the symmetry, we only need to consider them on the
fundamental simplex $\triangle_H$ defined in \eqref{triangle} or any other
simplex $\triangle_H \sigma$, $\sigma \in \CG$, that makes up
$\overline \Omega_H$. We shall work with $\triangle_H$ below and
recall that
\begin{align*}
   \triangle_H = \left\{ \tb \in \RR_H^{d+1}:
               0\leq t_i-t_j \le 1,\, 1\le i\le j\le d+1\right\}.
\end{align*}
In the case of $d =2$ and $d=3$, it is an equilateral triangle and a regular
tetrahedron, respectively; these regions are depicted in  Figure 4.1, in which
the corners are given in homogeneous coordinates.

\begin{figure}[htb]
\centering
\hfill\includegraphics[width=0.4\textwidth]{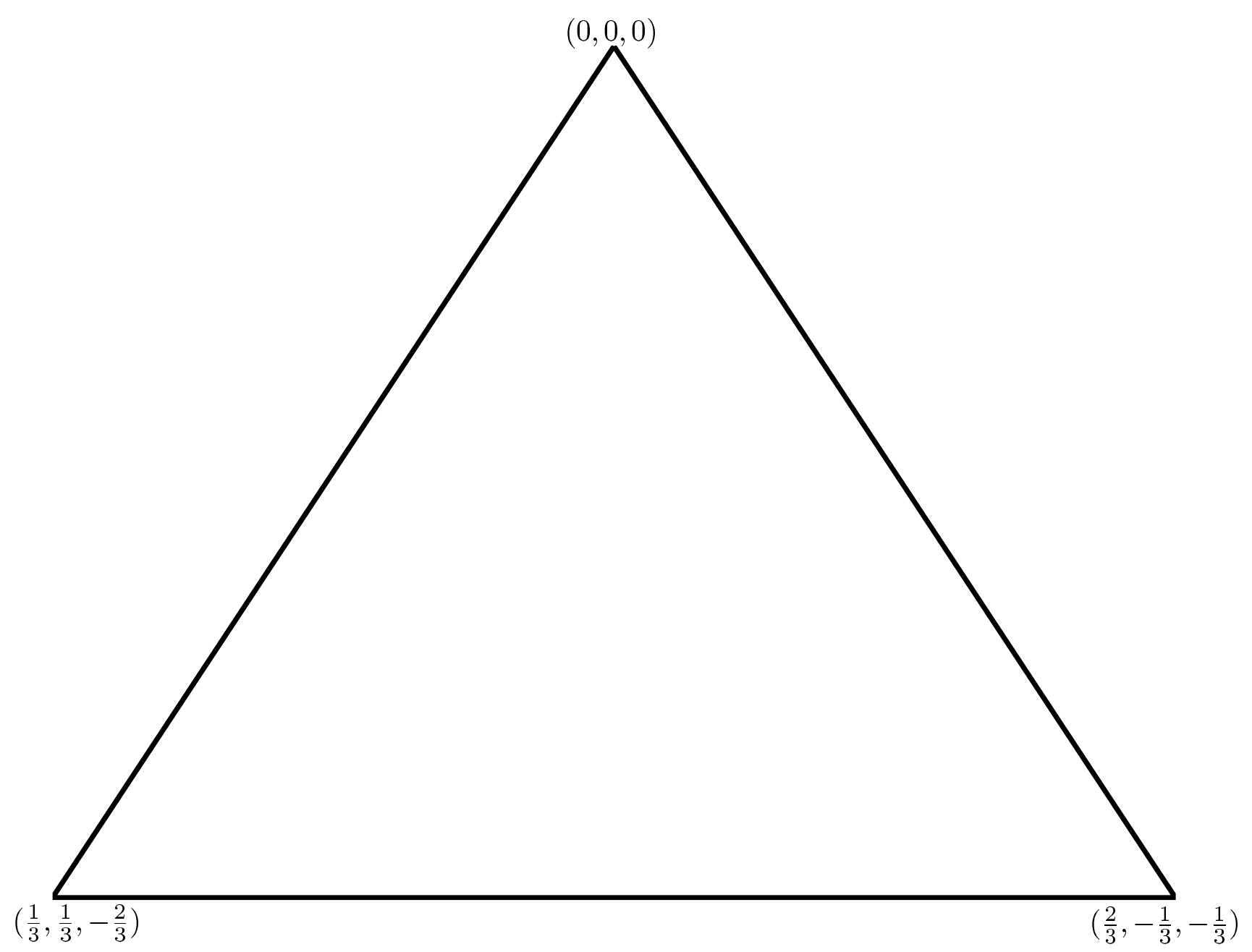}\hfill\includegraphics[width=0.55\textwidth]{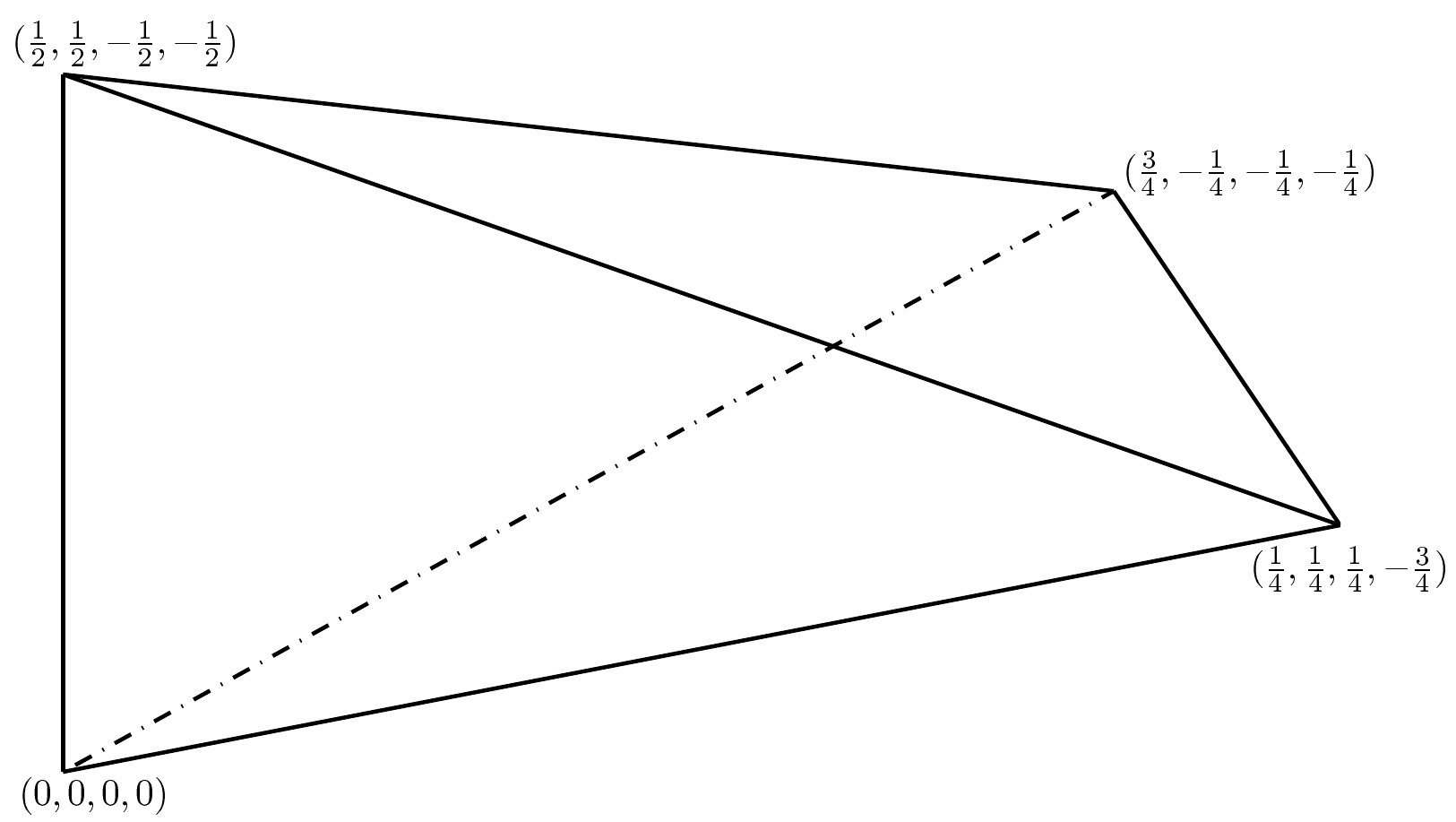}
\hspace*{\fill}
\caption{Reference triangle and tetrahedron.}
\end{figure}

By definition, $\phi_\kb$ and $\phi_{\kb \sigma}$ with $\sigma \in \CG$ lead to
the same $\TC_\kb$. In fact, we have
\begin{align}
\label{eq:TCs}
     \TC_{\kb \sigma}(\tb) = \TC_{\kb}(\tb\sigma) = \TC_{\kb}(\tb)
\qquad \hbox{for $\tb\in\triangle_H$ and $\sigma\in\G$}.
\end{align}
Thus, when working with $\TC_\kb$, we can restrict $\kb$ to the index set
\begin{align*}
\Lambda:= \left\{ \kb\in \HH: k_1  \geq k_2 \ge \cdots \ge k_{d+1} \right\}.
\end{align*}
As for $\TS_{\kb}$, it is easy to see that $\TS_{\kb\sigma}(\tb)=
\TS_{\kb}(\tb\sigma) = \TS_{\kb}(\tb)$ for $\sigma\in\mathcal{G}^{+}$
and  $\TS_{\kb\sigma}(\tb)= \TS_{\kb}(\tb\sigma) = -\TS_{\kb}(\tb) $ for
$\sigma\in\mathcal{G}^{-}$. In particular, $\TS_{\kb}(\tb)=0$
whenever two or more components of $\kb$ are equal. Thus, when working
with $\TS_{\kb}$, we only need to consider $\kb\in \Lambda^{\circ}$, where
\begin{align*}
\Lambda^{\circ}:= \left\{ \kb\in \HH: k_1 > k_2 > \cdots > k_{d+1}\right\},
\end{align*}
which is the set of the interior points of $\Lambda$. To describe the points
on the boundary of $\Lambda$, we need to consider the compositions of
the integer $d+1$. A composition $p$ of $d+1$ is a decomposition of $d+1$,
such that
$$
p = (p_1,\ldots, p_\ell) \in \ZZ^\ell, \quad p_i > 0, \quad1 \le i \le \ell \le d+1,
    \quad |p| := p_1+\ldots + p_{\ell} = d+1,
$$
where $\ell = \ell(p)$ is the length of the composition. Notice that the order
of $p_i$ matters, different orderings are deemed to be different compositions,
which is the difference between a composition and a partition. We denote
the collection of compositions of $d+1$ by $\mathcal{C}_{d+1}$; that is,
\begin{align*}
   \mathcal{C}_{d+1}:= \left\{ p \in \ZZ^{\ell(p)}:
         1 \le p_i \le |p| = d+1 \text{ for } 1\le i\le \ell(p) \le d+1 \right\}.
 \end{align*}
For $p \in \mathcal{C}_{d+1}$ we further define
\begin{align*}
\Lambda^{p} & := \left\{ \kb =\left ( \left\{k_1 \right\}^{p_1},  \left\{k_2 \right\}^{p_2},
\cdots, \left\{k_{\ell} \right\}^{p_{\ell}} \right) \in \HH: k_i>k_j, 1\le i < j \le \ell \right\}.
\end{align*}
Then evidently $\Lambda^\circ =  \Lambda^{\left\{1\right\}^{d+1}}$ and
$\Lambda = \bigcup_{p\in \mathcal{C}_{d+1}} \Lambda^p$.

Recall that $\kb\mathcal{G}= \left\{\kb\sigma :\ \sigma\in
\mathcal{G}\right\}$ is the orbit of $\kb $ under $\mathcal{G}$. The
definition of $\Lambda$ implies that, for $\kb, \jb \in \Lambda$,
$\kb\mathcal{G}\cap \jb\mathcal{G} = \emptyset$ whenever $\kb\neq
\jb$. It follows that
\begin{align*}
  |\kb\mathcal{G}| = \binom{d+1}{p}:=  \frac{(d+1)!}{p_1!p_2!\cdots p_{\ell}!}
     \qquad \text{ if }\kb \in  \Lambda^{p}.
\end{align*}
Since $\G/\G_\kb$ is isomorphic to $\kb \G$, where $\G_\kb$ is the
stabilizer of $\kb$, we also have
\begin{align} \label{eq:TC-redef}
 \TC_{\kb} (\tb) =  \frac{1}{|\mathcal{G}|}
  \sum_{\jb\in  \kb \mathcal{G}} \sum_{\sigma \in \G_\jb} \phi_{\jb \sigma}(\tb)
 =  \frac{1}{|\kb\mathcal{G}|}
  \sum_{\jb\in  \kb \mathcal{G}} \phi_{\jb}(\tb).
\end{align}

The length $\ell = \ell(p)$ of $p \in \mathcal{C}_{d+1}$ determines how many
indices in $\kb \in \Lambda^p$ are repeated, which determines the dimension
of the boundary elements of $\Lambda$. For instance, if $d=2$ then $\CG_3=
\{(1,1,1),(1,2),(2,1), (3)\}$ and
\begin{align*}
|\kb\mathcal{G}| =
  \begin{cases}
        6, & \kb\in \Lambda^{\circ}:= \Lambda^{1,1,1},\\
        3, & \kb\in \Lambda^{e}:= \Lambda^{1,2}\cup \Lambda^{2,1},\\
        1, & \kb\in \Lambda^{v}:= \Lambda^{3},\\
  \end{cases}
\end{align*}
where $\Lambda^e$ and $\Lambda^v$ consist of points in $\Lambda$
that are on the edges and vertices of $\Lambda$, respectively. If
$d=3$ then $\CG= \{(1,1,1,1),(1,1,2),(1,2,1),(2,1,1), (1,3),(3,1),
(2,2),(4)\}$ and
\begin{align*}
|\kb\mathcal{G}| =
  \begin{cases}
        24, & \kb\in \Lambda^{\circ}:= \Lambda^{1,1,1,1},\\
        12, & \kb\in \Lambda^{f}:= \Lambda^{1,1,2}\cup \Lambda^{1,2,1}\cup
           \Lambda^{2,1,1} ,\\
        6, & \kb\in \Lambda^{e,1}:= \Lambda^{2,2},\\
        4, & \kb\in \Lambda^{e,2}:= \Lambda^{1,3}\cup \Lambda^{3,1},\\
        1, & \kb=0\in \Lambda^{v}:= \Lambda^{4},
  \end{cases}
\end{align*}
where $\Lambda^f$, $\Lambda^e$ and $\Lambda^v$ consist of points in
$\Lambda$ that are on the faces, edges and vertices of $\Lambda$,
respectively.

We define an inner product on $\triangle_H$ by
\begin{align*}
\langle f, g \rangle_{\triangle_H}:= \frac{1}{|\triangle_H|}
\int_{\triangle_H} f(\tb) \overline{g({\tb})} d\tb
   = \frac{(d+1)!}{\sqrt{d+1}} \int_{\triangle_H} f(\tb) \overline{g({\tb})} d\tb.
\end{align*}
If $f\bar g$ is invariant under $\mathcal{G}$, then it follows
immediately that $\langle f, g \rangle = \langle f, g \rangle_{\triangle_H}$,
where $\langle \cdot,\cdot\rangle$ is the inner product defined in
\eqref{eq:ip-H} over $\Omega_H$.  The generalized cosine and sine
functions are orthogonal with respect to this inner product.

\begin{thm}\label{prop:trig-ortho}
For $\kb,\jb\in \Lambda$,
 \begin{align}
 \label{eq:orth-TC}
   \langle \TC_{\kb}, \TC_{\jb}\rangle_{\triangle_H}  =
     \frac{\delta_{\kb,\lb}}{|\kb\mathcal{G}|}
   =  \frac{p_1!p_2!\cdots p_{\ell}!}{(d+1)!}\delta_{\kb,\lb}  , \qquad
   \kb \in \Lambda^{p},
\end{align}
and for $\kb, \jb \in \Lambda^{\circ}$,
 \begin{align}
 \label{eq:orth-TS}
 \langle \TS_{\kb}, \TS_{\jb}\rangle_{\triangle_H}  = \frac1{(d+1)!}\delta_{\kb,\jb}.
 \end{align}
\end{thm}

\begin{proof}
Both of these relations follow from the identity $\langle f, g
\rangle = \langle f, g \rangle_{\triangle_H}$ for invariant
functions. For \eqref{eq:orth-TC}, the invariance is evident and we
only have to use the orthogonality of $\phi_\kb$ in
\eqref{eq:orth-H} and \eqref{eq:TC-redef}. For \eqref{eq:orth-TS},
we use the fact that $\TS_{\kb}(\tb)\overline{\TS_{\jb}(\tb)}$ is
invariant under $\G$ and the orthogonality of $\phi_\kb$ on
$\Omega_H$.
\end{proof}

The definition and orthogonality of these trigonometric functions have
appeared in the literature, we refer to \cite{Be} and its extensive references.
However, the study in the literature is more on the side of algebraic
polynomials, as will be discussed in Section 5 below.


\subsection{Discrete inner product on the simplex}

By Theorem \ref{ipdH} and \eqref{eq:orth-H}, $\{\phi_\kb:\kb \in
\HH_n\}$ is an orthonormal set with respect to the symmetric inner
product $\langle \cdot, \cdot\rangle_n^*$. Using the symmetry and
the invariance of $\TC_{\kb}$ and $\TS_{\kb}$ under $\G$, we can
deduce a discrete orthogonality for the generalized cosine and sine
functions. We define
\begin{align} \label{eq:Lambda_n}
\Lambda_n : = \HH_n^* \cap \Lambda =
 \left\{\kb\in \HH:  k_{d+1}\leq k_d\leq \ldots \leq k_1 \le  k_{d+1}+(d+1)n  \right \}.
\end{align}
 For $n=4$, the point sets for $d =2$ and $d=3$ are depicted in Figure 4.1.

\begin{figure}[htb]
\centering
\hfill\includegraphics[width=0.4\textwidth]{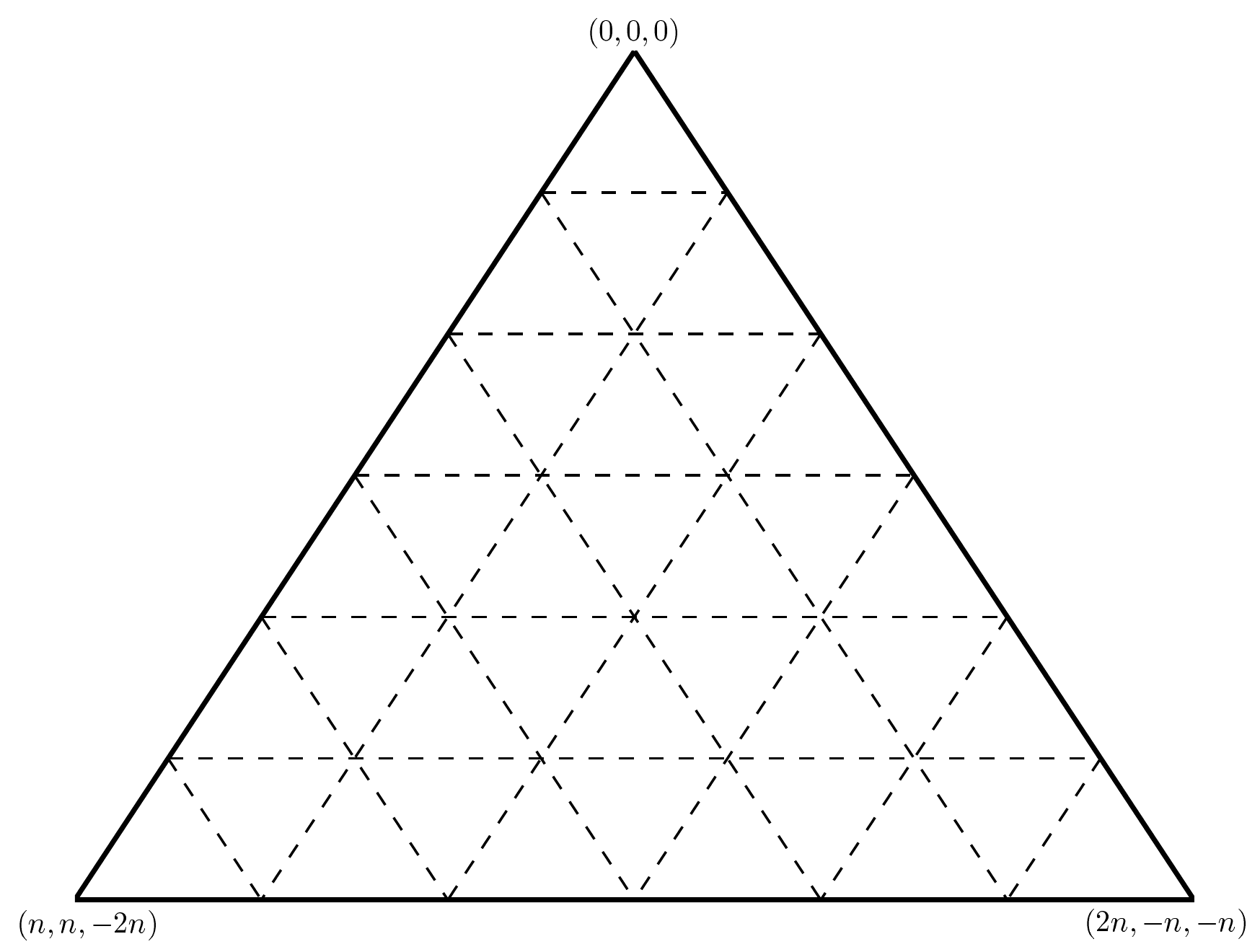}\hfill
\includegraphics[width=0.55\textwidth]{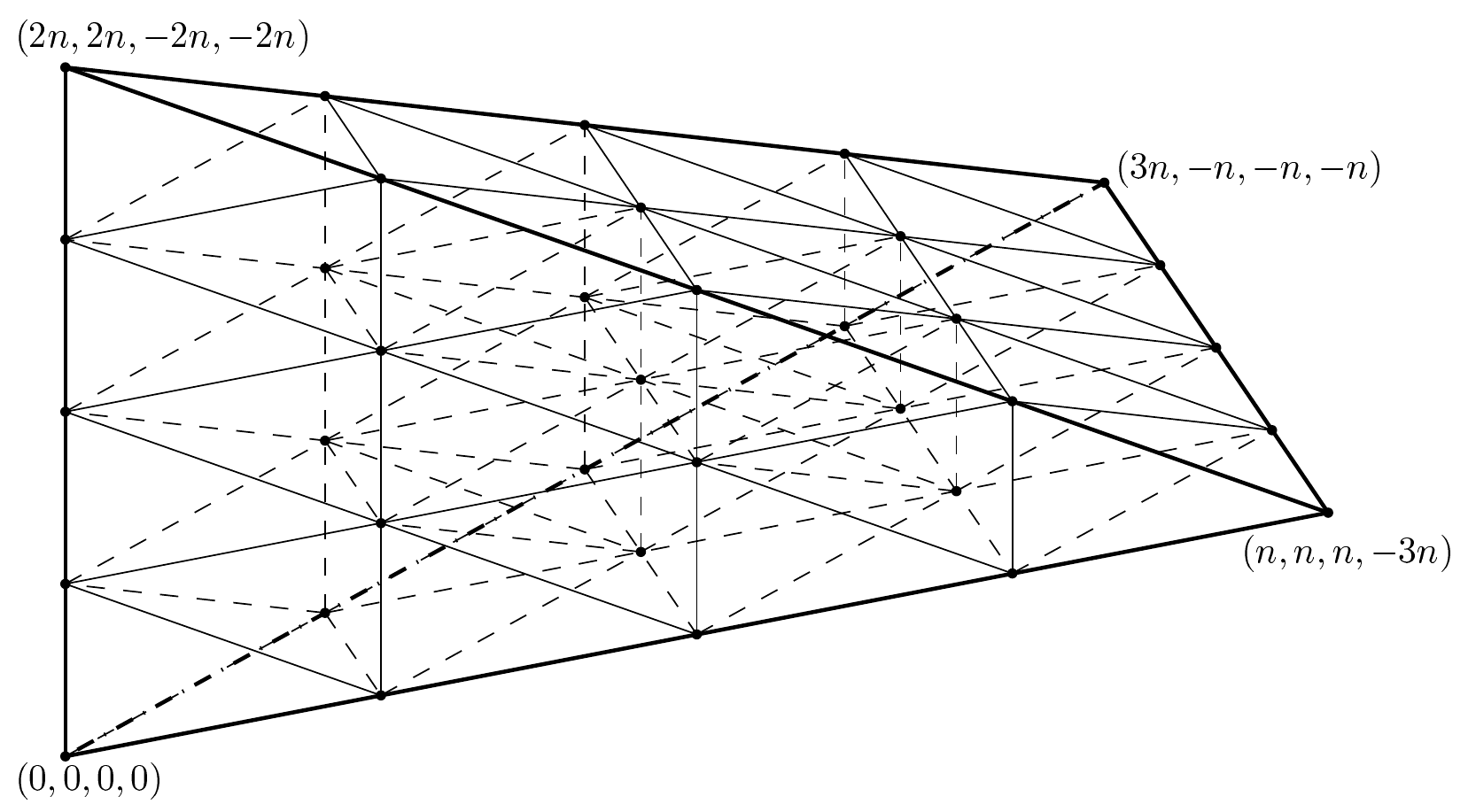}\hspace*{\fill}
\caption{$\Lambda_n$ for $d=2$ and $d =3$}
\end{figure}

By the definition of $\HH_n^*$,  the set $\{\tfrac{\kb}{(d+1)n}: \kb
\in \Lambda_n\}$ contains points inside $\triangle_H$. We will also
need to understand points on the boundary of $\triangle_H$, which
can come from two types of points in $\HH_n^*$. One part of the
boundary points in $\triangle_H$ are also the boundary points of
$\HH_n^*$, whereas another part of the boundary points in
$\triangle_H$ are points inside $\HH_n^\circ$ but on the faces of
$\triangle_H$. Accordingly, for $p\in \mathcal{C}_{d+1}$, we define
\begin{align*}
 \Lambda_n^{\circ,p} := \left\{( \{k_1\}^{p_1}, \{k_2\}^{p_2},\cdots, \{k_{\ell}\}^{p_{\ell}})
   \in \HH:\ k_{\ell}+(d+1)n>  k_1>k_2>\cdots>k_{\ell} \right\},\\
 \Lambda_n^{\partial,p} := \left\{( \{k_1\}^{p_1}, \{k_2\}^{p_2},\cdots,
  \{k_{\ell}\}^{p_{\ell}}) \in \HH: k_{\ell}+(d+1)n=  k_1>k_2>\cdots>k_{\ell} \right\}.
\end{align*}
Evidently, $\Lambda_n^{\circ}= \Lambda^{\circ,\{1\}^{d+1}}_n$ and
$$
\Lambda_n = \bigcup_{p\in \mathcal{C}_{d+1}} \left( \Lambda_{n}^{\circ,p}\cup
  \Lambda_{n}^{\partial,p} \right).
$$
If $p \ne \{1\}^{d+1}$, then $\Lambda_n^{\circ,p}$ is a $(\ell-1)$-face of
$\Lambda_n$ which, however, is a subset of $\HH_n^\circ$, whereas
$\Lambda_n^{\partial,p}$ is a $(\ell-2)$-face of $\Lambda_n$ which is also
a face of $\HH_n^*$. More precisely, considering the orbits of the points in
$\Lambda_n$, we see that
\begin{align} \label{H=Lambda}
  \HH_n^\circ = \bigcup_{p \in \mathcal{C}_{d+1}}  \Lambda_n^{^\circ,p} \G
    \qquad \hbox{and} \qquad
  \HH_n^*\setminus \HH_n^\circ = \bigcup_{p \in \mathcal{C}_{d+1}}
    \Lambda_n^{\partial,p} \G.
\end{align}

For $d=2$, the sets of interior points, edge points and vertices of $\triangle_H$ are
given, in this order, explicitly by
\begin{align*}
& \Lambda^{\circ,1,1,1}_n=\left\{\kb\in \HH:\  k_3+3n>k_1>k_2>k_3\right\},\\
& \Lambda^{\circ,1,2}_n \cup \Lambda^{\circ,2,1}_n \cup \Lambda^{\partial,1,1,1}_n
 =\left\{  (2k, -k,-k), (k, k, -2k), (3n-k, 2k-3n, -k):\   0<k<n\right\},\\
& \Lambda^{\circ,3}_n\cup \Lambda^{\partial,1,2}_n
\cup\Lambda^{\partial,2,1}_n =\left\{(0,0,0), (2n,-n,-n), (n,n,-2n)
\right\}.
\end{align*}
In this case, the geometry in Figure 4.3 shows clearly that two of the three edges,
$\Lambda^{\circ,1,2}_n \cup \Lambda^{\circ,2,1}_n$, of the triangle $\triangle_H$
are not edges of the hexagon, and the vertex $\Lambda^{\circ,3}_n =\{(0,0,0)\}$
is not an vertex of the hexagon.

\begin{figure}[htb]
\centering
 \includegraphics[width=0.5\textwidth]{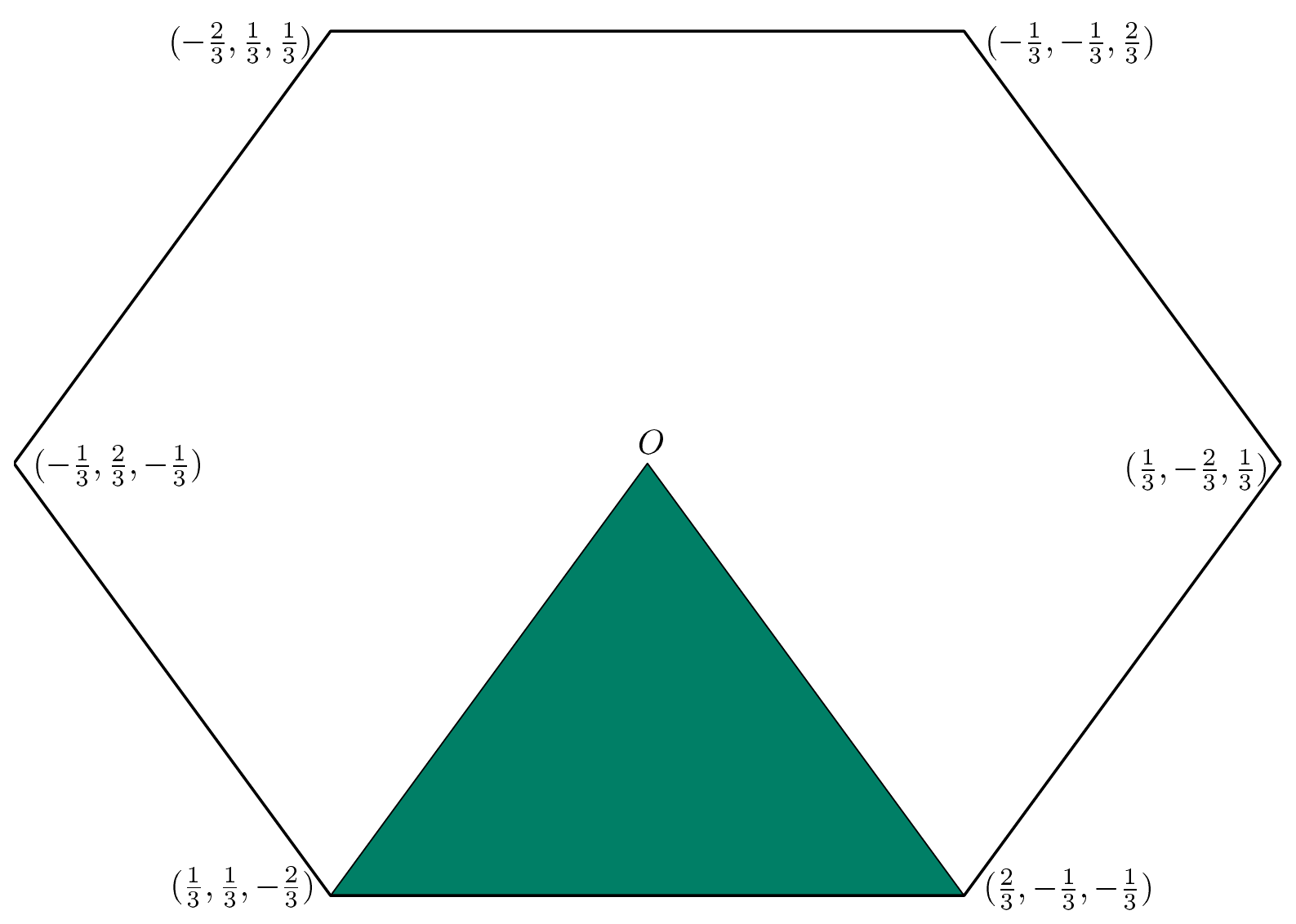}
\caption{Reference triangle $\triangle_H$ inside the hexagon}
\end{figure}

For $d=3$, the boundary sets are given in detail in \cite{LX08} and they are
denoted by $\Lambda_n^{\circ}$, $\Lambda_n^{f}$, $\Lambda_n^{e,1}$,
$\Lambda_n^{e,2}$ and  $\Lambda_n^{v}$, which correspond to points
of interior, faces, two type of edges, and vertices. In terms of the sets
$\Lambda_n^{\circ,p}$ and $\Lambda_n^{\partial,p}$, they are expressed by
\begin{align*}
&  \Lambda_n^{\circ} = \Lambda^{\circ,1,1,1,1}_n,\\
& \Lambda_n^{f}   = \Lambda^{\circ,1,1,2}_n\cup \Lambda^{\circ,1,2,1}_n \cup \Lambda^{\circ,2,1,1}_n \cup \Lambda^{\partial,1,1,1,1}_n,\\
& \Lambda_n^{e,1} = \Lambda^{\circ,2,2}_n  \cup \Lambda^{\partial,1,2,1}_n, \\
& \Lambda_n^{e,2} = \Lambda^{\circ,1,3}_n\cup \Lambda^{\circ,3,1}_n \cup \Lambda^{\partial,1,1,2}_n \cup \Lambda^{\partial,2,1,1}_n,\\
& \Lambda_n^{v}  = \Lambda^{\circ,4}_n \cup \Lambda^{\partial,1,3}_n \cup \Lambda^{\partial,2,2}_n  \cup \Lambda^{\partial,3,1}_n.
\end{align*}

We now define a discrete inner product $\langle \cdot,\cdot \rangle_{\triangle,n}$ by
\begin{align*}
    \langle f,g \rangle_{\triangle,n} = \frac{1}{(d+1)n^{d}} \sum_{\jb\in \Lambda_n}
    \lambda_{\jb}^{(n)} f(\tfrac{\jb}{(d+1)n}) \overline{g(\tfrac{\jb}{(d+1)n})}.
\end{align*}
where, with $c^{(n)}_{\jb}$ defined in Definition \eqref{defn:Sym-ipd},
\begin{align} \label{lambdaj}
  \lambda_{ \jb}^{(n)} :=  c^{(n)}_{\jb} \binom{d+1}{p}= \begin{cases}
     \dfrac{(d+1)!}{p_1!p_2!\cdots p_{\ell-1}!p_{\ell}!},
       \qquad \jb \in \Lambda^{\circ,p}_n,\\[1em]
     \dfrac{(d+1)!}{(p_{\ell}+p_1)!p_2!\cdots p_{\ell-1}!},
       \qquad \jb \in \Lambda^{\partial, p}_n.
\end{cases}
\end{align}
Let us verify the second equal sign in \eqref{lambdaj}. The case
$\jb \in\Lambda^{\circ,p}_n$ is easy, since then $\jb \in
\HH_n^\circ$ so that $c_\jb^{(n)} =1$ and $\lambda_\jb^{(n)} =
\binom{d+1}{p}$. In the case $\jb \in \Lambda^{\partial, p}_n$, we
have $j_1 - j_\ell = (d+1)n$ which implies $j_i - j_l = (d+1)n$ for
$i \in I_{p_1}$ and $l \in I_{p_\ell}$, where $I_{p_1} = \{1,\ldots,
p_1\}$ and $I_{p_\ell} =
 \{d+1-p_\ell +1, \ldots, d+1\}$. Consequently, $\frac{\jb}{(d+1)n} \in
 B^{p_1,p_\ell}$ by the definition of $B^{i,j}$, which implies that
 $\jb \in \HH^{p_1,p_\ell}$ so that $c_\ib^{(n)} = 1/\binom{p_1+p_\ell}{p_1}$
and \eqref{lambdaj}.

In the case of $d =2$ and $d=3$, the values of $\lambda_\jb^{(n)}$ are given by
\begin{align*}
\lambda_{ \jb}^{(n)} =
\begin{cases}
          6, & \jb\in\Lambda^{\circ}_n,\\
          3, & \jb\in\Lambda^{e}_n,\\
          1, &  \jb\in\Lambda^{v}_n,
\end{cases} \quad \text{ if } d=2 \qquad \text{ and } \qquad
\lambda_{ \jb}^{(n)} =
\begin{cases}
          24, & \jb\in\Lambda^{\circ}_n,\\
          12, & \jb\in\Lambda^{f}_n,\\
           6,  &  \jb\in\Lambda^{e,1}_n,\\
           4,  &  \jb\in\Lambda^{e,2}_n,\\
           1,  &  \jb\in\Lambda^{v}_n, \end{cases}
\quad \text{ if } d=3.
\end{align*}

We denote by $\CTC_n$ and $\CTS_n$ the spaces of the trigonometric
polynomials
\begin{align*}
\CTC_n:=\mathrm{span}\left\{\TC_{\kb}:\ \kb\in \Lambda_n  \right\},
\qquad \CTS_n:=\mathrm{span}\left\{\TS_{\kb}:\ \kb\in \Lambda_n^{\circ} \right\},
\end{align*}
respectively. Since $\Lambda_n$ contains integer points in a regular
simplex, it is easy to see that
\begin{align} \label{dimCT}
  \dim \CTC_n = |\Lambda_n| = \binom{n+d}{d} \qquad \hbox{and}\qquad
       \dim \CTS_n =  |\Lambda_n^\circ| = \binom{n -1}{d}.
\end{align}

\begin{thm} \label{thm:4.4}
For $f \bar{g} \in \CTC_{2n-1}$,
\begin{align}  \label{eq:orth-Tn}
   \la f, g \ra_{\triangle_H} =  \la f, g \ra_{\triangle,n}.
\end{align}
Moreover, the following cubature formula is exact for all $f\in
\CTC_{2n-1}$,
\begin{align} \label{eq:cubature-T}
   \frac{1}{|\triangle_H|} \int_{\triangle_H}f(\tb) d\tb =
       \frac{1}{(d+1)n^{d}} \sum_{\jb\in \Lambda_n} \lambda_{\jb}^{(n)} f(\tfrac{\jb}{(d+1)n}).
\end{align}
In particular,
\begin{align}  \label{eq:orth-TCn}
   \langle \TC_{\kb},\TC_{\jb} \rangle_{\triangle,n} =
    \frac{ \delta_{\kb,\jb}}{\lambda_{\kb}^{(n)}}, \quad \kb,\jb \in \Lambda_n.
\end{align}
\end{thm}

\begin{proof}
Let $f$ be a function invariant under $\G$. According to \eqref{H=Lambda},
\begin{align*}
 \sum_{\jb\in \HH_n^{\circ}} c_{\jb}^{(n)} f(\tfrac{\jb}{(d+1)n})
&   = \sum_{p\in \mathcal{C}_{d+1}} \binom{d+1}{p}\sum_{\jb\in \Lambda_n^{\circ,p}}
  c_{\jb}^{(n)} f(\tfrac{\jb}{(d+1)n}),\\
  \sum_{\jb\in\HH_n^{*} \setminus \HH_n^{\circ}} c_{\jb}^{(n)}f(\tfrac{\jb}{(d+1)n})
& = \sum_{p\in \mathcal{C}_{d+1}} \binom{d+1}{p}\sum_{\jb\in \Lambda_n^{\partial,p}}
  c_{\jb}^{(n)} f(\tfrac{\jb}{(d+1)n}).
\end{align*}
Adding these two expressions together, we conclude that
\begin{align} \label{eq:T2H}
       \sum_{\jb\in \HH_n^*} c_{\jb}^{(n)}f(\tfrac{\jb}{(d+1)n})
 & = \sum_{p\in \mathcal{C}_{d+1}} \binom{d+1}{p}
\left[ \sum_{\jb\in \Lambda_n^{\circ,p}} c_{\jb}^{(n)} f(\tfrac{\jb}{(d+1)n})
+  \sum_{\jb\in \Lambda_n^{\partial,p}} c_{\jb}^{(n)} f(\tfrac{\jb}{(d+1)n})
\right],
\\ &= \sum_{\jb\in \Lambda_n} c_{\jb}^{(n)}  \binom{d+1}{p} f(\tfrac{\jb}{(d+1)n})
= \sum_{\jb\in \Lambda_n} \lambda_{\jb}^{(n)} f(\tfrac{\jb}{(d+1)n}). \notag
\end{align}
Replacing $f$ by $f\bar g$, we have proved that $\la f, g \ra_n^* =
\la f, g \ra_{\triangle,n}$ whenever $f \bar g$ is invariant. Hence,
\eqref{eq:orth-Tn} follows from Theorem \ref{ipdH}. Furthermore,
since  $\frac{1}{|\Omega_H|}\int_{\Omega_H} f(\tb) d\tb =
\frac{1}{|\triangle|_H}\int_{\triangle_H} f(\tb) d\tb$ for all
invariant $f$, \eqref{eq:cubature-T} follows from Theorem
\ref{th:cubature-H}.

Furthermore, replacing $f$ by $\TC_{\kb}\overline{\TC_{\jb}}$ in
\eqref{eq:T2H}, we derive by \eqref{eq:TCs} that
\begin{align*}
\langle \TC_{\kb},\TC_{\jb} \rangle_{\triangle,n}
   & = \frac{1}{(d+1)n^{d}}\sum_{\lb\in \HH_n^*}c_{\lb}^{(n)}\TC_{\kb}(\tfrac{\lb}{(d+1)n})\overline{\TC_{\jb}(\tfrac{\lb}{(d+1)n})} \\
   & = \frac{1}{(d+1)n^{d}}\frac{1}{|\G|}\sum_{\lb\in \HH_n^*}c_{\lb}^{(n)}\sum_{\sigma\in \G}\phi_{\kb}(\tfrac{\lb \sigma}{(d+1)n}) \overline{\TC_{\jb}(\tfrac{\lb}{(d+1)n})} \\
   & = \frac{1}{(d+1)n^{d}}\frac{1}{|\G|}\sum_{\lb\in \HH_n^*}c_{\lb}^{(n)}\sum_{\sigma\in \G}\phi_{\kb}(\tfrac{\lb}{(d+1)n}) \overline{\TC_{\jb}(\tfrac{\lb\sigma}{(d+1)n})} \\
  & = \frac{1}{(d+1)n^{d}}\frac{|\G|}{|\G|}\sum_{\lb\in \HH_n^*}c_{\lb}^{(n)}\phi_{\kb}(\tfrac{\lb}{(d+1)n}) \overline{\TC_{\jb}(\tfrac{\lb}{(d+1)n})}
     = \langle \phi_{\kb}, \TC_{\jb} \rangle_{n}^*.
\end{align*}
Using \eqref{k=jmod} and abbreviating  $\kb \equiv \jb \mod (d+1)
\ZZ_\HH^{d+1}$ as $\kb \equiv \jb$, we further deduce that
\begin{align*}
 \langle \TC_{\kb},\TC_{\jb} \rangle_{\triangle,n}
  & = \frac1{|\G|}\sum_{\sigma\in \G} \langle \phi_{\kb}, \phi_{\jb\sigma} \rangle_{n}^*
     = \frac1{(d+1)!} \Big |\left \{\sigma\in \G:  \jb \sigma \equiv \kb \right\}\Big| \\
  & = \frac{\delta_{\jb,\kb}}{(d+1)!} \Big|
           \left\{\sigma\in \G:  \kb \sigma \equiv \kb \right\}\Big|
      = \frac{\delta_{\jb,\kb}}{\lambda^{(n)}_{\kb}},
\end{align*}
where the last equality follows readily from the definition of
$\Lambda_n$ and \eqref{lambdaj}.
\end{proof}

Since the above proof applies to invariant $f$, it also applies to $f, g \in \CTS_n$
since $f \bar g $ is invariant if both $f$ and $g$ are anti-invariant. Recall that
$\TS_{\kb}(\frac{\jb}{(d+1)n}) =0$ when $\jb \in  \Lambda_n\setminus
\Lambda_n^{\circ}$, we have also proven the following result.

\begin{thm}
Let the discrete inner product $\langle \cdot, \cdot
\rangle_{\triangle^{\circ},n}$ be defined by
\begin{align*}
\langle f, g \rangle_{\triangle^{\circ},n} = \frac{d!}{n^{d}}
\sum_{\jb \in \Lambda_n^{\circ}} f(\tfrac{\jb}{(d+1)n})
\overline{g(\tfrac{\jb}{(d+1)n})}.
\end{align*}
Then
\begin{align*}
\langle f,g \rangle_{\triangle^{\circ},n} =  \langle f,g
\rangle_{\triangle_H},\qquad
                f, g \in \CTS_{n}.
\end{align*}
\end{thm}


\subsection{Interpolation on the simplex}

Using invariance and the fact that the fundamental simplex $\triangle_H$ is
the building block of the fundamental domain, we can also deduce results on interpolation on the fundamental simplex. For $d =3$ the results in this
subsection have appeared in \cite{LX08} and the proof there can be followed
verbatim when $3$ is replaced by $d$. Thus, we shall omit the proof. Recall that
the operator  $\mathcal{P}^{\pm}$ is defined in \eqref{eq:Ppm}.

\begin{thm} \label{thm:4.6}
For $n>0$, and $f\in C(\Delta_H)$, define
\begin{align*}
   \mathcal{L}_nf(\tb):= \sum_{\jb\in \Lambda_n^{\circ}}
    f(\tfrac{\jb}{(d+1)n}) \ell_{\jb,n}^{\circ}(\tb), \qquad
  \ell_{\jb,n}^{\circ} (\tb):= \frac{d!(d+1)!}{n^{d}} \sum_{\kb\in \Lambda_n^{\circ}}
   \TS_{\kb}(\tb)    \overline{\TS_{\kb}(\tfrac{\jb}{(d+1)n})}.
\end{align*}
Then $\mathcal{L}_n f$ is the unique function in $\CTS_n$ that
satisfies
\begin{align*}
\mathcal{L}_nf( \tfrac{\jb}{(d+1)n}) =f( \tfrac{\jb}{(d+1)n}),
\qquad \jb\in \Lambda_n^{\circ}.
\end{align*}
Furthermore, the fundamental interpolation function
$\ell_{\jb,n}^{\circ}$ is real and satisfies
\begin{align*}
   \ell_{\jb,n}^{\circ} (\tb) = \frac{d!}{n^{d}} \mathcal{P}^{-}_{\tb} \left[
   \Theta_{n}(\tb-\tfrac{\jb}{(d+1)n}) - \Theta_{n-1}(\tb-\tfrac{\jb}{(d+1)n})
   \right],
\end{align*}
where $\mathcal{P}^{-}_{\tb}$ means that the operator
$\mathcal{P}^{-}$ is acting on the variable $\tb$ and $\Theta_{n}$
is defined in \eqref{DnTheta-n}.
\end{thm}

The function $\CL_n f$ interpolates at the interior points of $\Lambda_n$.
We can also consider interpolation on $\Lambda_n$ by working with the
operator $\CI^*_n f$ in Theorem \ref{interp-H}, which interpolates $f$ on
$\HH_n^\circ$ but interpolates a sum over congruent points on
$\HH_n^*\subset \HH_n^\circ$. Using invariance, however, $\CI_n^* f$
leads to a genuine interpolation operator on $\Lambda_n$.

\begin{thm} \label{thm:4.7}
For $n>0$ and $f\in C(\Delta_H)$ define
\begin{align*}
  \mathcal{L}^*_n f(\tb):= \sum_{\jb\in \Lambda_n} f(\tfrac{\jb}{(d+1)n})
     \ell^{\triangle}_{\jb,n}(\tb),
 \qquad  \ell^{\triangle}_{\jb,n}(\tb) := \frac{\lambda_{\jb}^{(n)}}{(d+1)n^{d}}
 \sum_{\kb\in \Lambda_n}\lambda^{(n)}_{\kb} \TC_{\kb}(\tb)
       \overline{\TC_{\kb}(\tfrac{\jb}{(d+1)n})}.
\end{align*}
Then $\mathcal{L}^*_n f$ is the unique function in $\CTC_n$ that
satisfies
\begin{align*}
\mathcal{L}^*_n f(\tfrac{\jb}{(d+1)n}) = f(\tfrac{\jb}{(d+1)n}),
\qquad \jb\in \Lambda_n.
\end{align*}
Furthermore, the fundamental interpolation function
$\ell^{\triangle}_{\jb,n}$ is given by
\begin{align*}
\ell^{\triangle}_{\jb,n}(\tb) = \lambda_{\jb}^{(n)} \mathcal{P}^{+}
\ell_{\jb,n}(\tb),
\end{align*}
where $\ell_{\jb,n}$ is defined in Theorem \ref{interp-H} and
has a compact formula.
\end{thm}

Let $\|\mathcal{L}_n\|$ and $\|\mathcal{L}_n^*\|$ denote the
operator norms of $\mathcal{L}_n$ and $\mathcal{L}_n^*$,
respectively, both as operators from $C(\triangle_H)\mapsto
C(\triangle_H)$. From Theorems \ref{thm:4.6} and \ref{thm:4.7},
an immediate application of Theorem \ref{LebesgueH} yields
the following theorem.

\begin{thm}
There is a constant $c$ independent of $n$, such that
\begin{align*}
  \|\mathcal{L}_n\| \leq c (\log n)^d \quad\hbox{and}\quad
    \|\mathcal{L}_n^*\| \leq c (\log n)^d.
\end{align*}
\end{thm}

It should be pointed out that the interpolation functions defined in
these theorems are analogous of trigonometric polynomial
interpolation on equally spaced points \cite[Chapt. X]{Z}. These
are trigonometric interpolation on equal spaced points in the
simplex $\triangle_H$, which can be easily transformed to
interpolation on regular simplex, say $\{y: 0 \le y_d \le \ldots \le y_1 \le 1\}$
in $\RR^d$. These interpolation functions are real, easily computable
from their compact formulas, and have small Lebesgue constants.
They should be the ideal tool for interpolation on the simplex
in $\RR^d$.


\section{Generalized Chebyshev polynomials and their zeros}\label{Sec:Chebyshev}
\setcounter{equation}{0}

The generalized sine and cosine functions can be used to define analogues of
Chebyshev polynomials of the first and the second kind, respectively, which are
algebraic orthogonal polynomials of $d$-variables, just as in the classical case
of one variable. These polynomials have been defined in the literature, as noted
in the Introduction. In the first subsection we define these polynomials and
present their basic properties. Most of the results in this subsection are not
new,  however, we shall present a coherent and independent treatment,  and
some of the results on recurrence relations appear to be new. In the second
subsection, we study the common zeros of these polynomials and use them
to establish a family of Gaussian cubature formulas, which exist rarely.

\subsection{Generalized Chebyshev polynomials}
The generalized trigonometric functions can be transformed into polynomials
under a change of variables $z: \RR_H^{d+1}\rightarrow \CC^d$, defined by
\begin{align} \label{z}
\begin{split}
& z_k = z_k(\tb):=   \TC_{\vb^k}(\tb)  =
\frac{1}{|\vb^k\CG|}\sum_{\jb\in \vb^k\CG} \e^{\frac{2\pi}{d+1}
\jb\, \cdot \tb},  \quad k=1,2,\dots,d,
\end{split}
\end{align}
where $\vb^k= (\{d+1-k\}^k, \{-k\}^{d+1-k})$, and
$\frac{\vb^k}{d+1}\ ( k=1,2,\dots,d)$ are vertices of the
fundamental triangle $\triangle_H$ defined in Section 2. It is easy
to see that $\overline{z}_k = z_{d+1-k}$. The homogeneity of $\tb$
shows that $\vb^k \cdot \tb = (d+1) (t_1 +\ldots + t_k)$ and,
consequently,
\begin{align}
\begin{split}
  z_k = \frac{1}{\binom{d+1}{k}}\sum_{J\subset \NN_{d+1} \atop |J|=k }
       \e^{2  \pi i \sum_{j\in J} t_j },
\end{split}
\end{align}
which shows that $z_1,\ldots, z_d$ are the first $d$ elementary
symmetric polynomials of $\e^{2\pi i t_1}, \ldots, \e^{2\pi i
t_{d+1}}$. The same change of variables are used in \cite{Be}.

Since $\TC_\kb(\tb)$, $\kb \in \Lambda$, is evidently a symmetric
polynomial in $\e^{2\pi i t_1}, \ldots, \e^{2\pi i t_{d+1}}$, it is
a polynomial in $z_1,\ldots, z_d$. Some of its properties can be
derived from the recursive relations given below.

\begin{lem}
\label{lm:recur}
The generalized sine and cosine functions satisfy the recurrence relations,
\begin{align}
\label{eq:recurTCTC}
 &\TC_{\jb}(\tb) \TC_{\kb}(\tb) = \frac{1}{|\CG|} \sum_{\sigma\in \CG}
    \TC_{\kb+\jb\sigma}(\tb),\qquad \jb,\kb \in \Lambda,\\
\label{eq:recurTCTS}
 &\TC_{\jb}(\tb) \TS_{\kb}(\tb) = \frac{1}{|\CG|} \sum_{\sigma\in \CG}
    \TS_{\kb+\jb\sigma}(\tb),\qquad \jb,\kb \in \Lambda,\\
\label{eq:recurTSTS}
 &\TS_{\jb}(\tb) \TS_{\kb}(\tb) = \frac{1}{|\CG|} \sum_{\sigma\in \CG}(-1)^{|\sigma|}\,
     \TC_{\kb+\jb\sigma}(\tb),\quad \jb,\kb \in \Lambda.
 \end{align}
\end{lem}

\begin{proof}
From the definition of $\TC_{\kb}$, we obtain that
\begin{align*}
 \TC_{\kb}(\tb) \TC_{\kb}(\tb) &=\frac{1}{|\CG|} \sum_{\sigma\in \CG} \phi_{\kb \sigma}( \tb)  \times
  \frac{1}{|\CG|} \sum_{\sigma\in \CG}\phi_{\kb\sigma}( \tb)
  =\frac{1}{|\CG|^2} \sum_{\tau\in\CG}  \sum_{\sigma\in \CG}
      \phi_{\kb\tau+\jb\sigma}( \tb)\\
   &= \frac{1}{|\CG|^2} \sum_{\tau\in\CG}  \sum_{\sigma\in \CG}
       \phi_{(\kb+\jb \sigma\tau^{-1}) \tau} ( \tb) =
        \frac{1}{|\CG|^2} \sum_{\tau\in\CG}  \sum_{\sigma\in
       \CG} \phi_{(\kb+\jb \sigma) \tau} ( \tb)\\
   &=\frac{1}{|\CG|} \sum_{\sigma \in \CG} \TC_{\kb+\jb\sigma}(\tb),
 \end{align*}
which proves \eqref{eq:recurTCTC}. The other two relations,
\eqref{eq:recurTCTS} and \eqref{eq:recurTSTS}, can be established
similarly.
\end{proof}

The polynomials defined by $\TC_\kb(\tb)$ under \eqref{z} are
analogue of Chebyshev polynomials of the first kind, to be defined
formerly below. We will also define Chebyshev polynomial of the
second kind, for which we need the following lemma.

\begin{lem} \label{lem:zkTSv0}
Let $\vb^{\circ}: =\big(\frac{d(d+1)}{2},\frac{(d-2)(d+1)}{2}, \dots,
  \frac{(2-d)(d+1)}{2},\frac{-d(d+1)}{2}\big)$, then
\begin{align} \label{TSv0}
  \TS_{\vb^{\circ}}(\tb)=\frac{1}{(d+1)!} \prod_{1\le \mu < \nu \le d+1}
      \left(\e^{2 i\pi t_{\mu}}-\e^{2 i\pi t_{\nu}} \right)
  = \frac{(2i)^{\frac{d(d+1)}{2}}}{(d+1)!} \prod_{1\le \mu < \nu \le d+1}
    \sin  \pi (t_{\mu}-t_{\nu}).
 \end{align}
Furthermore,
\begin{align} \label{zkTSv0}
 z_k \TS_{\vb^{\circ}}(\tb)  = \frac{k!(d+1-k)!}{(d+1)!}\TS_{\vb^{\circ}+\vb^k}(\tb),
     \qquad 1\le   k\le d.
\end{align}
\end{lem}

\begin{proof}
Let $\beta : = (d, d-1, \ldots, 1,0) \in \NN_0^{d+1}$. Then
$$
  \e^{\frac{2\pi i }{d+1} \vb^\circ \cdot \tb} = \e^{- \pi i d ( t_1+\ldots +t_{d+1})}
      \e^{2 \pi i \beta \cdot \tb}  =    \e^{2 \pi i \beta \cdot \tb}.
$$
Using the well-known Vandermond determinant (see, for example,
\cite[p. 40]{Mac}),
$$
  \sum_{\sigma \in S_{d+1}} \rho (\sigma) \sigma (x^\beta)
     = \det (x_i^{d+1-j})_{1 \le i,j \le d+1} = \prod_{1 \le i < j \le d+1}(x_i-x_j),
$$
and setting $x_j = \e^{2\pi i t_j}$, we obtain
$$
\TS_{\vb^\circ}(\tb) = \frac{1}{|\G|} \sum_{\sigma \in \G}
           \rho(\sigma)  \e^{\frac{2\pi i }{d+1} \vb^\circ \cdot \tb}
  =\frac{1}{(d+1)!} \prod_{1\le \mu < \nu \le d+1}
      \left(\e^{2 i\pi t_{\mu}}-\e^{2 i\pi t_{\nu}} \right).
$$
Furthermore, since $t_1+\ldots + t_{d+1} =0$, the above equation immediately gives
$$
\TS_{\vb^\circ}(\tb)
  = \frac{1}{(d+1)!} \prod_{1\le \mu < \nu \le d+1}
      \left(\e^{i\pi t_{\mu - \nu}}-\e^{i\pi t_{\nu - \mu}} \right),
$$
from which the second equal sign of \eqref{TSv0} follows.

Next we prove \eqref{zkTSv0}. Using \eqref{eq:recurTCTS}, we have
\begin{align*}
    z_k \TS_{\vb^{\circ}} (\tb) = \TC_{\vb^{k}} (\tb) \TS_{\vb^{\circ}} (\tb)
    = \frac{1}{(d+1)!}\sum_{\sigma\in \CG}
    \TS_{\vb^{\circ}+\vb^{k}\sigma}(\tb).
\end{align*}
Assume $\lb := \vb^{k}\sigma \neq \vb^{k}$ for some $\sigma \in \G$.
By the definition of $\vb^k$, there exists an integer $1 \le j\le d$
such that $l_j= -k$ and $l_{j+1}=d+1-k$.  Consequently,
$v^{\circ}_j+l_j = \frac{(d+2-2j)(d+1)}{2}-k =
\frac{(d-2j)(d+1)}{2}+(d+1-k)=v^{\circ}_{j+1}+l_{j+1}$, which
implies that  $\vb^{\circ}+\vb^{k}\sigma\in \partial\Lambda$ and
$\TS_{\vb^{\circ}+\vb^{k}\sigma}(\tb)=0$. Thus, by the definition of
the stabilizer $\G_{\vb^k}$,
\begin{align*}
  z_k \TS_{\vb^{\circ}} (\tb)
    = \frac{1}{(d+1)!}\sum_{\sigma\in \CG_{\vb^k}}
    \TS_{\vb^{\circ}+\vb^{k}\sigma}(\tb)
    = \frac{k!(d+1-k)!}{(d+1)!} \TS_{\vb^{\circ}+\vb^{k}}(\tb),
\end{align*}
as the stabilizer of $\vb^k$ has cardinality $k!(d+1-k)!$.
\end{proof}

The same argument that proves \eqref{TSv0} also shows that, for
each $\kb \in \Lambda$,
$$
\TS_{\kb + \vb^\circ}(\tb) = \det \left(x_i^{\lambda_j + \beta}\right)_{1 \le i,j\le d+1},
 \qquad \lambda := (k_1 - k_{d+1}, k_2 - k_{d+1}, \ldots,
  k_d - k_{d+1},0),
$$
where $x_i = \e^{2\pi i t_i}$. The definition of $\Lambda$ shows
that $\lambda$ is a partition. Hence, according to \cite[p.
40]{Mac},  $\TS_{\kb + \vb^\circ}(\tb)$ is divisible by
$\TS_{\vb^\circ}$ in the ring $\ZZ[x_1,\ldots, x_d]$, and the
quotient
$$
  s_\lambda (x_1,\ldots, x_d)= \TS_{\kb + \vb^\circ}(\tb) / \TS_{\vb^\circ}(\tb)
$$
is a symmetric polynomial in $x_1,\ldots, x_d$, which is the Schur function
in the variables $x_1,\ldots, x_d$ corresponding to the partition $\lambda$.
Since this ratio is a symmetric polynomial, it is then a polynomial in the
elementary symmetric polynomials $z_1,\ldots, z_d$; more precisely, it
is a polynomial in $z$ of degree $(k_1-k_{d+1})/(d+1)$ as shown by the
formula \cite[(3.5)]{Mac}. This is our analogue of Chebyshev polynomials
of the second kind.

To simplify the notation, we find it convenient to change the index and
define the Chebyshev polynomials of the first and the second kind formally
as follows:

\begin{defn}
Define the index mapping $\a: \Lambda \rightarrow \NN_0^d$,
\begin{align}
\label{eq:kappa}
     \a_{i} =\a_i(\kb):= \frac{k_i-k_{i+1}}{d+1},\qquad 1\le i\le d.
\end{align}
Under the change of variables \eqref{z} and \eqref{eq:kappa}, define
\begin{align*}
   T_{\a}(z) : = \TC_{\kb}(\tb) \quad \hbox{and} \quad
    U_{\a}(z) := \frac{\TS_{\kb+\vb^{\circ}}(\tb)}{\TS_{\vb^{\circ}}(\tb)},
     \qquad \a \in \NN_0^d.
\end{align*}
We call $T_{\a}(z)$ and $U_{\a}(z)$ Chebyshev polynomials of the first
and the second kind, respectively.
\end{defn}

It is easy to see that the mapping \eqref{eq:kappa} is an isomorphism. Indeed,
since $k_1 +\ldots + k_{d+1} =0$ for $\kb \in \Lambda$, it is easy to see that
the inverse of $\alpha$ is given by
\begin{align*}
 k_i = \sum_{\mu=1}^d \sum_{\nu =1}^{\mu}
     \a_{\nu} - (d+1)\sum_{\nu=1}^{i-1}  \a_{\nu}, \qquad 1\le i\le d+1.
\end{align*}
We also note that $\alpha(\vb^k) = \epsilon_k : = (\{0\}^{k-1},1, \{0\}^{d-k})$, the
$k$-th element of the standard basis for the Euclidean space $\RR^d$.

In the recurrence relation that we shall state in a moment, we will
need the definition of $T_\a(z)$ in which $\a$ can have negative
components. If $\a$ has a negative component, say $\alpha_i < 0$,
then $k_i < k_{i+1}$, so that $\kb$ in \eqref{eq:kappa} does not
belong to $\Lambda$.  If $\kb \in \HH$, we define by $\kb^+$ the
rearrangement of $\kb$ such that $\kb^+ \in \Lambda$. By the
definition of $\TC_\kb$, we have $\TC_{\kb^+}(\tb) = \TC_\kb (\tb)$
for all $\kb \in \HH$. Thus, if $\a$ has negative component, then we
define
$$
    T_\a(z) = T_{\a^+}(z)  \qquad \hbox{where $\a^+$ corresponds to $\kb^+$
       by \eqref{eq:kappa}}.
$$

Both $T_{\a}(z)$ and $U_{\a}(z)$ are polynomials of degree $|\a| = \a_1+\ldots
+ \a_d$ in $z$. Moreover, both of them satisfy a simple recursive relation,
which we  summarize in the following theorem.

\begin{thm}
Let $P_{\a}$ denote either $T_{\a}$ or $U_{\a}$. Then
\begin{equation}\label{conjugT}
  \overline{P_{\a}(z)}  = P_{\a_d,\a_{d-1},\dots, \a_1}(z),\quad   \a \in \NN^{d}_0,
\end{equation}
and they satisfy the recursion relation
\begin{align} \label{recurT}
    \binom{d+1}{i} z_i P_{\a} (z) = \sum_{\jb\in \vb^i\CG } P_{\a+\a(\jb)}(z),
        \quad \a \in \NN^d_0,
\end{align}
in which the components of $\a(\jb)$, $\jb \in \vb^i \G$, have values in
$\{ -1,0,1\}$, $U_\a(z) = 0$ whenever $\a$ has a component $\a_i =-1$,
and
\begin{align*}
& T_0(z)=1, \quad T_{\epsilon_k}(z)= z_k, \qquad 1\le k\le d\\
& U_0(z)=1, \quad U_{\epsilon_k}(z)= \binom{d+1}{k} z_k, \quad 1\le k\le d.
\end{align*}
\end{thm}

\begin{proof}
The relation \eqref{conjugT} follows readily form the fact that
$-(k_i-k_j) = k_j -k_i$ and \eqref{eq:kappa}. The relation
\eqref{recurT} follows immediately from \eqref{eq:recurTCTC} and
\eqref{eq:recurTCTS}. The values of $P_0$ and $P_{\epsilon_k}$
follow from definitions and \eqref{zkTSv0}. If $\a$ has a component
$\a_i = -1$ then, by \eqref{eq:kappa}, $k_i(\a) = k_{i+1}(\a) -
(d+1)$. A quick computation shows then $k_i(\a) + v^\circ_i =
k_{i+1}(\a) + v^\circ_{i+1}$, which implies that $\kb + \vb^\circ
\in \partial \Lambda$, so that $\TS_{ \kb + \vb^\circ}(\tb)=0$ and
$U_\a(z) =0$.
\end{proof}

It is easy to see that $|\alpha(\jb)|$, $\jb \in \vb^i$, also takes value in
$\{-1,0,1\}$. As a result, in terms of the total degree $|\a|$ of $P_\a$,
the right hand side of \eqref{recurT} contains only polynomials of degree
$n-1, n$ and $n+1$, so that it is a three-term relation in that sense.
For $d =2$, the relation \eqref{recurT} can be rewritten as a recursive
relation,
\begin{align*}
 P_{\a+\epsilon_1}(z) & = 3 z_1 P_\a(z) - P_{\a+(-1,1)}(z) - P_{\a-\epsilon_2}(z),\\
 P_{\a+\epsilon_2}(z) & = 3 z_2 P_\a(z) - P_{\a+({1,-1})}(z) - P_{\a-\epsilon_1}(z)
\end{align*}
which can then be used to generate a polynomial $P_\a$ from lower degree
polynomials recursively. The same, however, cannot be said for $d \ge 3$.
For example, if $d =3$ then \eqref{recurT}  for $k=2$ is
$$
 6 z_2 P_\alpha(z) = P_{\a+ \epsilon_2}(z) + P_{\a+(1,-1,1)}(z) +
   P_{\a+(1,0,-1)}(z)+ P_{\a+(-1,0,1)}(z) + P_{\a- \epsilon_1} +
   P_{\a+(-1,1,-1)}(z),
$$
which has two polynomials of $|\a| +1$ in the right hand side,
$P_{\a+ \epsilon_2}(z)$ and $P_{\a+(1,-1,1)}(z)$.  It is possible to combine
the relations in \eqref{recurT} to write $P_{\a + \varepsilon_i}(z) =
\sum_{i=1}^d \binom{d+1}{i} z_i P_\a(z) + Q(z)$, where $Q(z)$
contains only linear combinations of $\{P_\beta\}$ for $|\beta| = |\a|$
and $|\beta| = |\a|-1$. In lower dimension, the exact forms of $Q$
can be easily determined (for $d =3$ see \cite{Sun3}), but the general
formula for a generic $d$ appears to be complicated and we shall
nor pursuit it here.

Next we show that $T_\a$ and $U_\a$ are orthogonal polynomials. The integral
of the orthogonality is taken over the region that is the image of $\triangle_H$
with respect to a measure, or weight function, that comes from the Jacobian of
the change of variables. Furthermore, since $\bar z_k \mapsto z_{d+1-k}$, we
can consider real coordinates, denoted by $x$,
\begin{align} \label{realOP}
\begin{split}
& x_k = \frac{z_k+z_{d+1-k}}{2}, \quad x_{d+1-k} = \frac{z_k-z_{d+1-k}}{2i}
        \qquad \hbox{for} \quad 1\le k\le \lfloor \tfrac{d}{2}\rfloor, \\
 & \hbox{and}\quad  x_{\tfrac{d+1}{2}} = \frac1{\sqrt{2}}\, z_{\tfrac{d+1}{2}}, \qquad
     \hbox{if $d$ is odd}.
\end{split}
\end{align}
Combining \eqref{z} and \eqref{realOP}, our change variable becomes
$\RR_H^{d+1} \mapsto \RR^d$: $\tb \mapsto x$.  Let us define
$$
  w(x) = w(x(\tb)) := \prod_{1\le \mu<\nu\le d+1}\sin^2\pi (t_{\mu}-t_{\nu}).
$$

\begin{lem}
The Jacobian of the changing variable $\tb \mapsto x$ is given by
\begin{align} \label{eq:Jac1}
  \left|\det \left[\frac{\partial(x_1,x_2,\dots,x_{d})}{\partial(t_1,t_2,\dots,t_{d})} \right]
     \right | =  2^{\frac{d(d+2)}{2}}\prod_{k=1}^d \frac{\pi}{\binom{d+1}{k}}
        [w(x)]^{\frac12}.
\end{align}
\end{lem}

\begin{proof}
Let $\partial_i$ denote the partial derivative with respect to $t_i$ for
$1 \le i \le d$. We first prove that
\begin{align}\label{eq:det}
 \det  \left [ \frac{\partial(z_1,z_2,\dots,z_d)}{\partial(t_1,t_2,\dots,t_d)}\right]
       =\prod_{k=1}^d\frac{2 \pi i}{\binom{d+1}{k}}\times
           \prod_{1\le {\mu}<{\nu}\le d+1 } |\sin \pi (t_{\mu}-t_{\nu})|.
\end{align}
Using the Jacobian of the elementary symmetric polynomials with respect to
its variables and $t_1+\ldots +t_{d+1} =0$, this can be shown as in \cite[(5.9)]{Be}.
We give an inductive proof below.

Regarding $t_1,t_2,\dots,t_{d+1}$ as independent variables, one sees that
\begin{align*}
      \frac{\partial z_k}{\partial t_j} =\frac{ 2\pi i}{\binom{d+1}{k}}
     \sum_{j\in I\subseteq \NN_{d+1}, |I|=k} \e^{2\pi i \sum_{{\nu}\in I} t_{\nu}}.
\end{align*}
For each fixed $j$, split $I \subset \NN_{d+1}$ as two parts, one
contains $\{j, d+1\}$ and one does not, so that after canceling the
common factor, we obtain
\begin{align*}
 \frac{\partial z_k}{\partial t_j}-\frac{\partial z_k}{\partial t_{d+1}}
 & =  \frac{ 2\pi i}{\binom{d+1}{k}}  \bigg(
 \sum_{ j\in  I\subseteq \NN_{d+1}^{\{d+1\}}, |I|=k} \e^{2\pi i \sum_{{\nu}\in I} t_{\nu}}
 - \sum_{d+1\in I\subseteq \NN_{d+1}^{\{j\}}, |I|=k} \e^{2\pi i \sum_{{\nu}\in I} t_{\nu}}
    \bigg) \\
 &=  \frac{ 2\pi i}{\binom{d+1}{k}}  \big(\e^{2\pi i t_j}-\e^{2\pi i t_{d+1}}\big)
       \sum_{I\subseteq \NN^{\{j\}}_{d}, |I|=k-1} \e^{2\pi i \sum_{\nu \in I} t_{\nu}}.
\end{align*}
Hence, setting $f^{d}_{j,k} := \sum_{I\subseteq \NN^{\{j\}}_{d},
|I|=k-1} \e^{2\pi i \sum_{\nu \in I} t_{\nu}}$ and defining the
matrix  $F_d:=\big( f_{j,k}^d\big)_{1\le j,k\le d}$, we have shown
that
\begin{align*}
\det
\bigg(\frac{\partial(z_1,z_2,\dots,z_d)}{\partial(t_1,t_2,\dots,t_d)}
\bigg)  = \prod_{k=1}^d \frac{2\pi i}{\binom{d+1}{k}}\times
\prod_{{j}=1}^d \big(\e^{2\pi i t_{j}}-\e^{2\pi i t_{d+1}}\big)
\times \det(F_{d}).
\end{align*}
By definition, $f^d_{j,1}=1$; splitting $\NN^{\{j\}}$ according to if the part contains
$d$ or not, it follows that
\begin{align*}
 f^{d}_{j,k}-f^{d}_{d,k}
 = \big(\e^{2\pi i t_{d}}-\e^{2\pi i t_{j}}\big) \sum_{I\subseteq \NN_{d}^{\{j,d\}}, |I|=k-2} \e^{2\pi i \sum_{{\nu}\in I} t_{\nu}}
=  \big(\e^{2\pi i t_{d}}-\e^{2\pi i t_{j}}\big) f_{j,k-1}^{d-1},
\end{align*}
which allows us to use induction and show that
\begin{align*}
  \det(F_d) =  \prod_{1\le {\mu}<  d-1} \big(\e^{2\pi i t_{\mu}}-\e^{2\pi i t_{d}}\big)\times \det(F_{d-1}) = \ldots =
\prod_{1\le {\mu}<{\nu} \le d} \big(\e^{2\pi i t_{\mu}}-\e^{2\pi i
t_{\nu}}\big).
\end{align*}
This proves \eqref{eq:det} upon using \eqref{TSv0}. Now, a quick computation
using \eqref{realOP} shows that
\begin{align*}
 \left| \det  \left [ \frac{\partial(x_1,x_2,\dots,x_d)}
      {\partial(t_1,t_2,\dots,t_d)} \right] \right |
& =  \left| \det \left [\frac{\partial(x_1,x_2,\dots,x_d)}{\partial(z_1,z_2,\dots,z_d)} \right ]
 \right | \times
\left| \det \left[ \frac{\partial(z_1,z_2,\dots,z_d)}{\partial(x_1,x_2,\dots,x_d)} \right]
  \right | \\
 & =    2^{\frac{d(d+2)}{2}}\prod_{k=1}^d \frac{\pi}{\binom{d+1}{k}}
      \prod_{1\le {\mu}<{\nu}\le d+1 } |\sin \pi (t_{\mu}-t_{\nu})|,
\end{align*}
which is exactly \eqref{eq:Jac1}.
\end{proof}

Under this change of variables, the domain $\triangle_H$ is mapped to
\begin{align*}
 \Delta^* := \bigg \{x = x(\tb)\in \RR^d: \tb \in \RR_H^{d+1}, \prod_{1\le i< j\le d+1}
  \sin \pi (t_i-t_j) \ge 0 \bigg\},
\end{align*}
which is the image of $\triangle_H$ under the mapping $\tb \mapsto x$. In the
case of $d =2$ and $d=3$, $\Delta^*$ is depicted in the Figure 5.1.

\begin{figure}[htb]
\centering
\hfill\includegraphics[width=0.4\textwidth]{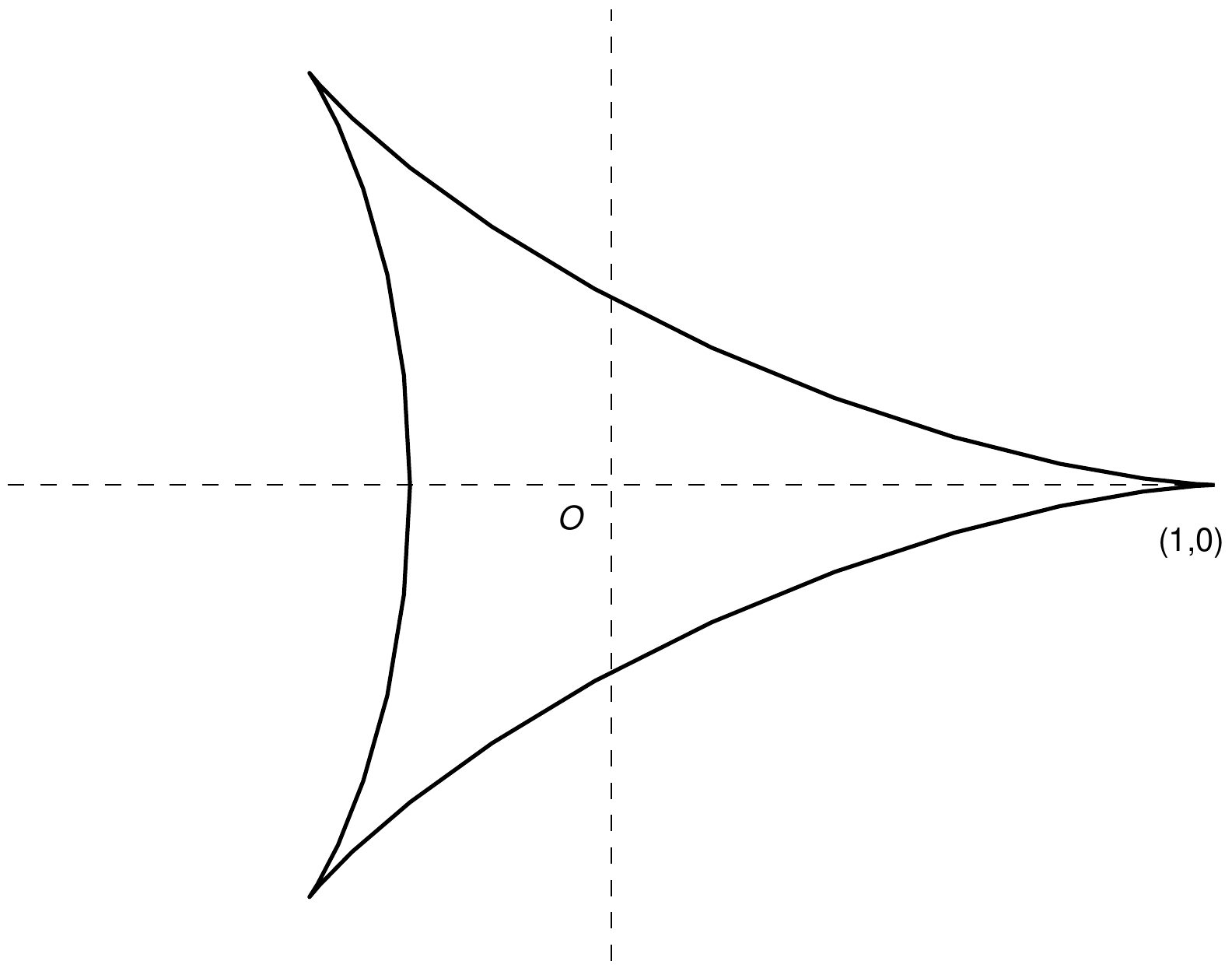}\hfill
\includegraphics[width=0.55\textwidth]{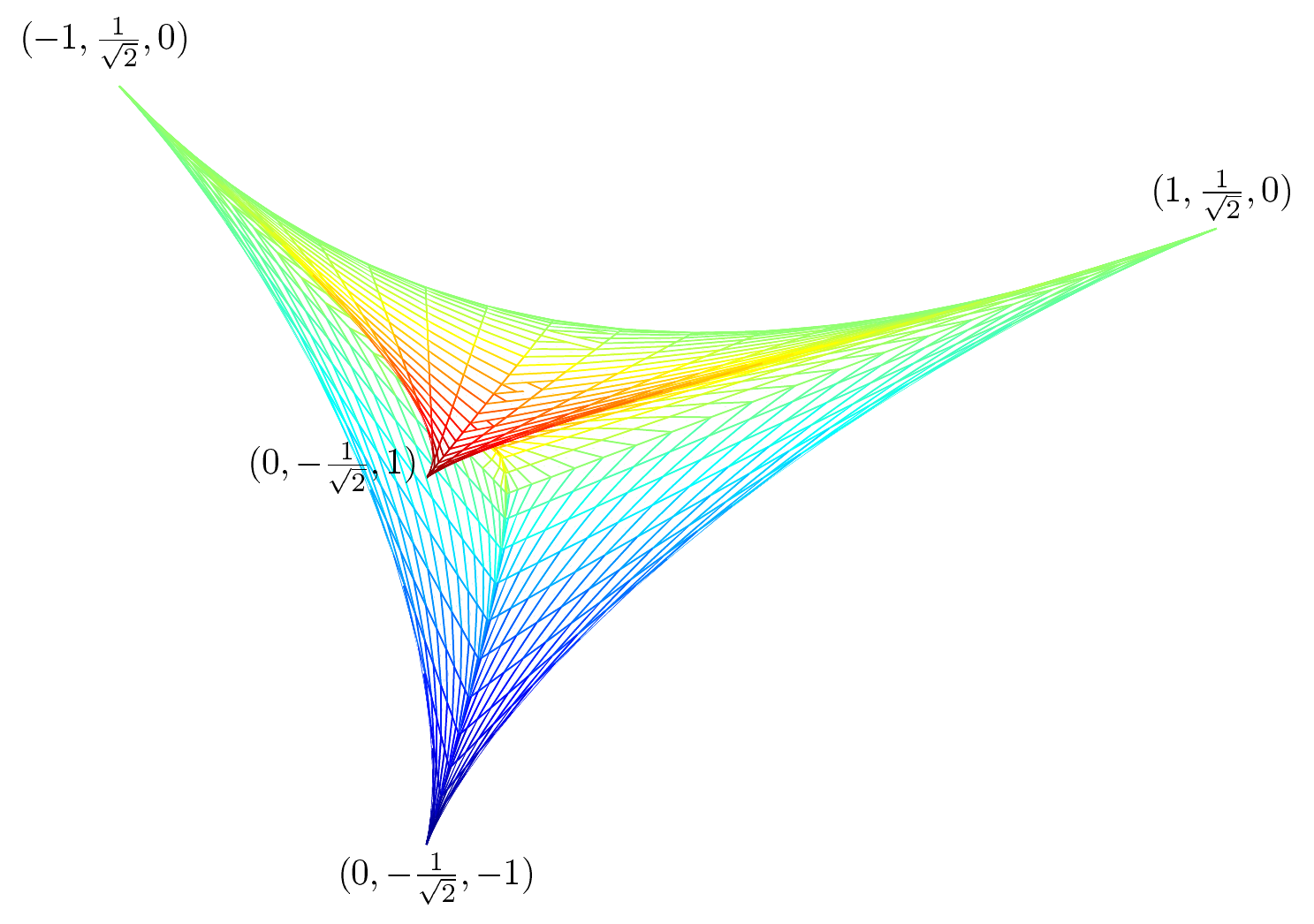}\hspace*{\fill}
\caption{The region $\Delta^*$ for $d=2$ and $d =3$}
\end{figure}

We will need the cases of $\mu = -1/2$ and $\mu = 1/2$ of the
weighted inner product
$$
 \la f, g \ra_{w^\alpha}: =
       c_\alpha \int_{\Delta^*} f(z)\overline{g(z)} w^\alpha(z) dx,
$$
where $c_\alpha$ is a normalization constant, $c_\alpha := 1/\int_{\Delta^*}
w^\alpha(z) dx$. The change of variables $\tb \mapsto x$ shows
immediately that
$$
 \frac{1}{|\Delta_H|} \int_{\Delta_H} f(\tb) d\tb = c_{-\frac{1}{2}}
    \int_{\Delta^*} f(z) w^{-\frac12}(z) dx,
$$
where
$$
c_{-\frac12}= \frac{1}{|\Delta_H|} \times \frac{\omega^{\frac12}(z)}{\Big| \det\big(\frac{\partial(x_1,x_2,\dots,x_d)}{\partial(t_1,t_2,\dots,t_d)}\big) \Big| }
= 2^{-\frac{d(d+2)}{2}} \sqrt{d+1}\, d! \prod_{k=1}^d\frac{\binom{d+1}{k}}{\pi}.
$$
Furthermore, by \eqref{TSv0} and \eqref{eq:recurTSTS},
\begin{align*}
w(z) =  \left(-\tfrac14\right)^{\frac{d(d+1)}{2}}
  \left[(d+1)!\right]^2 \left[\TS_{\vb^{\circ}}(\tb)\right]^2
  =  \left(-\tfrac14\right)^{\frac{d(d+1)}{2}}  (d+1)!\,
     \sum_{\sigma\in \CG} (-1)^{|\sigma|}\,\TC_{\vb^{\circ}+\vb^{\circ}\sigma}(\tb),
\end{align*}
which shows that $w$ is a polynomial in $z$ with a total degree of $2d$, and
furthermore, by \eqref{eq:orth-TS},
\begin{align*}
c_{\frac12} = 1 \Big / \int_{\Delta^*} w^{\frac{1}{2}}(x(\tb)) dx
       = 1 \Big/ \int_{\triangle_H} w(x(\tb)) d\tb =
   \sqrt{\frac{2^{d^2}}{d+1}}\, \prod_{k=1}^d \frac{\binom{d+1}{k}}{\pi}.
\end{align*}

The orthogonality of $T_{\kappa}$ and $U_{\kappa}$, respectively,
then follows from the orthogonality of $\TC_\kb$ and $\TS_\kb$ with
the change of variables. More precisely, Theorem \ref{prop:trig-ortho}
leads to the following theorem:

\begin{thm}
The polynomials $T_\a$ and $U_\a$ are orthogonal polynomials with
respect to $w^{-1/2}$ and $w^{1/2}$, respectively, and
\begin{align}
    \langle T_{\a}, T_{\beta}\rangle_{ w^{-\frac12}} &= d_{\a}\delta_{\a,\beta},
  \qquad  d_{\a} =\frac{1}{\binom{d+1}{p}}\,  \hbox{ if } \a \in
   \wedge^{p},  \label{TCnorm} \\
   \langle U_{\a}, U_{\beta}\rangle_{ w^{\frac12}} &= \frac{\delta_{\a,\beta}}{(d+1)!},
      \quad   \a, \beta \in \NN_0^d, \label{UCnorm}
\end{align}
where $\wedge^{p}$ is the image of $\Lambda^{p}$ under the index mapping \eqref{eq:kappa},
\begin{align*}
\wedge^{p} =\left\{\a \in \NN^d_0:  \a = \a (\kb), \quad  \kb \in \Lambda^{p}  \right\}.
\end{align*}
\end{thm}

The orthogonality of these polynomials has been established in the
literature; we refer to \cite{Be} and the references therein. By
\eqref{dimCT}, the set $\{T_\a: \a \in \NN_0^d\}$ and, respectively,
$\{U_\a: \a \in \NN_0^d\}$ forms a mutually orthogonal basis for the
space $\Pi^d$ of polynomials in $d$-variables with respect to
$w^{-\frac12}$ and, respectively $w^{\frac12}$.


\subsection{Zeros of Chebyshev polynomials of the second kind and
Gaussian cubature}

It is well known that zeros of orthogonal polynomials of one variable
are nodes of the Gaussian quadrature. We consider its extension in
several variables in this subsection. First we review the background
in several variables.

Let $w$ be a nonnegative weight function defined on a compact set
$\Omega$ in $\RR^d$.  A cubature formula of degree $2n-1$ for
the integral with respect to $w$ is a sum of point evaluations that
satisfies
$$
  \int_\Omega f(x) w(x) dx = \sum_{j=1}^N \lambda_j f(x^j),
   \qquad \lambda_j \in \RR
$$
for every $f \in \Pi_{2n-1}^d$. The points $x^j=(x_1^j,x_2^j,\dots,x_d^j)$
are called nodes of and the numbers $\lambda_j$ are called weights
of the cubature. It is well-known that a cubature formula of degree $2n-1$
exists only if $N \ge \dim \Pi_{n-1}^d$. A cubature that attains such a
lower bound is called a Gaussian cubature.

In one variable, a Gaussian quadrature always exists and its nodes
are zeros of orthogonal polynomials. For several variables, however,
Gaussian cubatures exist only rarely. In fact, the first family of
weight functions that admit Gaussian cubature were discovered only
relatively recent in \cite{BSX}.  When they do exist, however, their
nodes are common zeros of orthogonal polynomials. To be more
precise, let $\PP_{n} := \{P_{\kappa}: \kappa\in\NN_0^d,\, |\kappa|
= n\}$ denote a basis of orthonormal polynomials of degree $n$ with
respect to $w(x) dx$. Then, it is known that a Gaussian cubature
exists if and only if $\PP_n$ (every element in $\PP_n$) has $\dim
\Pi_{n-1}^d$ real common zeros; in other words, the polynomial ideal
generated by $\{P_\a: |\a| =n\}$ has a zero dimensional variety of
size $\dim \Pi_{n-1}^d$. See \cite{DX,My,St} for these results and
further discussions.

Below we shall show that a Gaussian cubature exists for $w^{\frac{1}{2}}$,
which makes  $w^{\frac12}$ the second family of examples that admits
Gaussian cubature. First we study the zeros of the Chebyshev polynomials
of the second kind. Because of $\bar z_k = z_{d+1-k}$ and \eqref{conjugT},
the real and complex parts of the polynomials in $\{U_\a(z): |\a| = n\}$ form
a real basis for the space of orthogonal polynomials of degree $n$. Hence,
we can work with the zeros of the complex polynomials $U_\a(z)$.

Let $Y_n$ and $Y_n^{\circ}$ be the image of $\left\{\tfrac{\jb}{(d+1)n}:
\ \jb \in\Lambda_n \right\}$ and $\left\{\tfrac{\jb}{(d+1)n}:\
\jb \in\Lambda_n^{\circ} \right\}$
under the mapping $\tb \mapsto x$ respectively,
\begin{align*}
Y_n :  =  \Big\{z\big(\tfrac{\jb}{(d+1)n} \big)  :\ \jb\in \Lambda_n\Big\}
 \quad \hbox{and} \quad
Y_n^{\circ} :  =  \Big\{z\big(\tfrac{\jb}{(d+1)n} \big)  :\ \jb\in \Lambda_n^{\circ} \Big\}.
\end{align*}

\begin{thm} \label{thm:zeroU}
The set $Y_{n+d}^\circ$ is the variety of the polynomial ideal
$\left \langle U_\alpha(x): |\alpha| = n \right \rangle$.
\end{thm}

\begin{proof}
The definition of $U_\alpha$ and \eqref{eq:kappa} shows that if $U_\a$ has
degree $n$, that is, $|\a| =n$, then $\kb \in \Lambda^\circ$ in
$\TS_{\vb^\circ + \kb}/\TS_{\vb^\circ}$ satisfies $k_1 - k_{d+1} = (d+1)n$.
Hence, in order to show that $U_\a$ vanishes on $Y_{n+d}^\circ$, it suffices to
show, by the definition of $\vb^\circ$, that
\begin{align*}
  \TS_{\kb} \left(\tfrac{\jb}{(d+1)n}\right) =0 \quad
   \hbox{for } \kb\in \Lambda^{\circ}\, \hbox{ and }  k_1-k_{d+1}=(d+1)n.
\end{align*}

Recall that $\sigma_{i,j}$ denotes the transposition that interchanges
$i$ and $j$. A simply computation shows that
\begin{align*}
 &\phi_\kb \left(\tfrac{\jb}{(d+1)n}\right) -
    \phi_\kb \left(\tfrac{\jb \sigma_{1,d+1}}{(d+1)n}\right) \\
& \qquad = 2 i \sin \frac{\pi (j_1-j_{d+1})(k_1-k_{d+1})}{(d+1)^2n}
      \e^{\frac{\pi i}{(d+1)^2n} ((j_1+j_{d+1})(k_1+k_{d+1}) +
        2\sum_{\nu=2}^d j_{\nu} k_{\nu}) },
\end{align*}
which is $0$ whenever $\jb \in \Lambda_n$ and $\kb \in \Lambda_n$ satisfies
$k_1-k_{d+1}=(d+1)n$.  Hence, it follows that
\begin{align*}
\TS_{\kb} \left(\tfrac{\jb}{(d+1)n}\right) = \frac{1}{2 |\CG|}\sum_{\sigma\in \CG}
  \varepsilon_\sigma \left[\phi_\kb \left(\tfrac{\jb}{(d+1)n}\right) -
    \phi_\kb \left(\tfrac{\jb \sigma_{1,d+1}}{(d+1)n}\right) \right]  =0,
\end{align*}
where $\varepsilon_\sigma$ takes value of either $-1$ or $1$.
\end{proof}

By \eqref{dimCT},  $|Y_{n+d}^\circ| = |\Lambda_{n+d}^\circ |=
\binom{n+d-1}{d} = \dim \Pi_{n-1}^d$. According to the general
theory, Theorem \ref{thm:zeroU} implies that the Gaussian cubature
exists for the weight function $w^{\frac{1}{2}}$ on $\Delta^*$. In
fact, the precise form of this cubature is known.

\begin{thm}
For $w^{\frac{1}{2}}$ on $\Delta^*$, a Gaussian cubature formula exists,
which is given explicitly by
\begin{equation} \label{GaussCuba}
  c_{\frac{1}{2}} \int_{\Delta^*} f(x) w^{\frac{1}{2}}(x) dx =
      \frac{2^{d(d+1)}}{(d+1)(n+d)^d} \sum_{y \in Y_{n+d}^{\circ}}
         w(y) f(y), \quad \forall f \in \Pi_{2n-1}^d,
\end{equation}
where the normalization constant $c_{\frac12}=
\sqrt{\frac{2^{d^2}}{d+1}}\, \prod_{k=1}^d
\frac{\binom{d+1}{k}}{\pi}$.
\end{thm}

\begin{proof}
The change of variables $\tb \mapsto x$ in \eqref{z} and \eqref{realOP} shows
\begin{align*}
c_{\frac12}\int_{\Delta^*} f(z) w^{\frac{1}{2}}(z) dx & =
   \frac{[(d+1)!]^2}{|\Delta_H |}  \int_{\Delta_H} f(z(\tb))
    \left|\TS_{\vb^{\circ}}(\tb)\right|^2 d\tb \\
 &   = \frac{d!(d+1)!}{(n+d)^d} \sum_{j \in \Lambda_{n+d}^\circ}
     \left|\TS_{\vb^{\circ}}(\tfrac{\jb}{(n+d)(d+1)})\right|^2
          f(z(\tfrac{\jb}{(n+d)(d+1)})),
\end{align*}
where the last step follows from the fact that
$\TS_{\vb^{\circ}}(\tb)$ vanishes on the boundary of $\Delta$. This
is exactly the cubature \eqref{eq:cubature-T} applied to the
function $f(z(\tb)) \left|\TS_{\vb^{\circ}}(\tb)\right|^2$.
\end{proof}

The set of nodes of a Gaussian cubature is poised for polynomial
interpolation, that is, there is a unique polynomial $P$ in
$\Pi_{n-1}^d$ such that $P(\xi) = f(\xi)$ for all $\xi \in
Y_{n+d}^{\circ}$ and a generic function $f$. In fact, let
$K_n(\cdot,\cdot)$ be the reproducing kernel of $\Pi_n^d$ in
$L^2_{w^{\frac12}}$, which can be written as
$$
    K_n(x,y): = \sum_{|\a|=n} (d+1)! U_\a(z)U_\a(w)
$$
by \eqref{UCnorm}, where $x,y$ is related to $z,w$ according to
\eqref{realOP}, respectively. Then it follows from the general theory
(cf.  \cite{DX}) that the following result holds:

\begin{prop}
The unique interpolation polynomial of degree $n$ on $Y_{n+d}^{\circ}$ in
$\Pi_{n-1}^d$ is given by
\begin{equation} \label{interp-tri-interior}
 \CL_n f(x(\tb)) =  \sum_{\jb\in \Lambda_{n+d}^{\circ}} f(x(\tfrac{\jb}{(n+d)(d+1)}))
                       w(x(\tfrac{\jb}{(n+d)(d+1)}))K_n^*(\tb, \tfrac{\jb}{(n+d)(d+1)}),
\end{equation}
where $K_n^*(\tb, \sb) = K_n (x,y)$ with $x, y$ related to $\tb,\sb$
by \eqref{z} and \eqref{realOP}, respectively.
\end{prop}

In fact, this interpolation operator is exactly the one in Theorem \ref{thm:4.6}
under the change of variables $\tb \mapsto x$. In particular, we can derive
a compact formula for $K_n(x,y)$ as indicated in that theorem.

\subsection{Cubature formula and Chebyshev polynomials of the first kind}

For the Chebyshev polynomials of the first kind, the polynomials in
$\{T_\a: |\a| = n\}$ do not have enough common zeros in general. In fact,
as shown in \cite{LSX}, in the case of $d=2$, the three orthogonal polynomials
of degree $2$ has no common zero at all. As a consequence, there is no
Gaussian cubature for the weigh function $w^{-\frac12}$ on $\Delta^*$ in
general.

However, the change of variables $\tb \mapsto x$ shows that
\eqref{eq:cubature-T} leads to a cubature of degree $2n-1$ with
respect to $w^{-\frac12}$ based on the nodes of $Y_n$.

\begin{thm}
For the weight function $w^{-\frac{1}{2}}$ on $\Delta^*$ the
cubature formula
\begin{equation} \label{GaussCuba2}
  c_{- \frac{1}{2}} \int_{\Delta^*} f(x) w^{-\frac{1}{2}}(x) dx =
     \frac{1}{(d+1)n^d}  \sum_{\jb\in \Lambda_n}\lambda^{(n)}_{\jb}
            f\big(z\big(\tfrac{\jb}{n(d+1)}\big)\big), \quad \forall f \in \Pi_{2n-1},
\end{equation}
holds, where $c_{-\frac12}=  2^{-\frac{d(d+2)}{2}} \sqrt{d+1}\, d!
\prod_{k=1}^d\frac{\binom{d+1}{k}}{\pi}$ and $\lambda_\jb^{(n)}$ is
given by \eqref{lambdaj}.
\end{thm}

It follows from \eqref{dimCT} that the $|Y_n| = \binom{n+d}{d} =
\dim \Pi_n^d$ and $Y_n$ includes points on the boundary of
$\Delta^*$, hence, the cubature in \eqref{GaussCuba2} is an analog
of the Gauss-Lobatto type cubature for $w^{-\frac12}$ on $w^{-\frac12}$.
The number of nodes of this cubature is more than the lower bound of
$\dim \Pi_{n-1}^d$. Such a formula can be characterized by the
polynomial ideal that has the set of the nodes as its variety. Indeed,
according to a theorem in \cite{X00}, the set of nodes $Y_n$ of
the formula \eqref{GaussCuba2} must be the variety of a polynomial
ideal generated by $\dim \Pi^{d}_{n+1}$ linearly independent
polynomials of degree $n+1$, and these polynomials are necessarily
quasi-orthogonal in the sense that they are orthogonal to all polynomials
of degree $n-2$. Giving the complicated relation between ideals and varieties,
it is of interesting to identify the polynomial idea that generates the
cubature \eqref{GaussCuba2}. This is given in the following theorem.

\begin{thm} \label{prop1}
Let $\vb^*=\kappa^{-1}(e^*)=(d+1,\{0\}^{d-1},-d-1)$.  For $\a = \a (\kb)$
as in \eqref{eq:kappa} and let $\a^* := \a\left((\kb - \vb^*)^+\right)$ .
Then $Y_n$ is the variety of the polynomial ideal
\begin{align} \label{ideal}
  \left \langle T_{\a}(x) - T_{(\a^*)}(x): \quad |\a| = n+1 \right \rangle.
\end{align}
Furthermore, the polynomial $T_{\a}(x) - T_{(\a^*)}(x)$ are of degree
$n+1$ and orthogonal to all polynomials in $\Pi_{n-2}$ with respect to
$w^{-\frac12}$.
\end{thm}

\begin{proof}
From the definition of $\TC_{\kb}$, we have
\begin{align*}
  \TC_{\kb}(\tb)-\TC_{\kb-\vb^*}(\tb) =
\frac{1}{|\CG|} \sum_{\sigma\in \CG}  \Big( \frac{\phi_{\kb}+\phi_{\kb\sigma_{1,d+1}}}{2}
-\frac{\phi_{\kb-\vb^*}+\phi_{\kb\sigma_{1,d+1}+\vb^*}}{2} \Big)(\tb\sigma).
\end{align*}
A direct computation shows that
\begin{align*}
 & \phi_{\kb}(\tb)+\phi_{\kb\sigma_{1,d+1}}(\tb)-\phi_{\kb-\vb^*}(\tb)-
  \phi_{\kb\sigma_{1,d+1}+\vb^*} (\tb) \\
  & \quad  =   \e^{ \frac{\pi i}{d+1} \left( \sum_{\nu=2}^d 2k_{\nu}
     t_{\nu} + (k_1+k_{d+1})  (t_1+t_{d+1})\right) }
   \left( \e^{\pi i (t_1-t_{d+1})} - \e^{\pi i(t_{d+1}-t_1)}\right) \\
  & \qquad\quad
   \times \left(  \e^{\frac{\pi i}{d+1} (k_1-k_{d+1}-d-1)(t_{1}-t_{d+1})}
  -\e^{\frac{\pi i}{d+1} (k_1-k_{d+1}-d-1)(t_{d+1}-t_{1})} \right) \\
 & \quad =\ -4\sin \pi  (t_1-t_{d+1})\, \sin \pi \left(\tfrac{k_1-k_{d+1}}{d+1}-1\right)
  (t_{1}-t_{d+1})\,
 \e^{ \frac{\pi i}{d+1} \left( \sum_{\nu=2}^d 2k_{\nu}t_{\nu} + (k_1+k_{d+1})(t_1+t_{d+1})\right) }.
\end{align*}
Hence, for any $\jb \in \Lambda_n$ and $\kb \in \HH$ with $k_1-k_{d+1}=(n+1)(d+1)$,
\begin{align*}
(\phi_{\kb}+\phi_{\kb\sigma_{1,d+1}}-\phi_{\kb-\vb^*}-\phi_{\kb\sigma_{1,d+1}+\vb^*})
   \left(\tfrac{\jb}{(d+1)n}\right) =0.
\end{align*}
which yields that
\begin{align*}
(\TC_{\kb}-\TC_{\kb-\vb^*})\left(\tfrac{\jb}{(d+1)n}\right) = 0 \qquad \hbox{for
$\kb \in \HH$ with $k_1-k_{d+1}=(n+1)(d+1)$}.
\end{align*}
As in the proof of Theorem \ref{thm:zeroU}, by \eqref{eq:kappa}, this shows
that $T_\a - T_{\a^*}$ vanishes on $Y_n$.

Moreover, suppose $\kb\in \Lambda$ and set $\jb:= (\kb-\vb^*)^+$;
since $k_1 \ge \ldots \ge k_{d+1}$ and $k_i = k_j \mod d+1$ for $i \ge j$,
it follows that $j_1 = k_1 -(d+1)$ if $k_1 > k_2$ and $j_1 = k_1$ if
$k_1 = k_2$, and $j_{d+1} = k_{d+1} +(d+1)$ if $k_{d+1} < k_d$ and
$j_{d+1} = k_{d+1}$ if $j_{d+1} = k_d$. Consequently,
$(j_1- j_{d+1})-(k_1-k_{d+1})\in \left\{0,-d-1,-2d-2\right\}$, which shows
that  $|\a^+| \in \left\{|\a|, |\a|-1, |\a|-2 \right\}$, so that $T_{\a^*}$ is
a Chebyshev polynomial of degree at least $n-1$ and $T_\a - T_{\a^*}$
is orthogonal to all polynomials in $\Pi_{n-2}^d$.
\end{proof}

We note that the general theory on cubature formulas in view of ideals
and varieties also shows that there is a unique interpolation polynomial
in $\Pi_n^d$ based on the points in $Y_n$. This interpolation polynomials,
however, is exact the interpolation trigonometric polynomial in Theorem
\ref{thm:4.7}. We shall not stay the result formally.

\end{document}